\documentclass[10pt]{article}  
\usepackage[utf8]{inputenc}
\usepackage{amsmath}
\usepackage{amsfonts}
\usepackage{cancel}
\usepackage{amssymb}
\usepackage{graphicx}
\usepackage{mathrsfs}
\usepackage{upref,amsthm,amsxtra,exscale}
\usepackage{stmaryrd}
\usepackage{cite}

\usepackage{lipsum}                     
\usepackage{xargs}                      
\usepackage[pdftex,dvipsnames]{xcolor}  
\usepackage[colorinlistoftodos,prependcaption,textsize=tiny,textwidth=1.2cm]{todonotes}
\newcommandx{\unsure}[2][1=]{\todo[linecolor=red,backgroundcolor=red!25,bordercolor=red,#1]{#2}}
\newcommandx{\change}[2][1=]{\todo[linecolor=blue,backgroundcolor=blue!25,bordercolor=blue,#1]{#2}}
\newcommandx{\info}[2][1=]{\todo[linecolor=OliveGreen,backgroundcolor=OliveGreen!25,bordercolor=OliveGreen,#1]{#2}}
\newcommandx{\improvement}[2][1=]{\todo[linecolor=Plum,backgroundcolor=Plum!25,bordercolor=Plum,#1]{#2}}
\newcommandx{\thiswillnotshow}[2][1=]{\todo[disable,#1]{#2}}

\usepackage[colorlinks=true,urlcolor=blue,citecolor=red,linkcolor=blue,linktocpage,pdfpagelabels,bookmarksnumbered,bookmarksopen]{hyperref}

\usepackage{todonotes}

\usepackage[cm]{fullpage}

\usepackage{subcaption}
\usepackage{caption}
\usepackage{cleveref}
\usepackage{relsize}

\usepackage{enumitem}

\def\N{\mathbb{N}}
\def\eps{\varepsilon}
\def\R{\mathbb{R}}

\def\cH{\mathcal{H}}
\def\ast{*}

\def\r{\mathbb{R}}

\def\rn{\mathbb{R}^N}

\DeclareMathOperator{\diam}{diam}

\newtheorem{theorem}{Theorem}[section]
\newtheorem{corollary}[theorem]{Corollary}
\newtheorem{definition}[theorem]{Definition}
\newtheorem{lemma}[theorem]{Lemma}
\newtheorem{proposition}[theorem]{Proposition}

\numberwithin{equation}{section}

\theoremstyle{definition}
\newtheorem{remark}[theorem]{Remark}

\title{Asymptotic analysis and uniqueness of positive solutions\\ to slightly subcritical fractional  problems\\
in symmetric domains}
\author{Francesca De Marchis, Isabella Ianni, and Alberto Saldaña}
\date{}

\begin{document}

\maketitle

\begin{abstract}
 We prove that  the fractional Lane-Emden problem
\[
(-\Delta)^s u = u^{p}\quad \text{ in }\Omega,
\qquad u>0\quad \text{ in }\Omega,
\qquad u=0 \text{ in }\mathbb R^N\setminus\Omega,
\]
admits a unique and nondegenerate solution when  $\Omega$ is a smooth and  bounded domain which is symmetric and convex with respect to $N$ orthogonal directions and the exponent $p$ is slightly subcritical. This result is new even in the case when $\Omega$ is a ball and $N\geq 2$. Our proof relies on a sharp characterization of the asymptotic behavior of low energy solutions, which we develop in general smooth bounded open sets. 

\medskip


\medskip


\end{abstract}

\tableofcontents

\section{Introduction}

Consider the following fractional Lane-Emden problem 
\begin{align}\label{equn}
(-\Delta)^s u =u^{p}\quad \text{ in }\Omega,
\qquad u>0\quad \text{ in }\Omega,
\qquad u= 0\ \mbox{ in }\mathbb R^N\setminus\Omega,
\end{align}
where $\Omega$ is a smooth bounded open set of $\mathbb{R}^N$ ($N\geq 2$), $s\in (0,1),$  
$p\in (1, 2^*_s-1)$ with $2^*_s:=\frac{2N}{N-2s}$ (the critical fractional Sobolev exponent), and $(-\Delta)^s$ is the fractional Laplacian given by
\begin{align*}
(-\Delta)^su(x):=c_{N,s}P.V. \int_{\mathbb R^N}\frac{u(x)-u(y)}{|x-y|^{N+2s}}\,dy,\qquad c_{N,s}:=s4^{s}\pi^{-\frac{N}{2}}\frac{\Gamma(\frac{N}{2}+s)}{\Gamma(1-s)}.
\end{align*}

The existence of a (least-energy) solution to \eqref{equn} follows by standard variational methods. A natural question is whether this is the only solution for this equation. This turns out to be a very challenging problem, even in the local case $s=1$. It is known that uniqueness does not hold in general (see Section \ref{Sec:MR}) and a full characterization of the circumstances under which uniqueness holds is currently a difficult open question. The nonlocal setting poses additional difficulties, since well established tools (for instance, ODE theory) is no longer available. Some partial results have been obtained under suitable assumptions on the exponent $p$, the geometry and the topology of the domain $\Omega$, and/or on the fractional exponent $s$. We give a survey of these results at the end of this introduction.  

In this paper, we show a new uniqueness and nondegeneracy result whenever the domain satisfies some partial symmetries and the exponent $p$ is sufficiently close to $2^*_s$.  Our result is new even in the case when $\Omega$ is a ball in $\R^N$ with $N\geq 2$ (uniqueness of solutions to \eqref{equn} in the one-dimensional case has been settled in the recent papers \cite{FallWeth_second_radial_preprint,DelaTorreParini} which we comment below). 

To be more precise, we say that a set $\Omega\subset \R^N$ is a Gidas-Ni-Nirenberg (GNN) type domain if
\begin{itemize}
    \item[\emph{(H1)}] $\Omega$ is convex in the $x_i$ direction, $i=1,\ldots,N$;
     \item[\emph{(H2)}] $\Omega$ is symmetric with respect to the hyperplanes $\{x_i=0\}$ for $i=1,\ldots,N.$
\end{itemize}

Our main result is the following. 
\begin{theorem}\label{thm:uniqueness1}
Let $s\in (0,1)$, $N\geq 2$, and let $\Omega\subset \R^N$ be a smooth bounded GNN type domain. Then there exists $p_0=p_0(N,s,\Omega)\in (1,2^*_s-1)$ such that, for any $p\in(p_0,2^*_s-1)$, problem \eqref{equn} has a unique solution and it is non-degenerate.
\end{theorem}

The proof of Theorem \ref{thm:uniqueness1} follows the strategy from \cite{G00}, where the local case $s=1$ is studied.  This strategy uses several ingredients from elliptic theory such as regularity results, asymptotic analysis, classification of blowing up points, Liouville theorems, Pohozaev identities, and more.  Many of these elements have been obtained recently in the nonlocal setting, whereas some others are shown in this work.  To be more precise, let us consider exponents of the type
\begin{equation}\label{defp_eps}
p=p_\eps:=2^*_s-1-\eps,
\end{equation} where  $\eps>0$ is a small parameter. We say that a family of solutions $u_\eps$ of \eqref{equn} with \eqref{defp_eps} is of \emph{low energy} if 
\begin{align}\label{low_energy_solutions_Intro}
Q_\eps(u_\eps):=\frac{
\|u_\eps\|_{\mathcal D^{s}}^2
}{\|u_\eps\|_{L^{p_\eps+1}(\Omega)}^2}\rightarrow  S_{N,s}\qquad \text{ as }\eps\to 0,\qquad \text{ where } 
\|u_\eps\|_{\mathcal D^{s}}:=\left(\frac{c_{N,s}}{2}\int_{\mathbb R^N}\int_{\mathbb R^N}\frac{\left(u_\eps(x)-u_\eps(y)\right)^2}{|x-y|^{N+2s}}\,dx dy\right)^\frac{1}{2}
\end{align}
and $S_{N,s}$ denotes the best Sobolev constant for the embedding $\mathcal D^{s}(\mathbb R^N)\hookrightarrow L^{2^*_s}(\R^N)$. In particular, if $u_\eps$ is a least energy solution of \eqref{equn}, \eqref{defp_eps}, then $u_\eps$ is a family of low energy solutions, but in general there are families of low energy solutions which do not coincide with the least energy solutions, see Section \ref{Sec:MR}.  Using this notation, two of the main steps in the proof of Theorem \ref{thm:uniqueness1} are
\begin{itemize}
    \item[$i)$] showing that, in GNN type domains, any family of solutions of \eqref{equn},  \eqref{defp_eps} is of low energy, 
    \item[$ii)$] a sharp characterization of the qualitative (blowing up) behavior of families of solutions of low energy  as $\eps\to 0$. 
\end{itemize}

Regarding the first point, we show the following. 
\begin{theorem}\label{PositiveAreLow:thm}
Let $s\in (0,1)$, $N\geq 2$, $\Omega\subset \R^N$ be a smooth bounded domain satisfying (H1) and (H2), and let $u_\eps$ be a family of solutions of \eqref{equn}-\eqref{defp_eps}, then $u_\eps$ is of low energy.  
\end{theorem}

The proof of Theorem~\ref{PositiveAreLow:thm} crucially depends on the symmetry assumptions on the domain, which guarantee that the family of solutions has only one blow up point (via a moving plane method). Afterwards, we characterize the blow up point as isolated and simple, using the ideas from \cite[Theorem 2.1]{G00} and the results from \cite{JLX14}.

Once Theorem~\ref{PositiveAreLow:thm} has been established, the proof of Theorem~\ref{thm:uniqueness1} is carried out by contradiction using blow-up techniques. At the core of the arguments lies a precise understanding of the asymptotic behavior of low-energy solutions  in the limit as $\eps\rightarrow 0$. However, the asymptotic results for equation \eqref{eq:generale_INTRO} available in the literature (see for instance \cite{PP2014, MosconiSquassina}) are rather partial, and not sufficient for our purposes.  For this reason, our next result provides a sharp  asymptotic characterization of families of solutions of low-energy as $\eps\rightarrow 0$. This analysis is performed   in general smooth bounded open sets $\Omega$, and it is of independent interest.

To state the theorem, we introduce some notation. We denote by $\mathcal G$ the Green's function of the Dirichlet fractional Laplacian in $\Omega$; in particular, $\mathcal G(x,y)=\Phi(|x-y|) - H(x,y)$, where $\Phi$ is the fundamental solution of the fractional Laplacian in $\rn$ and $H$ is the regular part of the Green's function. Let $\mathcal R(x):=H(x,x)$ denote the fractional Robin function in $\Omega$, and define
\begin{align*}
d(x):= \operatorname{dist}(x,\partial \Omega):=\min\{ |x-y|\::\:y\in \partial \Omega\}.
\end{align*}
By \cite[Lemma 6.5]{CS98}, the function $(x,y)\mapsto \frac{\mathcal G(x,y)}{d^s(x)}$ admits a continuous extension to a positive function on $\partial\Omega\times\Omega$. Furthermore, 
$\|\cdot\|_{\infty}$ denotes the essential supremum norm,  $\delta_{x_0}$ is the Dirac delta centered at a point $x_0$ and $\Gamma$ is the Gamma function. Finally, we denote by $U$  the unique positive radial solution of $(-\Delta)^s U=U^{\frac{N+2s}{N-2s}}$ in $\mathbb R^N$
 such that $U(0)=1$, namely,
 \begin{equation}\label{UIntro}
U(x):=\left(1+\tau_{N,s}|x|^2\right)^{-\frac{N-2s}{2}}\quad\mbox{ with } \tau_{N,s}:=\frac{1}{4}\left[\frac{\Gamma\left(\frac{N-2s}{2}\right)}{\Gamma\left(\frac{N+2s}{2}\right)}\right]^{\frac{1}{
s}},
\end{equation}
see \cite{CLO}.
\begin{theorem}\label{fractional:Han}
Let $s\in (0,1)$ and let $\Omega$ be a smooth bounded open set of $\R^N$, $N\geq 2$. Assume that $u_\eps$ is a family of solutions to \eqref{equn}-\eqref{defp_eps} satisfying \eqref{low_energy_solutions_Intro}. Then 
\begin{equation}\label{blow-upintro}\|u_\eps\|_{\infty}\rightarrow +\infty \qquad \text{ as }\eps\rightarrow 0.\end{equation}
Let $x_\eps\in\Omega$ be a point such that $\|u_\eps\|_{\infty}=u_\eps(x_\eps)$. Then there exists a sequence $(\eps_n)$ such that $\eps_n\to0$ as $n\to+\infty$ and $x_{\eps_n}\rightarrow x_0$, where $x_0$ is an interior point of $\Omega$ and satisfies
\[\nabla\mathcal R (x_0)=0.\]
Furthermore, let
$\gamma_{N,s}:=2^N \pi ^{N/2} \Gamma (s)\frac{\Gamma
   \left(\frac{N}{2}+s\right)^{\frac{N}{2 s}-1}}{ \Gamma \left(\frac{N}{2}-s\right)^{\frac{N}{2 s}}}$ and $S_{N,s}$ be the best fractional Sobolev constant (see \eqref{constant:SnsVALUE}),
then 
\begin{align}
    \lim_{n \rightarrow +\infty} \eps_n \|u_{\eps_n}\|_{\infty}^2
    &=\gamma_{N,s}^2 \frac{2N }{(N-2s)\left(S_{N,s}\right)^{\frac{N}{2s}}}\mathcal R(x_0); \label{eps_normueps2}\\
    \lim_{n \rightarrow +\infty}\|u_{\eps_n}\|_{\infty} {u_{\eps_n}}&=  
 \gamma_{N,s}{\mathcal G(\cdot,x_0)}
 \qquad\mbox{ in  }\ C_{loc}(\overline\Omega\setminus \{x_0\});\label{gamma0_Ga}\\
 \lim_{n \rightarrow +\infty}\|u_{\eps_n}\|_{\infty} \frac{u_{\eps_n}}{d^s}
 &=  
 \gamma_{N,s}\frac{\mathcal G(\cdot,x_0)}{d^s}
 \qquad\mbox{ in  }\ C_{loc}(\overline\Omega\setminus \{x_0\}).\label{gamma0_G_delta}
\end{align}
Finally, there exists $\Lambda >0$ such that
 \begin{equation}\label{BuonduIntro}
   u_\eps(x)\leq \Lambda \|u_\eps\|_{\infty}U\left(\frac{x-x_\eps}{\|u_\eps\|_{\infty}^{-\frac{p_\eps-1}{2s}} }\right)\qquad \text{ for all }x\in\rn,\end{equation}
and, setting 
\begin{equation}\label{rescaling_Intro}
v_\eps(x):=\|u_\eps\|_{\infty}^{-1}u_\eps(
\|u_\eps\|_{\infty}^{-
\frac{p_\eps-1}{2s}
}x+x_\eps),
\end{equation}
we have
\begin{align}\label{convvIntro}
    v_\eps\rightarrow U\qquad \text{  in  $C^{2s-\delta}_{loc}(\mathbb R^N)$ for any $\delta>0$ small as } \eps\rightarrow 0.
    \end{align}
\end{theorem}
It is not difficult to see that \eqref{blow-upintro} must hold for low energy solutions, since these are defined by condition \eqref{low_energy_solutions_Intro} and it is well known that the optimal  fractional Sobolev constant $S_{N,s}$ is not attained in domains  $\Omega\neq\mathbb R^N$. Theorem~\ref{fractional:Han} provides a finer description of the blow-up behavior. In particular, it shows that low-energy solutions blow up at exactly one point, which is a critical point of the fractional Robin function. This analysis is new even for least energy solutions, moreover, it is sharp, as all constants involved are explicit. These accurate results are made possible by recent Pohozaev-type identities for fractional Green's and Robin functions, established in  \cite{DjitteSueur,DiebIanni_poho} (see Lemmas~\ref{lemma:PohozaevGreen_DiebIanni}-\ref{lemma:PohozaevGreen_DjitteSeur}).
Finally,  the limit profile $U$ appearing in \eqref{convvIntro}  is given explicitly by a fractional Talenti bubble, namely one of the extremal functions for the fractional Sobolev inequality (which form an $(N+1)$-parameter family due to scaling and translation invariance, see Section~\ref{section:SobolevConstant_and_bubbles}). 
When $s=1$, Theorem~\ref{fractional:Han} 
directly recovers the classical results in \cite{H91,R89,AP87,BrezisPeletier}, see Remark~\ref{ball:rmk}.  We also refer to \cite{chokimlee}, where an analogue of Theorem~\ref{fractional:Han} is shown for the spectral fractional Laplacian.  Related asymptotic results for slightly subcritical Hartree equations involving the spectral fractional Laplacian were obtained in \cite{BCYZ}; in particular, for $2s<N<6s$ and families
satisfying a low-energy condition, they establish single-point
concentration at a critical point of the associated Robin function
and determine the precise blow-up rate and the limiting Green
function profile away from the concentration point.

When $\Omega$ satisfies  assumptions \emph{(H1)} and \emph{(H2)}, Theorem~\ref{fractional:Han} simplifies.  Indeed, by moving plane arguments (see \cite[Corollary 1.2]{JW14}), any solution  to \eqref{equn} is reflectionally symmetric with respect to the hyperplanes $\{x_i=0\}$ and strictly decreasing in $|x_i|$. 
    In particular, it follows that $x_\eps=x_0=0$, and consequently  the rescaled solutions $v_\eps$ in \eqref{rescaling_Intro} inherit the same symmetries as the original solutions $u_\eps$. In this symmetric setting, we can exploit all the information  provided by Theorem~\ref{fractional:Han} to establish our main uniqueness and nondegeneracy result.

Note that Theorem~\ref{fractional:Han} also shows that Theorem~\ref{thm:uniqueness1} fails in general if assumptions \emph{(H1)} and \emph{(H2)} are removed. Indeed, let $\Omega$ be the union of two identical disjoint balls, and consider a family of least energy solutions of \eqref{equn} in $\Omega$. By Theorem~\ref{fractional:Han}, for all $p$ sufficiently close to $2^*_s-1$, one finds that, up to a subsequence, these  solutions concentrate at one (and only one) of the two equal balls; in particular, they are  not invariant under reflections. Exploiting the symmetry of $\Omega$, one can then obtain a distinct family of least energy solutions by reflection. Hence, uniqueness of least-energy solutions fails in general. Moreover, using a Lyapunov-Schmidt reduction (see \cite[Theorem 1.2]{DRS17}) one can also prove the existence of multiple low-energy solutions; for instance
when $\Omega$ is the union of two disjoint balls with different radii (see \cite[Remark 4.9]{DIS23}). This last example also exhibits a family of low-energy solutions that are not of least-energy type.

When $\Omega$ is a general convex domain, without any symmetry assumption, we expect the uniqueness and nondegeneracy results to hold for any fractional power $s\in (0,1)$. However, the asymptotic analysis developed  in  Theorem~\ref{fractional:Han} is not  sufficient to establish such results in the general case, and a deeper analysis would be required. This remains an extremely challenging open problem for  further investigations.

We also mention that, if we let $\Omega=B_1(0)$ be the unit ball centered at zero, then the estimates in Theorem~\ref{fractional:Han} are significantly more explicit.

\begin{corollary}\label{ball:thm}
Let $\Omega=B_1\subset \rn$, $N\geq 2$ and $s\in (0,1)$. Let $u_\eps$ be a family of low energy solutions of \eqref{equn}-\eqref{defp_eps}. Then,
\begin{align}
\lim_{\eps \to 0} \eps \|u_\eps\|_{\infty}^2
 &=\frac{2^{N-2 s+2}N! }{(N-2 s)^2} \frac{\left(\Gamma
   \left(\frac{N}{2}+s\right)\right)^{\frac{N}{2 s}-2}}{\left(\Gamma \left(\frac{N}{2}-s\right)\right)^{\frac{N}{2 s}}}, \label{1ball}\\
 \lim_{\eps\to 0}\eps^{-\frac{1}{2}}\frac{u_\eps}{d^s}(z)
 &= \frac{2^{\frac{N}{2}-1} (N-2 s) \Gamma \left(\frac{N}{2}\right) 
  }{\Gamma(s+1) \sqrt{N!}}\left(\frac{\Gamma \left(\frac{N}{2}+s\right)}{\Gamma\left(\frac{N}{2}-s\right)}\right)^{\frac{N}{4 s}}\quad \text{ for all }z\in \partial B_1.
  \label{2ball}
\end{align}
\end{corollary}
 Similar explicit constants have been important tools in the asymptotic analysis of elliptic problems (see for instance, \cite{AP87,GST20,IS21,H91,R89}) and are therefore of independent interest.

\subsection*{Previously known results}

Next, we survey some of the previously known results regarding uniqueness and nondegeneracy results for fractional Lane-Emden type equations. 

 Establishing conditions which guarantee uniqueness and nondegeneracy results is a challenging and  largely unsolved open question for solutions of problem \eqref{equn}, or of the more general problem
\begin{align}\label{eq:generale_INTRO}
(-\Delta)^s u +\lambda u=u^{p}\quad \text{ in }\Omega,
\qquad u>0\quad \text{ in }\Omega,
\qquad u= 0\ \mbox{ in }\mathbb R^N\setminus\Omega,
\end{align}
with $\lambda>-\lambda_1^s(\Omega)$ to ensure existence. Here $\lambda_1^s(\Omega)$ denotes the first Dirichlet eigenvalue of $(-\Delta)^s$ in $\Omega$. 
Recently,  uniqueness and nondegeneracy have been obtained in \cite{DIS23} for general solutions of \eqref{eq:generale_INTRO} in smooth bounded open sets $\Omega$\footnote{In \cite{DIS23} $\Omega$ is assumed to be a domain, but the connectedness is not used in the proof}, under the assumption that  the power $p$ is close enough to $1$ (the almost linear regime). Let us observe that this is the best that one can expect under  no other additional assumptions on the set $\Omega$, indeed  one may find sets $\Omega$ for which  multiplicity results hold true for larger values of $p$, see the discussion before Corollary~\ref{ball:thm}. 

However, in \cite{DIS23} it has been proved that if $\Omega$, $p,$ and $\lambda$ are  such that uniqueness and nondegeneracy are known to hold true for the solutions to the corresponding local problem
\begin{align}\label{eq:generaleLOCAL_INTRO}
-\Delta u +\lambda u=u^{p}\quad \text{ in }\Omega,
\qquad u>0\quad \text{ in }\Omega,
\qquad u= 0\ \mbox{ in }\partial\Omega,
\end{align}
where the fractional Laplacian $(-\Delta)^s$ is replaced by the standard Laplacian $-\Delta$ (which formally corresponds to the case
$s=1$), then the uniqueness and nondegeneracy properties also hold true  for the solutions of \eqref{eq:generale_INTRO}, with the same choice of $\Omega$, $p$, and $\lambda$, \emph{provided the fractional power $s$ is close enough to $1$}, i.e., in the almost local regime.

Sufficient conditions for uniqueness and nondegeneracy  for the standard Laplacian case  \eqref{eq:generaleLOCAL_INTRO}  have been widely studied (see \cite{GNN1979, NiJDE1983, NiNussbaumCPAM1985,  KwongLiTAMS1992, ZhangCPDE1992, SrikanthDiffIntEq1993, AdimurthiYadavaARMA1994, AftalionPacellaJDE2003, ZouPisa1994, DancerJDE1988, DamascelliGrossiPacellaAIHP1999,  G00, DancerRockyMountain1995, McKennaPacellaPlumRothJDE2009, McKennaPacellaPlumRoth2012, lin1994, DeMarchisGrossiIanniPacellaJMPA2019, GrossiIanniLuoYanJMPA2022, LiWeiZou, DancerMa2003 }); hence, their counterparts for the nonlocal problem \eqref{eq:generale_INTRO} are now established  as a consequence of this general result in \cite{DIS23}, but only for $s$ close enough to $1$. These nonlocal results include cases when $\Omega$ is a ball, or a general convex set, or a general symmetric set of  Gidas-Ni-Nirenberg (GNN) type, namely when it satisfies \emph{(H1)} and \emph{(H2)}. In particular,  if $\Omega$ is a ball, one gets uniqueness and nondegeneracy for the solutions of \eqref{eq:generale_INTRO}  
for  any $\lambda,p$ admissible for the existence (see   \cite{GNN1979, NiJDE1983, NiNussbaumCPAM1985,  KwongLiTAMS1992, ZhangCPDE1992, SrikanthDiffIntEq1993, AdimurthiYadavaARMA1994, AftalionPacellaJDE2003} for the local case). In the case when $\Omega$ is a general convex set, they are proved when $\lambda=0$ and $p=p_\eps$ with $\eps$ small (for the local case, this has been first established in dimension $N=2$ {and $p$ large} in \cite{DeMarchisGrossiIanniPacellaJMPA2019, GrossiIanniLuoYanJMPA2022}, and later extended to $N\geq 3$ and $p=p_\eps$ in \cite{LiWeiZou}  under the additional assumption that the Robin function of $\Omega$ is Morse;  see also \cite{lin1994} for the case of least energy solutions in dimension $N=2$ and any admissible $p$). When $\Omega$ is a GNN type domain, uniqueness and nondegeneracy for the nonlocal problem \eqref{eq:generale_INTRO}
are deduced in the case   $\lambda=0$ and $p=p_\eps$  (see \cite{DancerJDE1988, DamascelliGrossiPacellaAIHP1999, G00} for the local case), or in the case  $\lambda>0$ and any admissible $p$,  under the additional assumption that  $\Omega$ is sufficiently  \emph{large}  (see \cite{DancerRockyMountain1995} for the local case).

We emphasize again that all these listed results obtained for problem \eqref{eq:generale_INTRO} are true for $s$ close to $1$. To extend them to \emph{any} fractional power $s\in (0,1)$ is a challenging and interesting problem. 

When $\Omega$ is a ball, $\lambda\geq 0$, and $p$ is admissible for the existence, this has been recently accomplished by Fall and Weth \cite{FallWeth_second_radial_preprint} for families of  solutions of \eqref{eq:generale_INTRO} having Morse index equal to $1$ (see also \cite{FallWeth_nondegeneracy_preprint}). We stress that even though in this case solutions of \eqref{eq:generale_INTRO} are radial, ODE techniques are not available due to the nonlocality of the operator. The core of the proof is to show  the nondegeneracy property, since the uniqueness then follows from it by a continuity argument in $p$ combined with the results in \cite{DIS23} for the almost linear problem $p\sim 1$. Their proof of the nondegeneracy consists in showing first a radial symmetry result for the solutions of the linearized problem and, eventually, both a Courant nodal type theorem  and  a new Hopf lemma for its radial solutions (when $\lambda=0$, the same result has been also obtained in \cite{DelaTorreParini} with a similar strategy). This approach, although working for any fractional power $s\in(0,1)$, on many points strongly relies on the  assumption that the set $\Omega$ is a ball, and it cannot be extended to general convex sets, or to general GNN symmetric domains. 

When $\Omega$ satisfies assumptions \emph{(H1)} and \emph{(H2)}, and for any $s\in (0,1)$,  uniqueness and nondegeneracy 
have been obtained recently  by Dieb, Ianni, and Saldaña \cite{DiebIanniSaldanaCalcVar} for families of least energy solutions to \eqref{eq:generale_INTRO}, for $\lambda>0$, for  any admissible exponent $p$, under the additional assumption that $\Omega$ is sufficiently large. 
Their proof combines a symmetry result for the linearized problem\textemdash established in general GNN type domains\textemdash with an asymptotic argument that ultimately relies on the celebrated uniqueness and nondegeneracy result for the fractional Schr\"odinger equation in the whole space $\mathbb R^N$, proved by  Frank-Lenzmann-Silvestre \cite{FrankLenzmannSilvestre}. While the symmetry result in \cite{DiebIanniSaldanaCalcVar} is also proved  in the case $\lambda=0$, the asymptotic argument no longer applies.

We also mention that in dimension $N=1$, where a GNN type domain is simply an interval, Fall and Weth \cite{FallWeth_nondegeneracy_preprint} proved uniqueness and nondegeneracy for \emph{general} solutions of \eqref{equn}.

Theorem~\ref{thm:uniqueness1} extends the uniqueness and nondegeneracy result known for the fractional Lane-Emden problem in a ball \cite{DIS23, FallWeth_second_radial_preprint, DelaTorreParini} to the more  general setting of GNN type domains and to general nonnegative solutions. However, unlike \cite{DIS23, FallWeth_second_radial_preprint, DelaTorreParini}, we are not able to treat all exponents 
$p$ for which existence holds; instead, our result applies when $p$ 
is sufficiently close to the critical exponent $2^*_s-1$. 
This restriction is intrinsically related to our proof, which relies heavily on the asymptotic description of families of low energy solutions in the limit  
as $\eps\rightarrow 0$.
We recall that uniqueness and nondegeneracy for problem \eqref{equn} have also been established for 
$p$ sufficiently close to $1$
 (see \cite{DIS23}), and we conjecture that these properties hold for all 
$p\in (1,2^*_s-1)$ in GNN type domains.

\medskip

The paper is organized as follows. Section~\ref{Section:preliminariesandnotations} collects the notation, functional framework, and preliminary results. In Section~\ref{section:asymptotic}, we prove Theorem~\ref{fractional:Han}. In  Sections~\ref{section:uniqueness} and~\ref{section:nondeg}, we focus on domains $\Omega$ satisfying assumptions \emph{(H1)} and \emph{(H2)} and establish  the uniqueness and nondegeneracy results, respectively. Finally, in Section~\ref{ball:sec}, we discuss the case of the unit ball.

\section{Notations and preliminaries}\label{Section:preliminariesandnotations}

Let  $v:\mathbb R^N\rightarrow\mathbb R$. For an open set $A\subseteq\mathbb R^N$,  we use the following  notations
\begin{align*}
&\|v\|_{L^\alpha(A)}:=\left(\int_A |v|^\alpha\right)^\frac{1}{\alpha}, \ \mbox{for }\alpha\in[1,\infty),\\
& \|v\|_{L^\infty(A)}:=
esssup_A |v|.
\end{align*}
When $A=\rn$, we simply write 
$\|v\|_\alpha:= \|v\|_{L^\alpha(\rn)}$,   $\|v\|_{\infty}:=\|v\|_{L^\infty(\rn)}$.
Moreover, we define
\begin{align*}
&  [v]_{C^\sigma(A)}:=\sup_{\substack{x,y\in A\\x\neq y}}\frac{|v(x)-v(y)|}{|x-y|^\sigma}\qquad \mbox{ for }\sigma\in (0,1)\\
& \|v\|_{C^m}(A):=
 \sum_{|\alpha|\leq m-1}\|\partial^\alpha v\|_{L^\infty(A)}+\sup_{\substack{x,y\in A\\x\neq y\\ |\alpha|=m-1}}\frac{|\partial^\alpha v(x)-\partial^\alpha v(y)|}{|x-y|} \qquad \mbox{ for } m\in \N\\
&\|v\|_{C^{m+\sigma}(A)}:= \sum_{|\alpha|\leq m}\|\partial^\alpha v\|_{L^\infty(A)}+\sup_{|\alpha|=m}[\partial^\alpha v]_{C^\sigma(A)}\qquad \mbox{ for } m\in \N,\quad \sigma\in (0,1).
\end{align*}
Here the space  $C^{m+\sigma}(A)$  denotes the space $C^{m,\sigma}(A)$ of $m$-times continuously differentiable functions in $A$  whose derivatives of order $m$ are $\sigma$-H\"older continuous in $A$.    As usual, we use $C_{loc}^\sigma(A)$, $C^m_{loc}(A)$, $C_{loc}^{m+\sigma}(A)$ for the local version of these spaces.

Moreover, similar notations are used also for functions $v:\mathbb R^{N+1}\rightarrow \mathbb R$, and sets $A\subseteq\mathbb R^{N+1}$.  Furthermore, let \begin{align*}
& \|v\|_{\mathcal D^{s}}:=[v,v]_s^\frac{1}{2},
\end{align*}
where  
\begin{align}\label{cns}
[w,v]_s:=\dfrac{c_{N,s}}{2}\iint_{\R^N\times \R^N}\dfrac{(w(x)-w(y))(v(x)-v(y))}{|x-y|^{N+2s}}dx\,dy
,\qquad c_{N,s}:=s4^{s}\pi^{-\frac{N}{2}}\frac{\Gamma(\frac{N}{2}+s)}{\Gamma(1-s)}.
\end{align} 
We also let
\begin{align*}
& H^s(\R^N):=\left\lbrace v \in
L^2(\R^N):\|v\|_{\mathcal D^{s}}<+\infty \right\rbrace,\quad 
\mathcal D^{s}(\mathbb R^N):=\left\{ v\in L^{2^*_s}(\rn): \|v\|_{\mathcal D^{s}}<\infty \right\},
\end{align*}
and
\begin{align*}
\cH^s_0(A):=\left\{ v \in H^{s}(\R^N):\, \, v=0 \, \text{ in } \R^N\backslash A \right\}.  
\end{align*}
We define the Hilbert space $\mathcal D_0^{s}(A)$ as the completion of $C^{\infty}_c(A)$ with respect to the norm $\|\cdot\|_{\mathcal D^s}$.  Note that if $A$ is a bounded Lipschitz set, then $\mathcal D_0^{s}(A)=\cH^s_0(A)$, furthermore $\cH^s_0(A)\subset \mathcal D^{s}(\mathbb R^N)=\mathcal D_0^{s}(\mathbb R^N)$.

\subsection{Regularity via harmonic extension}\label{subsec:ext}

In this subsection we use the well-known local extension of fractional problems (sometimes called the Caffarelli-Silvestre extension, $s$-harmonic extension, or Poisson extension) to show some crucial regularity results.  Let us begin by introducing some notation. Let
\[
\mathbb{R}_{+}^{N+1}:=\left\{(x, t)\ :\  x \in \mathbb{R}^{N}, t>0\right\}.
\] 
For $r>0$ and $x_0\in \R^{N}$, let
\begin{align}
B_r(x_0)&:=\{x\in \R^N\::\: |x-x_0|^2<r^2\},\\
B_r^+(x_0)&:=\{(x,t)\in \R^{N+1}_+\::\: |x-x_0|^2+t^2<r^2\},\\
Q_r(x_0)&:=B_r(x_0)\times (0,r)\subset \R^{N+1}_+,\label{Qdef}
\end{align}
and we simplify the notation as $B_r=B_r(0)$, $B_r^+:=B_r^+(0)$, $Q_r:=Q_r(0)$, when $x_0=0$.

For a set $D\subset \overline{\mathbb R^{N+1}_+}$ we write $\partial D=\partial 'D\cup\partial ''D$, where
\begin{align}\label{boundaries}
\partial 'D:=\operatorname{int}\left(\overline D\cap (\mathbb R^N\times \{0\})\right),\qquad 
\partial ''D:=\partial D\setminus \partial 'D.
\end{align}

Let  $V:\mathbb R^{N+1}\rightarrow\mathbb R^M$, $M\in\mathbb N$. For $D\subseteq\mathbb R^{N+1}$ open or the closure of an open set,  we 
denote by  $L^2(t^{1-2s},D)$  the Banach space of all measurable functions $V$ defined on $D$ such that
\begin{align*}
   & \|V\|_{L^2(t^{1-2s},D)}:=\left( \int_D t^{1-2s}|V(x,t)|^2\,dxdt\right)^{\frac12} < +\infty.
\end{align*}
We say that $V\in H(t^{1-2s},D)$ if $V\in L^2(t^{1-2s},D)$ and its weak derivatives $\nabla V$ exist and belong to $L^2(t^{1-2s},D)$. The norm of $V$ in $H(t^{1-2s},D)$ is given by 
\[
\|V\|_{H(t^{1-2s},D)}:=\left( \|V\|_{L^2(t^{1-2s},D)}^2+\| \nabla V\|_{L^2(t^{1-2s},D)}^2\right)^{\frac12}.
\]

Let $D^{1,2}\left(t^{1-2 s},\mathbb{R}_{+}^{N+1} \right)$ be the completion of $C^\infty_{c}\left(\overline{\mathbb{R}_{+}^{N+1}}\right)$ with respect to the norm $\|\nabla V\|_{L^2(t^{1-2s},\mathbb{R}_{+}^{N+1})}.$ For an open set $A\subseteq \mathbb R^N$ we also consider the space
\begin{align}\label{DAspace}
D_A^{1,2}\left(t^{1-2 s}, \mathbb{R}_{+}^{N+1} \right):=\left\{V \in D^{1,2}\left(t^{1-2 s}, \mathbb{R}^{N+1}_{+} \right): V(\cdot, 0)=0 \text { on } \mathbb{R}^N \backslash A\right\} .    
\end{align}

We write $V \in L^\infty_{\mathrm{loc}}(\overline{\mathbb{R}^{N+1}_+})$ if $V \in L^\infty(\overline{Q}_R)$ for any $R > 0$. Similarly, we write 
\begin{align}\label{Hlocdef}
V \in H_{\mathrm{loc}}(t^{1-2s}, \overline{\mathbb{R}^{N+1}_+})
\end{align}
if $V \in H(t^{1-2s}, \overline{Q_R})$ for any $R > 0$.   If $v \in \mathcal{D}^{s}_0(A)$ and 
\begin{align}\label{P}
V(x, t):=p_{N,s} t^{2 s}\int_{\mathbb{R}^N} \frac{v(y) }{\left(t^2+|x-y|^2\right)^{\frac{N+2 s}{2}}}\, dy\qquad \text{ for }x\in \R^N,\ t>0,
\end{align}
with
\begin{align}\label{def:pNs}
p_{N,s}:=\left(\int_{\mathbb{R}^N}\left(1+|z|^2\right)^{-\frac{N+2 s}{2}} d z\right)^{-1}
=\frac{\Gamma\left(\frac{N}{2}+s\right)}{\pi^{\frac{N}{2}}\Gamma(s)},    
\end{align}
then $V \in D_A^{1,2}\left(t^{1-2 s}, \mathbb{R}_{+}^{N+1} \right)$ and it is called the \emph{local extension} of $v$ (see, for instance, \cite{CS14,CS07}). This function solves
\begin{align}\label{harm}
\operatorname{div}\left(t^{1-2 s} \nabla V\right)=0\quad \text{ in }\R^{N+1}_+.    
\end{align}
It is known (see, e.g. \cite{LionsMagenes}) that there exists $C>0$ depending only on $N$ and $s$ such that for $V\in H(t^{1-2s}, \R^{N+1}_+)\cap C(\overline{\R^{N+1}_+})$, $\|V(\cdot,0)\|_{\mathcal D^s}\leq C \|V\|_{H(t^{1-2s},\R^{N+1}_+)}$. Hence, by a standard density argument, every $V\in H(t^{1-2s},\R^{N+1}_+)$ has a well-defined trace $v:=V(\cdot,0)\in \mathcal D^s(\R^N)$.  Moreover (see \cite{CS14}),
\begin{align*}
\int_{\mathbb{R}_{+}^{N+1}} t^{1-2 s} \nabla V \cdot \nabla \Phi\ dx\, dt=k_s [v,\varphi]_s
\end{align*}
for all $\Phi \in D_A^{1,2}\left(t^{1-2 s}, \mathbb{R}_{+}^{N+1} \right)$ with $\varphi=\left.\Phi\right|_{\mathbb{R}^N \times\{0\}} \in \mathcal{D}_0^{s}(A),$ where 
\begin{align}\label{dsdef}
k_s:= \frac{\Gamma(1-s)}{2^{2 s-1}\Gamma(s)},    
\end{align}
see \cite{CS14}. If $v \in \mathcal{D}_0^{s}(A) \cap C^\infty(A) \cap C^s(\mathbb R^N)$, then 
$$
V \in C_{l o c}^s\left(\overline{\mathbb{R}_{+}^{N+1}}\right) \cap C^\infty\left(\mathbb{R}_{+}^{N+1}\right)
$$
and 
$$
-\lim _{t \rightarrow 0} t^{1-2 s} \partial_t V(x, t)=k_s(-\Delta)^s v(x) \quad \text { for all } x \in A.
$$

\subsubsection{A reverse Hölder inequality}

We begin with an auxiliary lemma.
\begin{lemma}\label{est:lem}
Let $v\in L^{\beta}(\R^N)$, $\beta>1,$ and let $V$ denote its local extension. Then,
\begin{align*}
    \int_{Q_r} t^{1-2s}|V(x,t)|^{\beta} \,  dx\,dt 
    \leq \frac{r^{2(1-s)}}{2(1-s)}\int_{\R^N} |v(x)|^{\beta}\, dx\qquad  \text{ for all }r>0. 
    \end{align*}
\end{lemma}
\begin{proof}
Using \eqref{P}, a change of variables, and Jensen's inequality, 
\begin{align*}
|V(x,t)|^{\beta}
= \left(\int_{\mathbb{R}^N} v(tz+x) \, \frac{ p_{N,s}} {\left(1+|z|^2\right)^{\frac{N+2 s}{2}}}dz\right)^{\beta}\leq \int_{\mathbb{R}^N} |v(tz+x)|^{\beta} \, \frac{ p_{N,s} }{\left(1+|z|^2\right)^{\frac{N+2 s}{2}}}dz,
\end{align*}
then, using Fubini, and the definition of $p_{N,s}$ in \eqref{def:pNs}
\begin{align*}
\int_{Q_r}  t^{1-2s}|V(x,t)|^{\beta} \,  dx\,dt
&\leq \int_{\mathbb{R}^N}\frac{ p_{N,s}}{\left(1+|z|^2\right)^{\frac{N+2 s}{2}}}\int_0^r  t^{1-2s}\int_{B_r(0)} |v(tz+x)|^{\beta} \, dx\, dt\,dz\\
&\leq \int_0^r  t^{1-2s}\, dt\int_{\R^N} |v(x)|^{\beta}\, dx
=\frac{r^{2(1-s)}}{2(1-s)}\int_{\R^N} |v(x)|^{\beta}\, dx 
.
\end{align*}
\end{proof}
    
\begin{lemma}\label{lemma:Lemma_6_Han}
Let $\Omega$ be a domain in $\R^N,$ $2<\beta\leq 2^*_s+1$, and let $v\in \mathcal D_0^{s}(\Omega)\cap L^\beta(\R^N)
$ be a nonnegative smooth solution of 
\begin{align*}
    (-\Delta)^s v = b(x) v^{q-1}\quad \text{ in }\Omega,
\end{align*}
where $b\in L^\infty(\Omega)$ and $2<q_0\leq q\leq 2^*_s=\frac{2N}{N-2s}.$  Then there is $\Theta_0>0$ and $r_0>0$ depending on $N,s,\|b\|_{\infty},q_0,$ and $\beta$ such that, for any $x_0\in \R^N$ and any $0<r\leq r_0$ with 
\begin{align*}
\int_{\Omega \cap B_{2r}(x_0)} v^q \leq \Theta_0,
\end{align*}
we have that there exists $C=C(N,s,\beta)>0$ such that
\begin{align*}
\|v\|_{L^{\frac{\beta \,2^*_s}{2}}(\Omega \cap B_{r}(x_0))}
\leq Cr^{-\frac{2s}{\beta}} \|v\|_{L^{\beta}(\mathbb R^N)}.
\end{align*}
\end{lemma}
\begin{proof}
Let $\eta:\R^{N+1}\to \R$ be a smooth cutoff function such that
\begin{align*}
\eta(x,t)=1\quad \text{ if } \sqrt{|x-x_0|^2+t^2}<r,\quad
\eta(x,t)=0\quad \text{ if } \sqrt{|x-x_0|^2+t^2}\geq 2r,
\end{align*}
$0\leq \eta\leq 1$ and $|\nabla \eta|\leq \frac{2}{r}.$  Let $V$ denote the extension of $v$ and we assume without loss of generality that $x_0=0$. Being $\overline{B_{2r}}=\partial B_{2r}^+\cap \partial \R^{N+1}_+$, testing \eqref{harm} with $\eta^2 V^{\beta-1}$ and integrating by parts, we get
\begin{align*}
0
&=\int_{B^+_{2r}}t^{1-2s}\nabla V \nabla(\eta^2 V^{\beta-1})\,  dxdt
-\int_{ B_{2r}}\eta(x,0)^2 v(x)^{\beta+q-2} b(x)\, dx.
\end{align*}
Now we adapt the proof of \cite[Lemma 6]{H91}. Using
\begin{align*}
2\nabla V \nabla \eta V^{\beta-1} \eta
\geq -\frac{1}{2}\left(|(\beta-1)^\frac{1}{2}\nabla V V^{\frac{\beta-2}{2}}\eta|^2+|2(\beta-1)^{-\frac{1}{2}}\nabla \eta V^\frac{\beta}{2}|^2\right),
\end{align*}
we have
\begin{align*}
&\int_{B^+_{2r}}t^{1-2s}\nabla V \nabla(\eta^2 V^{\beta-1})\,  dxdt
=\int_{B^+_{2r}} t^{1-2s}\left( 2\nabla V \nabla \eta V^{\beta-1} \eta
+ (\beta-1) |\nabla V|^2 \eta^2 V^{\beta-2}\right)\,  dxdt\\
&\geq \int_{B^+_{2r}} t^{1-2s}\left(
\frac{\beta-1}{2} |\nabla V|^2 \eta^2 V^{\beta-2}
-\frac{2}{\beta-1} |\nabla \eta|^2  V^{\beta}\right)
\,  dxdt.
\end{align*}
In particular, 
\begin{align}\label{intermed}
\frac{\beta-1}{2} \int_{B^+_{2r}} t^{1-2s} \eta^2
|\nabla V|^2  V^{\beta-2}
\,  dxdt
\leq 
\frac{2}{\beta-1} \int_{B^+_{2r}} t^{1-2s}|\nabla \eta|^2  V^{\beta}
\,  dxdt+\int_{ B_{2r}}\eta(x,0)^2 v(x)^{\beta+q-2} b(x)\, dx.
\end{align}
By the trace Sobolev embedding (see, for instance, \cite[Corollary 5.3]{GQ13}), there is a constant $\kappa_{N,s}>0$ such that 
\begin{align*}
&\kappa_{N,s} \left(
\int_{ B_{2r}}
(\eta v^{\frac{\beta}{2}})^{2^*_s}\, dx
\right)^\frac{2}{2^*_s}
\leq \int_{B^+_{2r}} t^{1-2s}|\nabla (\eta V^{\frac{\beta}{2}})|^2\,  dxdt\\
&\leq 2\int_{B^+_2r} t^{1-2s}\left(
|\nabla \eta|^2 V^{\beta}
+
\frac{\beta^2}{4}
\eta^2|\nabla V|^2 V^{\beta-2}
\right)\,  dxdt\\
&=\frac{\beta^2}{\beta-1}
\left[\frac{\beta-1}{2}\int_{B^+_2r} t^{1-2s}
\eta^2|\nabla V|^2 V^{\beta-2}
\,  dxdt\right]
+
2\int_{B^+_2r
} t^{1-2s}|\nabla \eta|^2 V^{\beta}\,  dxdt\\
&\overset{\eqref{intermed}}{\leq} \left(\frac{2\beta^2}{(\beta-1)^2}+2\right)
\int_{B^+_{2r}} t^{1-2s}|\nabla \eta|^2  V^{\beta}
\,  dxdt
+\frac{\beta^2}{\beta-1}\int_{ B_{2r}}\eta^2 v^{\beta+q-2} b(x)\, dx\\
&\leq 10
\int_{B^+_{2r}} t^{1-2s}|\nabla \eta|^2  V^{\beta}
\,  dxdt
+4\|b\|_{\infty}(\beta-1)\int_{ B_{2r}}\eta^2 v^{\beta+q-2}\, dx\\
&\leq 10
\int_{B^+_{2r}} t^{1-2s}|\nabla \eta|^2  V^{\beta}
\,  dxdt
+4\|b\|_{\infty}(\beta-1)
\left(\int_{ B_{2r}}(\eta v^{\frac{\beta}{2}})^{2^*_s}\, dx\right)^{\frac{2}{ 2^*_s  }}
\left(\int_{ B_{2r}}v^{(q-2)\frac{N}{2s}} \, dx\right)^{\frac{2s}{N}}.
\end{align*}
Since $q\leq \frac{2N}{N-2s}$, we have that $2sq\geq N(q-2)$.  Choosing $\Theta_0$ and $r_0$ small enough so that 
\begin{align*}
    4\|b\|_{\infty}(\beta-1)\left(\int_{ B_{2r}}v^{(q-2)\frac{N}{2s}} \, dx\right)^{\frac{2s}{N}}
   & \leq 
4\|b\|_{\infty}(\beta-1)
\left(\int_{ B_{2r}}v^{q} \, dx\right)^{\frac{q-2}{q}}| B_{2r}|^{\frac{2s}{N}-\frac{q-2}{q}}\\
&\leq 4\|b\|_{\infty}(\beta-1)\Theta_0^{\frac{q-2}{q}}| B_{2r}|^{\frac{2s}{N}-\frac{q-2}{q}}\leq \frac{\kappa_{N,s}}{2},
\end{align*}
it follows that
\begin{align*}
\frac{\kappa_{N,s}}{2}
\left(
\int_{ B_{2r}}
(\eta v^{\frac{\beta}{2}})^{2^*_s}\, dx
\right)^\frac{2}{2^*_s}\leq 10\int_{B^+_{2r}} t^{1-2s}|\nabla \eta|^2  V^{\beta}
\,  dxdt,
\end{align*}
which implies that there exists $\tilde C=\tilde C_{N,s,\beta}$ such that 
\begin{align*}\|v\|_{L^{\frac{\beta \,2^*_s}{2}}(\Omega \cap B_{r}(x_0))}
\leq 
\left(
\int_{ B_{2r}}
(\eta v^{\frac{\beta}{2}})^{2^*_s}\, dx
\right)^\frac{2}{2^*_s \beta}\leq 
\frac{C}{r^\frac{2}{\beta}}\left(\int_{B^+_{2r}} t^{1-2s} V^{\beta}
\,  dxdt\right)^\frac{1}{\beta}\leq 
\frac{C}{r^\frac{2}{\beta}}\left(\int_{Q_{2r}} t^{1-2s} V^{\beta}
\,  dxdt\right)^\frac{1}{\beta}.
\end{align*}
The claim now follows from Lemma~\ref{est:lem}.
\end{proof}

\subsubsection{Weak solutions of degenerate elliptic equations}

Now we collect some regularity results from \cite{JLX14}. Let $D\subset \mathbb{R}^{N+1}_+$  be an open set with $\partial' D \neq \emptyset$ (as defined in \eqref{boundaries}) and let $a:\partial' D\to \R$ be a suitable function.  Next we introduce some definitions of weak solutions depending on the problem.

\textbf{In an open set without boundary conditions:} Consider
\begin{equation}\label{2.7}
\begin{cases}
\operatorname{div}(t^{1-2s} \nabla V(x,t)) = 0 & \text{in } D, \\
-\displaystyle\lim_{t \to 0^+} t^{1-2s} \partial_t V(x,t) = k_s a(x) V(x,0) & \text{on } \partial' D.
\end{cases}
\end{equation}

We say that $V \in H(t^{1-2s}, D)$ is a \emph{weak solution} (resp.\ \emph{supersolution}, \emph{subsolution}) of \eqref{2.7} in $D$ if 
\begin{equation*}
\int_D t^{1-2s} \nabla V \nabla \Phi 
= (\text{resp. } \ge, \le) \int_{\partial' D} aV\Phi\qquad \text{for every (nonnegative) $\Phi \in C_c^\infty(D \cup \partial' D)$.}
\end{equation*}

\textbf{In an open set with boundary conditions:} Let $A\subset \R^N$ be a Lipschitz open subset with $A\cap \partial' D\neq \emptyset$
and recall the definition of 
$D_A^{1,2}\left(t^{1-2 s}, \mathbb{R}_{+}^{N+1} \right)$
given in \eqref{DAspace}. We say that $V \in D_A^{1,2}\left(t^{1-2 s}, \mathbb{R}_{+}^{N+1} \right)$ is a weak solution of
\begin{equation*}
\begin{cases}
\operatorname{div}\left(t^{1-2s} \nabla V(x,t)\right)=0 & \text { for } (x,t)\in D,\\
-\lim_{t \rightarrow 0^{+}} t^{1-2s} \partial_t V(x, t) =  k_s a(x) V(x, 0) & \text { for }x\in \partial' D\cap A,\\
V(x, 0)=0 & \text{ for }x\in \partial' D\backslash A,
\end{cases}
\end{equation*}
if 
\begin{equation*}
\int_D t^{1-2s} \nabla V \nabla \Phi 
= \int_{\partial' D} aV\Phi\qquad \text{for every (nonnegative) $\Phi \in C_c^\infty(D \cup (\partial' D\cap A))$.}
\end{equation*}

\textbf{In the whole space:} We say that $V \in H_{\mathrm{loc}}(t^{1-2 s}, \overline{\mathbb{R}^{N+1}_+})$ is a \emph{weak solution} of 
\begin{equation}\label{crit}
\begin{cases}
\operatorname{div}(t^{1-2s} \nabla V(x,t)) = 0 & \text{in } \mathbb{R}^{N+1}_+, \\
-\displaystyle\lim_{t \to 0} t^{1-2s} \partial_t V(x,t) = k_s V(x,0)^{\frac{N+2s}{N-2s}} & \text{ for }x\in\mathbb{R}^N,
\end{cases}
\end{equation}
if 
\begin{equation*}
\int_{ \mathbb{R}^{N+1}_+} t^{1-2s} \nabla V \nabla \Phi 
= \int_{ \mathbb{R}^{N}} k_s V(x,0)^{\frac{N+2s}{N-2s}}\Phi(x)\, dx\qquad \text{for every $\Phi \in C_c^\infty( \overline{ \mathbb{R}^{N+1}_+} )$.  }
\end{equation*}
Observe that, in this setting, being a weak solution is a local property and we do not require that weak solutions belong to $H(t^{1-2 s}, \overline{\mathbb{R}^{N+1}_+}).$

\begin{theorem}[{\cite[Theorem 1.8]{JLX14}}]\label{1.8}
Let $V \in H_{\mathrm{loc}}(t^{1-2 s}, \overline{\mathbb{R}^{N+1}_+})$, with $V \ge 0$ in $\mathbb{R}^{N+1}_+$ and $V \not\equiv 0$, be a weak solution of \eqref{crit}. Then, there is $t>0$ and $a\in \rn$ such that 
\begin{align*}
V(x,0)=U_{t,a}(x)\qquad \text{ for all $ x\in \rn$,}
\end{align*}
where $U_{t,a}$ is given in \eqref{tutteLeSoluzEqLimite}.   
\end{theorem}

\begin{proposition}[{\cite[Proposition 2.6]{JLX14}}]\label{2.6}
\label{prop:yyLiProp2.6} Let $r>0$ and $a \in L^q(B_r)$ for some $q > \frac{N}{2s}$.
\begin{itemize}
\item[(i)] Let $V\in H(t^{1-2s}, Q_r)$ be a weak subsolution of \eqref{2.7} in $Q_r$ (defined in \eqref{Qdef}). Then, for all $\nu > 0$,
\[
\sup_{Q_{r/2}} V^+ \le C \left( \| V^+\|_{L^\nu(t^{1-2s}, Q_r)} \right),
\]
where $V^+ = \max(0,V)$, and $C > 0$ depends only on $N, s, q, \nu, r,$ and $\|a^+\|_{L^q(B_r)}$.
\item[(ii)] Let $V \in H(t^{1-2s}, Q_r)$ be a nonnegative weak solution of \eqref{2.7} in $Q_r$. Then we have the following Harnack inequality:
\begin{equation*}
\sup_{Q_{r/2}} V \le C \inf_{Q_{r/2}} V,
\end{equation*}
where $C > 0$ depends only on $N,s,q, \|a\|_{L^q(B_r)}$, and $r$. Consequently, there exists $\alpha \in (0,1)$ depending only on $N,s,q,\|a\|_{L^q(B_r)}$ such that any weak solution $V$ of \eqref{2.7} is in $C^\alpha(\overline{Q}_{r/2})$. Moreover,
\[
\|V\|_{C^\alpha(\overline{Q}_{r/2})} \le C \left( \|V\|_{L^\infty(Q_r)} \right),
\]
where $C > 0$ depends only on $N,s,q, \|a\|_{L^q(B_r)},$ and $r$.
\item[(iii)] 
Let $A\subset \R^N$ be a Lipschitz open subset with $A\cap B_r\neq \emptyset$ and let $V \in D_A^{1,2}\left(t^{1-2 s}, \mathbb{R}_{+}^{N+1} \right)$ be a weak solution of
\[
\begin{aligned}
\operatorname{div}\left(t^{1-2s} \nabla V(x,t)\right)&=0 &&\ \text { for } (x,t)\in Q_r,
\\
-\lim_{t \rightarrow 0^{+}} t^{1-2s} \partial_t V(x, t) &=  k_s a(x) V(x, 0) &&\ \text { for }x\in B_r\cap A,
\\
V(x, 0)&=0 && \ \text{ for }x\in B_r\backslash A,
\end{aligned}
\]
where $k_s$ is defined in \eqref{dsdef}.  Then,
$$
\sup _{Q_{\frac r2}} V^{+} \leq C \left(\int_{Q_r}t^{1-2s}(V^+)^{2^*_s}dx dt\right)^{\frac1{2^*_s}},
$$
where $V^{+}=\max (0, V)$ and $C>0$ depends only on $N, s, q, r$, and $\left\|a^{+}\right\|_{L^{q}\left(B_r\right)}$.
\end{itemize}    
\end{proposition}
\begin{proof}
This result follows from \cite[Proposition 2.6]{JLX14}, where only the case $r=1$ and $A=B_1(0)$ is considered, but the same proof also covers the case $r\neq 1$ and $A\neq B_1(0)$ with Dirichlet exterior boundary conditions.  Via a rescaling, the general case for $r>0$ can be easily derived.
\end{proof}

\begin{proposition}[{\cite[Proposition 2.19]{JLX14}}]\label{prop:Li2.19}
Let $\eps_0>0$, $\eta>0$, $ q\in\R^N,$ and let $\tilde K\in C^1(B_{2\eta}( q))$. Let $\tilde \Phi\in H(t^{1-2s}, Q_{2\eta}(q))$, satisfying $\tilde \Phi\geq0$ in $Q_{2\eta}( q)$, be a weak solution of
\begin{equation}\label{eq_Phi_tilde}
\left\{
\begin{array}{ll}
\operatorname{div}\left(t^{1-2 s} \nabla \tilde \Phi\right)=0
&\text{ in }Q_{2\eta}( q),\\
-\lim _{t \rightarrow 0} t^{1-2 s} \partial_t \tilde \Phi(y,t)=\tilde K(y)\tilde \Phi^p(y,0)&\text{ for }y\in B_{2\eta}( q),
\end{array}
\right.
\end{equation}
where $p\in(\frac{N+2s}{N-2s}-\eps_0,\frac{N+2s}{N-2s})$. Then there exist $C>0$ and $\alpha\in(0,1)$, both depending only on $N$, $s$, $\eta$, $\eps_0$, $\|\tilde \Phi\|_{L^\infty(Q_{2\eta}( q))},$ and $\|\tilde K\|_{C^1(B_{2\eta}}(q))$, such that $\nabla_x \tilde \Phi$ and $t^{1-2s}\partial_t\tilde \Phi$ are in $C^\alpha(\overline{Q_{\eta}( q)})$ and
\[
\|\nabla_x \tilde \Phi\|_{C^\alpha(\overline{Q_{\eta}( q)})}+\|t^{1-2s}\partial_t\tilde \Phi\|_{C^\alpha(\overline{Q_{\eta}( q)})}\leq C.
\]

\end{proposition}
\begin{proof}
This result follows from \cite[Proposition 2.19]{JLX14}, where only the case $\eta=\frac12$ is considered and observing that the dependence of $C$ on $p$ can be reduced to the dependence on $\eps_0$. Via a rescaling, the general case for $\eta>0$ can be easily derived.
\end{proof}

\begin{lemma}[{\cite[Lemma 4.10]{JLX14}}]\label{lemma:Li4.10}
Let $r>0$ and let $F$ be a function on $B_r^+$ such that,
for every $\gamma\in(0,r)$,
$F\in H(t^{1-2s},B_r^+\setminus\overline{B_\gamma^+})$,
$F>0$ in $B_r^+\setminus\overline{B_\gamma^+}$,
and $F$ is a weak solution of
\[
\left\{
\begin{array}{ll}
\operatorname{div}\left(t^{1-2 s} \nabla F\right)=0
&\text{ in }B^+_r\setminus\overline{B^+_\gamma},\\
-\lim _{t \rightarrow 0} t^{1-2 s} \partial_t F(y, t)=0 &\text{ for }y\in B_r\setminus\overline{B_\gamma}.
\end{array}
\right.
\]
Then,
\[
F(X)=A|(y,t)|^{2s-N}+W(y,t),
\]
where $y\in\R^N$ and $t\in \R$, $A$ is a nonnegative constant and $W\in  H(t^{1-2s},B^+_r)$ satisfies
\[
\left\{
\begin{array}{ll}
\operatorname{div}\left(t^{1-2 s} \nabla W\right)=0
&\text{ in }B^+_r,\\
-\lim _{t \rightarrow 0} t^{1-2 s} \partial_t W(y,t)=0 &\text{ for }y\in B_r.
\end{array}
\right.
\]
\end{lemma}


We conclude stating a Pohozaev identity.
\begin{proposition}[{\cite[Proposition 4.7]{JLX14}}]\label{prop:Li4.7}
Let $\delta>0$ and let $\Phi\in H(t^{1-2s},B_+^{2\delta})$ with $\Phi\geq0$ in $B_+^{2\delta}$
be a weak solution of
\begin{equation*}
\left\{
\begin{array}{ll}
\operatorname{div}\left(t^{1-2 s} \nabla \Phi\right)=0
&\text{ in }B^+_{2\delta},\\
-\lim _{t \rightarrow 0} t^{1-2 s} \partial_t \Phi(y,t)= k_s \Phi^p(x,0)&\text{ for }y\in B_{2\delta},
\end{array}
\right.
\end{equation*}
where $p>0$ and $k_s$ is defined in \eqref{dsdef}. Then,
\begin{equation}\label{Poh_for_ext}
\int_{B_\delta}k_sE_1(y,0,\Phi,\nabla\Phi,\delta,s)\,dy+\int_{\partial'' B_{\delta}^+}t^{1-2s}E_2(y,t,\Phi,\nabla\Phi,\delta,s)\,dydt=0,
\end{equation}
where $y\in\R^N$, $t\geq0$, and
\begin{align}
E_1(y,t,\Phi,\nabla\Phi,\delta,s)&:=\frac{N-2s}{2}\Phi^{p+1}+\langle (y,t),\nabla \Phi\rangle \Phi^p,\label{def:E1}\\
E_2(y,t,\Phi,\nabla\Phi,\delta,s)&:=\frac{N-2s}{2}\Phi\frac{\partial \Phi}{\partial\nu}-\frac \delta2|\nabla \Phi|^2+\delta\left|\frac{\partial \Phi}{\partial\nu}\right|^2. \label{def:E2}
\end{align}
\end{proposition}

\subsection{Bubble and limit problem}\label{section:SobolevConstant_and_bubbles}

Let $S_{N,s}$ denote the best constant in the fractional Sobolev inequality $S_{N,s}\|v\|_{2^*_s}^2\leq \|v\|^2_{\mathcal D^s}$, for  $v\in \mathcal D_0^{s}(\Omega)$. Since the ratio $\frac{\|v\|_{\mathcal D^s}^2}{\|v\|_{2^*_s}^2}$ is invariant under dilations,  $S_{N,s}$
is independent of $\Omega$ (and depends only on $N$ and $s$), namely, for $\Omega\subseteq\mathbb R^N$, 
\begin{equation}\label{SNS}
S_{N,s}:=\inf_{v\in \mathcal D_0^{s}(\Omega)\setminus\{0\}}
\frac{\|v\|_{\mathcal D^s}^2}{\|v\|_{2^*_s}^2}.    
\end{equation}

When $\Omega=\mathbb R^N$,  it is well known that the minimum $S_{N,s}$ is achieved and that the 
only minimizers are the  functions
\begin{equation}\label{minimizers}x\mapsto\frac{c}{\left(\lambda^2+|x-x_0|^2\right)^{\frac{N-2s}{2}}},
\qquad \text{where $c,\lambda\in\mathbb R\setminus\{0\}$, $x_0\in\mathbb R^N$,}
\end{equation}
see \cite{Lieb}. As a consequence, $S_{N,s}$ is never achieved when $\Omega\neq\mathbb R^N$. We also recall that $S_{N,s}$ has been explicitly computed 
\begin{equation}\label{constant:SnsVALUE}S_{N,s}
=
\frac{\pi^{s}}{(\tau_{N,s})^s}
\left[\frac{\Gamma\left(\frac{N}{2}\right)}{\Gamma\left(N\right)}\right]^{2s/N}
=2^{2s}\pi^{s}\frac{\Gamma\left(\frac{N+2s}{2}\right)}{\Gamma\left(\frac{N-2s}{2}\right)}
\left[\frac{\Gamma\left(\frac{N}{2}\right)}{\Gamma\left(N\right)}\right]^{2s/N},
\end{equation}
where 
\begin{equation}\label{def:constanTau}
\tau_{N,s}:=\frac{1}{4}\left[\frac{\Gamma\left(\frac{N-2s}{2}\right)}{\Gamma\left(\frac{N+2s}{2}\right)}\right]^{\frac{1}{
s}},
\end{equation}
see \cite{CotsiolisTavoularis}. Furthermore, choosing the minimizer in \eqref{minimizers} with $x_0=0$, $\lambda=\tau_{N,s}^{-\frac{1}{2}} 
$ and $c=\lambda^{N-2s}$, namely 
\begin{equation}\label{Bubble}
U(x):=\frac{1}{\left(1+\tau_{N,s}|x|^2\right)^{\frac{N-2s}{2}}},\end{equation}
 it is not difficult to see that $U$ is radial, $U\leq 1$, $U(0)=1$, and $U$
solves
\begin{equation}\label{equationBubble}
\left\{
\begin{array}{ll}
(-\Delta)^s U=U^{\frac{N+2s}{N-2s}} & \mbox{in }\mathbb R^N,\\
 U>0 & \mbox{in }\mathbb R^N.
\end{array}
\right.\end{equation}
Moreover, by \cite{CLO}, the unique regular solutions of \eqref{equationBubble} are
\begin{equation}\label{tutteLeSoluzEqLimite}
U_{t,a}(x):=t^{-\frac{N-2s}{2}}U\left(\frac{x-a}{t}\right)\qquad \text{ for any } a \in\mathbb R^N,\, t>0.\end{equation}
Let us observe that, for any $\eps\geq 0$ small, by direct computations, one has
\begin{equation}\label{U_in_L2*-1}
    \int_{\mathbb R^N} U^{p_\eps}(y)\, dy <+\infty.
\end{equation}
In particular, for $\eps=0$, 
\begin{equation}
\int_{\mathbb R^N} U^{\frac{N+2s}{N-2s}}(y)\,dy
=\sigma_N\int_{0}^{+\infty}\frac{\rho^{N-1}}{\left(1+\tau_{N,s}\rho^2\right)^{\frac{N+2s}{2}}}\, d\rho =\frac{\sigma_N}{2} (\tau_{N,s})^{-\frac{N}{2}} B\Big(\frac{N}{2},s\Big),
\label{ExactValue}
\end{equation} 
where the constant $\tau_{N,s}$ is defined in \eqref{def:constanTau}, $\sigma_N=|\partial B_1|=\frac{2\pi^\frac{N}{2}}{\Gamma(\frac{N}{2})},$ and $B(a,b)=\frac{\Gamma (a)\Gamma (b)}{\Gamma (a+b)}$ is the Beta function.
Similarly, 
\begin{align}
\int_{\mathbb R^N} U^{\frac{2N}{N-2s}}(y)\,dy
=\sigma_N\int_{0}^{+\infty}\frac{\rho^{N-1}}{\left(1+\tau_{N,s}\rho^2\right)^{N}}\, d\rho
 = \pi^{\frac{N}{2}}  \frac{\Gamma\left(\frac{N}{2}\right)}{\Gamma(N)}\tau_{N,s}^{-\frac{N}{2}}= \left(S_{N,s}\right)^{\frac{N}{2s}}.
\label{ExactValue2}
\end{align}

\subsection{Known results on blow up points}

In this subsection we collect some definitions and results from \cite{JLX14}. We write them according to our notation and for the special case of our interest.  Let $u_\eps\in L^\infty(\Omega)\cap \mathcal D^s(\mathbb{R}^N)$, $u_\eps\geq0$ in $\mathbb{R}^N$ be a family of solutions to 
\[
(-\Delta)^s u_\eps=u_\eps^{p_\eps}\qquad\text{in $\Omega$,}
\]
where $p_\eps$ is as in \eqref{defp_eps}.

\begin{definition}
  We say that ${\bar x_\eps}\in\Omega$ is a blowing up sequence if $u_\eps({\bar x_\eps})\to+\infty.$ We say that $\bar x\in \overline\Omega$ is a blow up point of $u_\eps$ if there exists a blowing up sequence $\bar x_\eps$ such that $\bar x_\eps\to\bar x.$ 
\end{definition}

\begin{definition}\label{def:isolated_bup}
We say that  $\bar x\in\Omega$ is an isolated blow up point of $u_\eps$ if there exist $r_0\in(0,\operatorname{dist}(\bar x,\partial\Omega))$, $C_0>0,$ and a blowing up sequence ${\bar x_\eps}$ converging to $\bar x$ such that ${\bar x_\eps}$ is a local maximum point of $u_\eps$, and 
\begin{align}\label{almostisolated}
u_\eps(x)\leq\frac{C_0}{|x-{\bar x_\eps}|^\frac{2s}{p_\eps-1}}\qquad \text{ for all } x\in B_{r_0}({\bar x_\eps})\setminus \{\bar x_\eps\}.
\end{align}
\end{definition}

\begin{remark}\label{rem:r1C1}
If $\bar x$ is an isolated blow up point of $u_\eps$   with blowing up sequence  $\bar x_\eps$, then there exist $r_1\in(0,\operatorname{dist}(\bar x, \partial\Omega))$, $C_1>0$ such that 
\begin{equation}\label{C1}
u_\eps(x)\leq \frac{C_1}{|x-{\bar x_\eps}|^\frac{2s}{p_\eps-1}}\qquad \text{ for all } \, x\in B_{r_1}(\bar x)\backslash \{\bar x_\eps\}.
\end{equation}
\end{remark}

Let $\bar x$ be an isolated blow up point of $u_\eps$ with blowing up sequence $\bar x_\eps$, and define
\begin{align}\label{not}
\bar u_\eps(r):=\frac{1}{|\partial  B_r|}\int_{\partial B_r({\bar x_\eps})}u_\eps,\qquad r>0,    
\end{align}
and 
\[
\overline \psi_\eps(r):=r^{\frac{2s}{p_\eps-1}}\bar u_\eps(r),\qquad r>0.
\]

\begin{definition}\label{def:isol_simple_bup}
We say that $\bar x\in\Omega$ is an isolated and simple blow up point of $u_\eps$, if $\bar x$ is an isolated blow up point
such that, for some $\rho>0$ (independent of $\eps$), $\overline \psi_\eps$ has precisely one critical point in $(0,\rho)$ for small $\eps$.
\end{definition}

\begin{remark}
    By the translation invariance of the problem, without loss of generality, we may always assume that  $\Omega\subset\R^N$, $N\geq 2$, is a smooth bounded open set containing the origin and that $\bar x=0$.
\end{remark}

\subsubsection{Local estimates for isolated,
and isolated and simple, blow up points}

In this section we recall some results from \cite{JLX14} that are useful to show that isolated blow up points are simple (see Theorem~\ref{prop:Li4.16}) and to prove Theorem~\ref{PositiveAreLow:thm}.  

For all the results in this subsection, we let $A\subset \R^N$ be an open set and let $u_\eps$ be a family of nonnegative solutions of $(-\Delta)^s u_\eps =u_\eps^{p_\eps}$ in $A.$

\begin{lemma}[{Harnack inequality, \cite[Lemma 4.3]{JLX14}}]\label{lemma:Li4.3}
 Let us suppose that $\bar x$ is an isolated blow up point of $u_\eps$, with blowing up sequence $\bar x_\eps$, let $r_1$ and $C_1$ be those in \eqref{C1} and let $U_\eps$ be the local extension of $u_\eps$ (see Subsection~\ref{subsec:ext}). Then, for any $r\in(0,\frac13 r_1)$, 
 \[
 \sup_{B^+_{2r}({\bar x_\eps})\setminus \overline {B^+_{\frac r2}({\bar x_\eps})}} U_\eps \leq C\inf_{B^+_{2r}({\bar x_\eps})\setminus \overline {B^+_{\frac r2}({\bar x_\eps})}} U_\eps,
 \]
 where $C$ is a positive constant depending only on $N$, $s$, $r_1$, and $C_1$.
 \end{lemma}

\begin{proposition}[{\cite[Proposition 4.4]{JLX14}}]\label{prop:Li4.4}
 Suppose that $\bar x$ is an isolated blow up point of $u_\eps$, with blowing up sequence $\bar x_\eps$.  Then for any $R_\eps\to+\infty$ and $\gamma_\eps\to 0^+$ we have, after passing to a subsequence (keeping the notations $u_\eps$ and ${\bar x_\eps}$), 
\begin{align*}
\|u_\eps^{-1}({\bar x_\eps})u_\eps({\bar x_\eps}+u_\eps^{-\frac{p_\eps-1}{2s}}({\bar x_\eps})\, \cdot\, )-U\|_{C^2(B_{2R_\eps}(0))}\leq \gamma_\eps,
\end{align*}
\[
R_\eps u_\eps^{-\frac{p_\eps-1}{2s}}({\bar x_\eps})\to 0 \qquad \text{as $\eps\to0$},
\]
where $U$ is defined in \eqref{Bubble}.
\end{proposition}

\begin{remark}\label{blowup:rmk}
Observe that  the blowing up sequence $\bar x_\eps$ in Definition~\ref{def:isolated_bup}  consists of local maximum points for $u_\eps$.  However, if $\bar x$ is an isolated blow up point of $u_\eps$, and if  $\bar x_\eps\rightarrow \bar x$ is a blowing up sequence which  satisfies the local inequality \eqref{almostisolated} (but may not be of local maxima), then it is not difficult to show, following the proofs in \cite{JLX14}, that there exists $t>0$ and $a\in \mathbb R^N$ such that 
for any $R_\eps\to+\infty$ and $\gamma_\eps\to 0^+$, up to a subsequence,
\[
\|u_\eps^{-1}({\bar x_\eps})u_\eps({\bar x_\eps}+u_\eps^{-\frac{p_\eps-1}{2s}}({\bar x_\eps})\,\cdot\,)-U_{t,a}(\cdot)\|_{C^2(B_{2R_\eps}(0))}\leq \gamma_\eps,
\]
where $U_{t,a}$ is defined in \eqref{tutteLeSoluzEqLimite}.
Furthermore, since $U_{t,a}(0)=1$ then $t\leq 1$ and $|a|^2=\frac{t(1-t)}{\tau_{N,s}}\leq (\tau_{N,s})^{-1}.$
\end{remark}

\begin{proposition}[{\cite[Proposition 4.5]{JLX14}}]\label{prop:Li4.5}
 Suppose that $\bar x$ is an isolated blow up point of $u_\eps$, with blowing up sequence $\bar x_\eps$ and let $C_1$ the constant in \eqref{C1}.  Then there is $c=c(N,s,C_1)>0$ such that, for any $e\in\R^N$, $|e|=1$, and for any $\eps$ sufficiently small
\[
u_\eps(\bar x_\eps) u_\eps(\bar x_\eps+e)\geq c.
\]
\end{proposition}


\begin{lemma}[{\cite[Lemma 4.8]{JLX14}}]\label{lemma:Li4.8}
 Suppose that $\bar x=0$ is an isolated and simple blow up point of $u_\eps$, with blowing up sequence $\bar x_\eps$. Then, 
 \begin{equation}\label{property+1}
\eps=O\!\left(
u_\eps(\bar x_\eps)^{
-\min\left\{2,\frac{2}{N-2s}\right\}+o(1)}
\right),
\qquad\text{and hence}\qquad
u_\eps(\bar x_\eps)^\eps=1+o(1).
\end{equation}
\end{lemma}
This follows from the proof of \cite[Lemma 4.8]{JLX14},
retaining both error terms in estimates (4.14)--(4.17).

\begin{proposition}[{\cite[Proposition 4.9]{JLX14}}]\label{prop:Li4.9}
 Suppose that $\bar x=0$ is an isolated and simple blow up point of $u_\eps$, with blowing up sequence $\bar x_\eps$. Then there exist $C>0$ such that, for any $\eps$ sufficiently small, 
\begin{align*}
u_\eps(x)\leq \frac{C}{u_\eps({\bar x_\eps})|x-{\bar x_\eps}|^{N-2s}}\qquad \text{ for all }\, x\in \overline{B_{1}({\bar x_\eps})}\backslash \{\bar x_\eps\}.
\end{align*}
\end{proposition}

\subsubsection{Isolated blow up points are simple}

In this section we prove that every isolated blow-up point of a solution to \eqref{equn}$-$\eqref{defp_eps} is simple. Although this property is established in \cite[Proposition 4.16]{JLX14}, our equation does not fall directly within the scope of their result. Therefore, for the reader's convenience, we present a detailed proof following the same strategy as in \cite[Proposition 4.16]{JLX14}.

\begin{proposition}[{\cite[Proposition 4.16]{JLX14}}]\label{prop:Li4.16}
Let $u_\eps$ be a family of solutions of \eqref{equn}$-$\eqref{defp_eps}. If ${\bar x}= 0$ is an isolated blow up point of $u_\eps$, then ${\bar x}= 0$ is an isolated and simple blow up point.
\end{proposition}
\begin{proof} 
Let $\bar x_\eps$ be the sequence in Definition~\ref{def:isolated_bup} and $r_1>0$ and $C_1>0$ the constants of Remark~\ref{rem:r1C1}; namely, 
\begin{equation}\label{C1def}
u_\eps(x)\leq \frac{C_1}{|x-{\bar x_\eps}|^\frac{2s}{p_\eps-1}}\qquad \text{ for all } \, x\in B_{r_1}(0)\backslash\{x_\eps\}.
\end{equation}
Let   $M_\eps:=u_\eps(\bar x_\eps)$, and $R_\eps\to \infty$ be an arbitrary sequence. Passing to a subsequence (we omit the relabeling), we may assume that
\begin{align}\label{repsdef}
r_\eps:= \frac{R_\eps}{M_\eps^{\frac{p_\eps-1}{2s}}}\to 0.     
\end{align}
Then, using the notation \eqref{not},
\begin{align*}
\overline  \psi_\eps(r):=r^{\frac{2s}{p_\eps-1}}\bar u_\eps(r)
=r^{\frac{2s}{p_\eps-1}} \frac{1}{|\partial  B_r|}\int_{\partial B_r({\bar x_\eps})}u_\eps
=\frac{1} {\sigma_N}\int_{\partial B_1} r^{\frac{2s}{p_\eps-1}}u_\eps(\bar x_\eps+r x)\, d\sigma_x.
\end{align*}
In particular, 
\[\overline  \psi_\eps(M_\eps^{-\frac{p_\eps-1}{2s}} r)
=
\frac{1 }{\sigma_N} \int_{\partial B_1}r^{\frac{2s}{p_\eps-1}}
M_\eps^{-1} u_\eps(\bar x_\eps+ M_\eps^{-\frac{p_\eps-1}{2s}} rx)
\, d\sigma_x=
\frac{1 }{\sigma_N} \int_{\partial B_1} r^{\frac{2s}{p_\eps-1}}
 v_\eps(rx)
\, d\sigma_x
,\]  where $v_\eps(y):=M_\eps^{-1} u_\eps(\bar x_\eps+ M_\eps^{-\frac{p_\eps-1}{2s}} y)$. Hence, 
\begin{align*}
\partial_r (\overline  \psi_\eps(M_\eps^{-\frac{p_\eps-1}{2s}} r))
=
\frac{1}{\sigma_N} \int_{\partial B_1}\partial_r(r^{\frac{2s}{p_\eps-1}}
v_\eps(rx) )\, d\sigma_x.
\end{align*}
By Proposition~\ref{prop:Li4.4}, 
\begin{align}\label{vepsconv}
\|v_\eps-U\|_{C^2(B_{2R_\eps}(0))}\to 0\qquad \text{as $\eps\to 0.$ }
\end{align}

The critical points of
$r\mapsto\overline\psi_\eps
(M_\eps^{-\frac{p_\eps-1}{2s}}r)$
are precisely the positive zeros of
\[
\frac{1}{\sigma_N}\int_{\partial B_1}
\left(
\frac{2s}{p_\eps-1}v_\eps(rx)
+r\nabla v_\eps(rx)\cdot x
\right)\,d\sigma_x.
\]
By \eqref{vepsconv}, this expression converges in $C^1([0,R])$,
for every fixed $R>0$, to
\[
\frac{N-2s}{2}
\frac{1-\tau_{N,s}r^2}
{(1+\tau_{N,s}r^2)^{\frac{N-2s}{2}+1}}.
\]
The limit is positive before $r_*=\tau_{N,s}^{-1/2}$,
negative afterward, and has strictly negative derivative at $r_*$.
The $C^1$ convergence therefore implies that, for every fixed
$R>r_*$ and sufficiently small $\eps$, the rescaled average
has exactly one critical point in $(0,R)$.
This remains true on $(0,R_\eps)$, while \eqref{repsdef} still
holds. Returning to the original radial variable, we obtain that

\begin{align}\label{ki}
\text{ $r \mapsto \overline  \psi_\eps(r)$ has precisely one critical point in the interval $(0,r_\eps)$. }
\end{align}
Next, assume, by contradiction, that the isolated blow up point $\bar x=0$  of $u_\eps$ is not an isolated and simple blow up point.  Then, by Definition \ref{def:isol_simple_bup} and \eqref{ki}, there exists a second critical point $\mu_\eps$ of $\overline  \psi_\eps(r)$ such that $\mu_\eps\geq r_\eps$ and $\lim_{\eps\to0}\mu_\eps=0.$ 
Let  
\begin{equation}\label{def_phi_eps}
\phi_\eps(y):=\mu_\eps^{\frac{2s}{p_\eps-1}}u_\eps(\mu_\eps y+\bar x_\eps)\qquad \text{ for  } y\in\R^N.
\end{equation}
By scalings, we know that $\phi_\eps$ is a nonnegative function which solves  $(-\Delta)^s \phi_\eps = \phi^{p_\eps}_\eps$ in $\mu_\eps^{-1}(\Omega - \bar x_\eps)$, moreover $0$ is a local maximum point of $\phi_\eps$ and, by \eqref{C1def},
\begin{align}\label{bounphiepsutile}
    \phi_\eps(y)\leq \frac{C_1}{|y|^{\frac{2s}{p_\eps-1}}} \quad \mbox{ for }y\in B_{\frac{r_1}{\mu_\eps}}(0)\backslash \{0\}.
\end{align}
Using that $\mu_\eps\geq r_\eps$ and \eqref{repsdef}, it also follows that
\begin{align}\label{phi:infty}
\lim\limits_{\eps\to0}\phi_\eps(0)
=\lim\limits_{\eps\to0} \mu_\eps^{\frac{2s}{p_\eps-1}}u_\eps(\bar x_\eps )
\geq \lim\limits_{\eps\to0} (r_\eps M_\eps^{\frac{p_\eps-1}{2s}})^{\frac{2s}{p_\eps-1}}
=\lim\limits_{\eps\to0} R_\eps^{\frac{2s}{p_\eps-1}}=\infty,
\end{align}
namely we have proved that $0$ is an isolated blow up point of $\phi_\eps$ (with blowing up sequence the one  identically equal to $0$).
Let us consider $\bar\phi_\eps(r):=|\partial B_r|^{-1}\int_{\partial B_r} \phi_\eps$. Moreover, by \eqref{ki} (recalling that $\overline  \psi_\eps(r)=r^{\frac{2s}{p_\eps-1}}\bar u_\eps(r)$), it follows that
\begin{equation}\label{phi_eps_0_to_0}
r^{\frac{2s}{p_\eps-1}}\bar \phi_\eps(r)
=(\mu_\eps r)^{\frac{2s}{p_\eps-1}}\bar u_\eps(\mu_\eps r)
\text{ has precisely one critical point in $(0,1)$}
\quad \text{ and }\quad  \frac{d}{dr}\{r^{\frac{2s}{p_\eps-1}}\bar \phi_\eps(r)\}{\mathlarger{\mathlarger{\mathlarger{|}}}}_{r=1}=0.
\end{equation}
Therefore $0$ is an isolated and simple blow up point of $\phi_\eps$, with blowing up sequence identically equal to $0$.

Hence, by Propositions~\ref{prop:Li4.5} and~\ref{prop:Li4.9}, there are $c,C>0$ such that, for any $e\in\R^N$, $|e|=1$, and for any $\eps$ sufficiently small,
\begin{equation}\label{phi:bd}
c\leq\phi_\eps(0)\phi_\eps(e)\leq C.
\end{equation}

 Let $\Phi_\eps$ be the local extension of $\phi_\eps$ (see Section~\ref{subsec:ext}). Then, 
 \begin{align*}
 \operatorname{div}(t^{1-2s} \nabla \Phi_\eps) = 0\ \text{in } \mathbb{R}^{N+1}_+, \quad 
 \Phi_\eps(x,0) = \phi_\eps(x)\  \text{ for }x\in \R^N,\quad 
 -\displaystyle\lim_{t \to 0} t^{1-2s} \partial_t \Phi_\eps(x,t) = k_s \phi_\eps(x)^{p_\eps} \ \text{ for }x\in \mu_\eps^{-1}(\Omega-\bar x_\eps),
 \end{align*}
  where $k_s$ is defined in \eqref{dsdef}.
Let us fix $r_1'\in (0, \frac{r_1}{\mu_\eps})$,  then by Lemma~\ref{lemma:Li4.3}, there exists  a positive constant $C_1'$ independent of $\eps$ such that
  \begin{equation}\label{supbd}
\sup_{B^+_{2r}\setminus\overline{B^+_{r/2}}}
\Phi_\eps(0,0)\Phi_\eps
\le C\left(\frac{2}{r}\right)^{N-2s},
\qquad
0<r<\min\left\{\frac{r_1'}{3},\frac12\right\}.
\end{equation}
The estimate follows from Proposition~\ref{prop:Li4.9} and
Lemma~\ref{lemma:Li4.3}. Repeated application of
Lemma~\ref{lemma:Li4.3} along overlapping annuli gives uniform
bounds for $\Phi_\eps(0,0)\Phi_\eps$ on every compact subset of
$\overline{\mathbb R^{N+1}_+}\setminus\{0\}$.

 Let $\eps_0>0$, $\eta>0$, and $q\in\R^N$ such that $0\not\in B_\eta(q)$. By \eqref{supbd}, we have that $\|\Phi_\eps(0,0)\Phi_\eps\|_{L^\infty(Q_{2\eta}( q))}$ is bounded independently of $\eps$. By Proposition~\ref{prop:Li2.19}, there exist $C''>0$ and $\alpha\in(0,1)$ independent of $\eps$,  such that 
 \[
\|\nabla_x \Phi_\eps(0,0)\Phi_\eps\|_{C^\alpha(\overline{Q_{\eta}( q)})}+\|t^{1-2s}\partial_t \Phi_\eps(0,0)\Phi_\eps\|_{C^\alpha(\overline{Q_{\eta}( q)})}\leq C''.
 \]
 By a diagonalization argument and Lemma~\ref{lemma:Li4.10}, we obtain that
 \begin{equation}\label{conv_Phi_tilde_eps}
 \Phi_\eps(0,0)\Phi_\eps(y,t)\to F(y,t)\quad \text{in $C^\alpha_{loc}(\overline{\R^{N+1}}\setminus\{0\})\cap C^2_{loc}(\R^{N+1}_+)$}
 \end{equation}
 and
 \begin{equation}\label{phi_to_F}
 \phi_\eps(0)\phi_\eps(y)\to F(y,0)\quad \text{in $C^2_{loc}(\R^{N}\setminus\{0\})$},
\end{equation}
 as $\eps \to0$, where
 \[F(y,t)=A|(y,t)|^{2s-N}+W(y,t),\] with $A\geq 0$, and $W(y,t)$ satisfies
 \[
 \left\{
 \begin{array}{ll}
 \operatorname{div}\left(t^{1-2 s} \nabla W\right)=0
 &\text{ in }\R^{N+1}_+,\\
 -\lim _{t \rightarrow 0} t^{1-2 s} \partial_t W(y, t)=0 &\text{ for }y\in\R^N.
  \end{array}
 \right.
 \]
 Noting that $F(y,t)$ is nonnegative, we have $\liminf_{|(y,t)|\to+\infty}W(y,t)\geq0$. It follows from the weak maximum principle and the Harnack inequality that $W(y,t)\equiv B\geq0$. Thus, being
 \[
 \frac{d}{dr}[r^{\frac{2s}{p_\eps-1}}\phi_\eps(0)\bar\phi_\eps(r)]{\mathlarger{\mathlarger{\mathlarger{|}}}}_{r=1}=0,\qquad\text{and}\qquad \frac{d}{dr}[r^{\frac{N-2s}{2}}(A r^{2s-N}+B)]{\mathlarger{\mathlarger{\mathlarger{|}}}}_{r=1}=\frac{N-2s}{2}(B-A),
 \]
 making use of \eqref{phi_to_F} we can deduce that $A=B$ and hence
 \begin{equation}\label{Fform}
     F(y,t)=A(|(y,t)|^{2s-N}+1).
 \end{equation}

 Let us now show that $A>0$. Fix $e\in\R^N$, $|e|=1$. By \eqref{phi:bd} we derive that (up to a subsequence) there exists $K>0$ such that $\phi_\eps(0)\phi_\eps(e)\to K$ as $\eps\to0$. This, combined with \eqref{phi_to_F} and \eqref{Fform}, yields $K=2A$,
 %
%
which immediately implies that $A>0$.
 
 Since $\Phi_\eps$ is a solution to 
 \[
 \left\{
 \begin{array}{ll}
 \operatorname{div}\left(t^{1-2 s} \nabla \Phi_\eps\right)=0
 &\text{ in }\R^{N+1}_+,\\
 -\lim _{t \rightarrow 0} t^{1-2 s} \partial_t \Phi_\eps(y, t)=k_s\phi_\eps^{p_\eps}(y) &\text{ for }y\in \mu_\eps^{-1}(\Omega-\bar x_\eps),
 \end{array}
\right.
 \]
it must satisfy the Pohozaev identity \eqref{Poh_for_ext} in Proposition~\ref{prop:Li4.7}.
Next we show  that, for small  $\delta>0$,
 \begin{equation}\label{E1}
\limsup\limits_{\eps\to0}\Phi_\eps^2(0,0)\int_{B_\delta} k_sE_1(y,0,\Phi_\eps,\nabla \Phi_\eps, \delta,s)\,dy\leq0,
 \end{equation}
 and
 \begin{equation}\label{E2}
{\lim\limits_{\eps\to0}}\Phi_\eps^2(0,0)\int_{\partial'' B_{\delta}^+} t^{1-2s}E_2(y,t,\Phi_\eps,\nabla \Phi_\eps, \delta,s)\,dydt<0,
 \end{equation}
 where $E_1$ and $E_2$ are defined in \eqref{def:E1} and \eqref{def:E2} respectively, which gives a contradiction to the Pohozaev identity,  thus  concluding the proof.\\
In order to prove \eqref{E2} we set $\widetilde \Phi_\eps(y,t):=\Phi_\eps(0,0)\Phi_\eps(y,t)$ and $\eta:=\frac\delta4$, then there exist $q_1,\ldots,q_\ell\in\partial B_\delta$ such that
 \[
\partial'' B_{\delta}^+\subset (\overline B_\delta\times [\eta,\delta])\cup(\cup_{i=1}^\ell Q_\eta(q_i)).
 \]
 Observing that there exists $B>0$ such that $||\widetilde \Phi_\eps||_{L^\infty(Q_{2\eta}(q_i))}\leq B$ for any $\eps$ sufficiently small, and that $\widetilde \Phi_\eps$ solves \eqref{eq_Phi_tilde} with $\widetilde K(y)=k_s\phi_\eps^{1-p_\eps}(0)$, by \eqref{phi:infty} and Proposition~\ref{prop:Li2.19} it follows that there exists $C>0$ such that
\begin{equation}\label{bound_Phi_tilde_eps}
 \|\nabla_x \widetilde \Phi_\eps\|_{C^\alpha(\overline{Q_{\eta}( q)})}+\|t^{1-2s}\partial_t\widetilde \Phi_\eps\|_{C^\alpha(\overline{Q_{\eta}( q)})}\leq C.
 \end{equation}
As a consequence, recalling the definition of $E_2$ in \eqref{def:E2}
 \begin{align*}
 \lim\limits_{\eps \to0}&\Phi_\eps^2(0,0)\int_{\partial'' B_{\delta}^+}t^{1-2s}E_2(y,t,\Phi_\eps,\nabla \Phi_\eps, \delta,s)\,dydt\\
 &=\lim\limits_{\eps \to0}\Phi_\eps^2(0,0)\int_{\partial'' B_{\delta}^+}t^{1-2s}\left(\frac{N-2s}{2}\Phi_\eps\frac{\partial\Phi_\eps}{\partial\nu}-\frac\delta2|\nabla\Phi_\eps|^2+\delta|\frac{\partial\Phi_\eps}{\partial\nu}|^2\right)\,dydt\\
 &\overset{\star}{=}\int_{\partial'' B_{\delta}^+}t^{1-2s}\left(\frac{N-2s}{2}F\frac{\partial F}{\partial\nu}-\frac\delta2|\nabla F|^2+\delta|\frac{\partial F}{\partial\nu}|^2\right)\,dydt
 =-\frac{(N-2s)^2}{2}A^2\int_{\partial'' B_{1}^+}t^{1-2s}\,dydt<0,
 \end{align*}
 where in $\star$ we have used that $t^{2s-1}$ is integrable and then that the bound \eqref{bound_Phi_tilde_eps} allows to apply the dominated convergence theorem and to pass to the limit, using \eqref{conv_Phi_tilde_eps} and \eqref{Fform}.
\\
 Finally we prove \eqref{E1}. Recalling the definition of $E_1$ in \eqref{def:E1}, and  integrating by parts,
 \begin{align*}
\Phi_\eps^2(0,0)\int_{B_\delta}&k_sE_1(y,0,\Phi_\eps,\nabla \Phi_\eps, \delta,s)\,dy
=\phi_\eps^2(0)k_s\left[\frac{N-2s}{2}\int_{B_\delta}\phi_\eps^{p_\eps+1}\,dy+\int_{B_\delta}\langle y,\nabla\phi_\eps\rangle\phi_\eps^{p_\eps}\,dy\right]\\
&=\phi_\eps^2(0)k_s\left[\left(\frac{N-2s}{2} -\frac{N}{p_\eps+1}\right)
\int_{B_\delta}\phi_\eps^{p_\eps+1}\,dy
 +\frac{\delta}{p_\eps+1}\int_{\partial B_\delta}\phi_\eps^{p_\eps+1}\,d\sigma\right]\leq C\phi^{1-p_\eps}_\eps(0)=o_\eps(1),
 \end{align*}
 where the inequality follows from \eqref{def_phi_eps} and Proposition~\ref{prop:Li4.9}, while the conclusion is a consequence of \eqref{phi_eps_0_to_0}.
\end{proof}

\subsection{Pohozaev-type identities}
From  \cite{Ros-OtonSerra2014} we know that, for  $u\in \mathcal H^s_0(\Omega)$, 
\begin{align*}
\text{if $(-\Delta)^{s}u\in L^{\infty}(\Omega)$ in the case $2s>1$\qquad  or\qquad  $(-\Delta)^{s}u\in C^{\alpha}_{loc}(\Omega)\cap L^{\infty}(\Omega)$  in the case $2s\leq 1$ (with $\alpha>1-2s$),}
\end{align*}
then the following fractional Pohozaev identity holds 
\begin{align}\label{GeneralPohozaev}
\int_{\Omega}\langle \,\cdot-y_0,\nabla u\rangle (-\Delta)^{s}u\,dx=-\frac{N-2s}{2}\int_{\Omega}u(-\Delta)^{s}u\,dx -\frac{\Gamma^{2}(1+s)}{2}\int_{\partial\Omega}\left(\frac{u}{d^s}\right)^2\langle \,\cdot-y_0,\nu\rangle \,d\sigma
\end{align}
for any\footnote{Actually, in \cite{Ros-OtonSerra2014,RS14}, this is stated in the case $y_0=0$, but the same proof can be done for any  $y_0\in\mathbb R^N$ (see also \cite{DjitteFallWeth2024}).} $y_0\in\mathbb R^N$.

If $u$ is  a solution  of \eqref{equn}, then by easy calculations we have that $(p+1)\int_{\Omega}\langle \,\cdot-y_0,\nabla u\rangle u^p\,dx=-N\int_{\Omega}u^{p+1}\,dx$ and 
\eqref{GeneralPohozaev} becomes
\begin{equation}\label{Pohozaev}\left[\frac{2N}{p+1}-(N-2s)\right]\int_{\Omega}u^{p+1}\,dx=\Gamma^2(1+s)\int_{\partial\Omega}\left(\frac{u}{d^s}\right)^2\langle \,\cdot\ {-y_0},\nu \rangle d\sigma\qquad \text{for any $y_0\in\mathbb R^N$.}
\end{equation}

We also need the following bilinear version of the fractional Pohozaev identity which is an immediate consequence of \eqref{GeneralPohozaev} (with $y_0=0$)
\begin{align}\label{bilinearPohozaev}
\int_{\Omega}\langle \,\cdot\, ,\nabla u\rangle (-\Delta)^{s}\varphi\,dx+\int_{\Omega}\langle\,\cdot\, ,\nabla \varphi\rangle (-\Delta)^{s}u\,dx=\frac{2s-N}{2}\int_{\Omega}\left[u(-\Delta)^{s}\varphi+\varphi(-\Delta)^{s}u\right]\,dx -\Gamma^{2}(1+s)\int_{\partial\Omega}\frac{u}{d^s}\frac{\varphi}{d^s}\langle \,\cdot\, ,\nu\rangle \,d\sigma.
\end{align}
This identity holds true for  $u,\varphi\in \mathcal H^s_0(\Omega)$ such that 
\begin{align*}
    \text{$(-\Delta)^{s}u, (-\Delta)^{s}\varphi\in L^{\infty}(\Omega)$ if $2s>1$
    \qquad or \qquad $(-\Delta)^{s}u, (-\Delta)^{s}\varphi\in C^{\alpha}_{loc}(\Omega)\cap L^{\infty}(\Omega)$ if $2s\leq 1$ (with $\alpha>1-2s$).
}
\end{align*}
A special case of \eqref{bilinearPohozaev} is given in the following result, where $u$ and $\varphi$ are chosen to be a solution of \eqref{equn} and of the linearized problem at $u$, respectively.

\begin{lemma}\label{Lemma:boundaryTerm=0} Let $\Omega\subset \R^N$, $N\geq2 $, be an open smooth bounded set, and let $u$ be a solution to \eqref{equn} and $\varphi$ be a solution to the linearized problem at $u$, namely 
\begin{equation}\label{LinearizedProblem}
\left\{
\begin{array}{ll}
(-\Delta)^s \varphi -pu^{p-1}\varphi =0& \mbox{in }\Omega,\\
\varphi=0 & \mbox{in }\mathbb R^N\setminus\Omega.
\end{array}
\right.\end{equation}
 Then,
\begin{equation}\label{boundaryTerm=0}
\int_{\partial\Omega}\frac{u}{d^s}\frac{\varphi}{d^s}\langle x,\nu\rangle \,d\sigma=0.
\end{equation} 
\end{lemma} 

\begin{proof} Taking $u$ as a test function in \eqref{LinearizedProblem}, $\varphi$ as a test function in \eqref{equn}, and subtracting both equations, we get
\begin{equation}\label{upv}
(1-p)\int_{\Omega}u^p\varphi=0.
\end{equation}
As a consequence of \eqref{upv}, taking $u$ as a test function in \eqref{LinearizedProblem},  $\varphi$ as a test function in \eqref{equn}, and adding  both equations, we get
\begin{equation}\label{upvBIS}\int_{\Omega}\left[u(-\Delta)^{s}\varphi+\varphi(-\Delta)^{s}u\right]\,dx =\int_{\Omega}\left[u(pu^{p-1}\varphi)+\varphi u^{p}\right]\,dx=(p+1)\int_{\Omega}u^{p}\varphi\,dx=0.
\end{equation}
Let us now denote
\begin{align}\label{defIuvarphi}
\mathcal{I}(u,\,\varphi):=\int_{\Omega}\langle x,\nabla u\rangle (-\Delta)^{s}\varphi\,dx+\int_{\Omega}\langle x,\nabla \varphi\rangle (-\Delta)^{s}u\,dx.
\end{align}
By the bilinear version  \eqref{bilinearPohozaev} of the fractional Pohozaev identity, we have
\begin{align}\label{uv}
\mathcal{I}(u,\,\varphi)=\frac{2s-N}{2}\int_{\Omega}\left[u(-\Delta)^{s}\varphi+\varphi(-\Delta)^{s}u\right]\,dx -\Gamma^{2}(1+s)\int_{\partial\Omega}\frac{u}{d^s}\frac{\varphi}{d^s}\langle x,\nu \rangle\,d\sigma,
\end{align}
hence, \eqref{upvBIS} implies that
\begin{align}\label{uvFINAL}
\mathcal{I}(u,\,\varphi)=
-\Gamma^{2}(1+s)\int_{\partial\Omega}\frac{u}{d^s}\frac{\varphi}{d^s}\langle x, \nu \rangle \,d\sigma.
\end{align}
Let us now compute the quantity $\mathcal I(u,\varphi)$ defined in \eqref{defIuvarphi}, using the equations \eqref{equn} and \eqref{LinearizedProblem} for $u$ and $\varphi$ respectively.  After an  integration by parts, we get
\begin{eqnarray}\label{calcoloI}
\mathcal{I}(u,\,\varphi)&=& \int_{\Omega}\langle x,\nabla u\rangle pu^{p-1}\varphi\,dx+\int_{\Omega}\langle x,\nabla \varphi\rangle  u^p\,dx=\int_\Omega \langle x, \nabla(u^p)\rangle \varphi\,dx+\int_{\Omega}\langle x,\nabla \varphi\rangle u^p\,dx  \nonumber\\
&=&-\sum_{i=1}^N\int_\Omega u^p\frac{\partial (x_i\varphi)}{\partial x_i}\,dx+\int_{\Omega}\langle x,\nabla \varphi\rangle u^p\,dx 
=-N\int_\Omega u^p\varphi\,dx -\int_{\Omega}\langle x,\nabla \varphi\rangle u^p\,dx+\int_{\Omega}\langle x,\nabla \varphi\rangle u^p\,dx  \nonumber\\
&=&-N\int_\Omega u^p\varphi\,dx =0, 
\end{eqnarray}
where the last equality follows from \eqref{upv}. Hence, \eqref{boundaryTerm=0} follows substituting \eqref{calcoloI} into  the left-hand side of \eqref{uvFINAL}:
\[
0= \Gamma^{2}(1+s)\int_{\partial\Omega}\frac{u}{d^s}\frac{\varphi}{d^s}\langle x,\nu\rangle \,d\sigma.
\]
\end{proof}

\subsection{A Kelvin-type transform}\label{subsec:Kelvin}
First, let us recall the following properties of the point inversion transformation, see, for instance, \cite{FallWeth_2012_JFA,ADFJS19,b16}. Let $u:\rn\to \r$ be a measurable function, $x_0 \in \mathbb{R}^N$, $c>0$, and, for $x\in \mathbb{R}^N \setminus \{-x_0\}$, let 
\[
\kappa x := c \, \frac{x+x_0}{|x+x_0|^2} - x_0
\quad \text{and} \quad
K_s u(x) := |x+x_0|^{2s-N} \, u(\kappa x)
\quad \text{for } x \in \mathbb{R}^N \setminus \{-x_0\}.
\]
Then, for $x \in \mathbb{R}^N \setminus \{-x_0\}$,
\begin{align}\label{general:Kelvin}
(-\Delta)^s (K_s u)(x)
= c^{2s} \frac{K_s\big((-\Delta)^s u\big)(x)}{|x+x_0|^{4s}},\qquad \text{ i.e., }\qquad 
(-\Delta)^s \left( \frac{u \circ \kappa(x)}{|x+x_0|^{N-2s}} \right)
= c^{2s} \frac{(-\Delta)^s u(\kappa x)}{|x+x_0|^{N+2s}},
\end{align}
see, for example, \cite[Proposition 2]{ADFJS19}. In particular, if $x_0=0$ and $c=1,$ we have the Kelvin function $\hat u$ of $u$ given by 
\[
\hat u(x):=|x|^{2s-N}u\left(\frac{x}{|x|^2}\right)\qquad \text{for a.e. $x\in\mathbb R^N\setminus\{0\}$}.
\]

Let $D\subset\mathbb R^N$ be an arbitrary open set and define
\[
\widehat D
:=
\left\{
\frac{x}{|x|^2}:x\in D\setminus\{0\}
\right\}
\subset\mathbb R^N\setminus\{0\}.
\]
If
$u\in\mathcal D^s(\mathbb R^N)$, then
\[
\widehat u\in\mathcal D^s(\mathbb R^N),
\qquad
\|\widehat u\|_{\mathcal D^s}
=
\|u\|_{\mathcal D^s},
\qquad
\widehat{\widehat u}=u
\quad\text{a.e. in }\mathbb R^N.
\]
Moreover, if $u\in\mathcal D_0^s(D)$, then $\widehat u\in\mathcal D_0^s(\widehat D);$ see, for instance,
\cite[Lemma~2.2 and Corollary~2.3]{FallWeth_2012_JFA}.
In addition, a change of variables gives 
\begin{equation}\label{normInvarianceKelvin}
\|u\|_{2_s^*}=\|\widehat u\|_{2_s^*}.
\end{equation}

Let now $f\in L^1_{\mathrm{loc}}(D)$ and suppose that
$u\in\mathcal D_0^s(D)$ satisfies $(-\Delta)^s u=f$ in $D$ in the weak sense. Then $\widehat u$ satisfies 
\begin{equation}\label{sharm}
\left\{
\begin{array}{ll}
\displaystyle
(-\Delta)^s\widehat u(x)
=
|x|^{-N-2s}
f\left(\dfrac{x}{|x|^2}\right)
=
|x|^{-4s}\widehat f(x)
& \text{in }\widehat D,\\[3mm]
\widehat u=0
& \text{in }\mathbb R^N\setminus\widehat D,
\end{array}
\right.
\end{equation}
where
\[
\widehat f(x)
:=
|x|^{2s-N}
f\left(\frac{x}{|x|^2}\right),
\qquad x\in\widehat D.
\]

\subsection{Fractional Green's function}

We denote by $\mathcal G$ the Green's function related to the 
fractional Laplacian in $\Omega$, namely, for every fixed $y\in\Omega$, 
$ \mathcal G(\cdot,y)$ is the unique solution of
\[\left\{\begin{array}{ll}
(-\Delta)^s \mathcal G(\cdot, y)=\delta_{y} &\mbox{ in } \Omega,\\
\hspace{1cm}\mathcal G(\cdot, y) = 0 &\mbox{ in } \mathbb R^N\setminus\Omega,
\end{array}\right.
\]
where $\delta_{y}$ denotes the Dirac mass at the point $y$. 
It is well known that the integral representation formula
    \[
w(x)=\int_{\Omega}\mathcal G(x,y)f(y)\,dy
    \]
holds true for solutions of
\[\left\{\begin{array}{ll}
(-\Delta)^s w =f  &\mbox{ in } \Omega,\\
\hspace{1cm}w = 0 &\mbox{ in } \mathbb R^N\setminus\Omega,
\end{array}\right.
\]
see \cite[Proposition 1.2.2]{N15}.

By \cite[Theorem 1.1]{CS98}, there is $C=C(\Omega,s)>0$ such that 
\begin{equation}\label{Green0}
\mathcal G(x,y)\leq C\frac{1}{|x-y|^{N-2s}}\qquad \text{ for all } x,y\in\Omega,\ x\neq y,
\end{equation}
and
\begin{equation}\label{Green1}
\mathcal G(x,y)\leq C\frac{d^s(x)}{|x-y|^{N-s}}\qquad\text{ for all } x,y\in\Omega,\ x\neq y.
\end{equation}
As a consequence, given $x_0\in\Omega$ and $\Omega'$ a neighborhood of $\partial\Omega$ such that $x_0\notin\overline{\Omega'}$, there exists $K=K(x_0,\Omega',s)>0$ such that
\begin{equation}\label{Green2}
0<\frac{\mathcal G(x,x_0)}{d^s(x)}\leq K \qquad \text{ for all } x\in\Omega'.
\end{equation}

Let us denote the regular part of the Green's function in $\Omega$ by
\[H(x,y):=\Phi(x,y)-\mathcal G(x,y),\] 
where
$\Phi(x,y)$ is the Green's function of $(-\Delta)^s$ in the whole space $\mathbb R^N$, which is explicit
\[\Phi(x,y):=\frac{a_{N,s}}{|x-y|^{N-2s}},\]
with $a_{N,s}=\frac{\Gamma\big(\frac{N}{2}-s\big)}{4^s\pi^{N/2}\Gamma(s)}$. Recall that $N\geq 2$ and $s\in (0,1)$, so $N-2s>0$.
Observe  that, for any fixed $y\in\Omega$, $H(\cdot,y)$ satisfies
\[\left\{\begin{array}{ll}
(-\Delta)^s H(\cdot, y)= 0&\mbox{ in } \Omega,\\
\hspace{1cm}H(\cdot, y) = \Phi(\cdot, y) &\mbox{ in } \mathbb R^N\setminus\Omega.
\end{array}\right.
\]
   
We define the fractional Robin function as
\begin{equation}\label{def:Robin}
\mathcal R(x):=H(x,x), \qquad x\in \Omega.
\end{equation}
By the fractional weak maximum principle (see, for instance, \cite{FallJarohs2015}) it follows that
\begin{equation}\label{RobinStrictlyPositive}\mathcal R(x) >0 \quad \text{ for all } x\in \Omega.\end{equation}
Let us recall the following  Pohozaev-type  identities involving the fractional  Green's and Robin functions recently proved in \cite{DjitteSueur,DiebIanni_poho}.  Let $\nu$ denote the outward unit normal to $\partial\Omega$ and recall that $d:=\operatorname{dist} (\ \cdot\ , \partial\Omega)$.
\begin{lemma}[{\cite[Theorem 1.1]{DiebIanni_poho}}]
\label{lemma:PohozaevGreen_DiebIanni}Let $ N\geq 2$, $s\in (0,1)$, 
and let $\Omega$ be a bounded open set of $\mathbb R^N$ of class $C^{1,1}$. 
Then, 
\begin{equation*}
\mathcal R(x_0)=\frac{\Gamma^2(1+s)}{N-2s}\int_{\partial\Omega}\left( \frac{\mathcal G(x_0, \cdot)}{d^s(\cdot)}\right)^2 \langle \cdot-x_0,\nu \rangle \, d\sigma\qquad \text{for every $x_0\in\Omega$}.
\end{equation*}
\end{lemma}
\begin{lemma}[{\cite[Corollary 1.7]{DjitteSueur}}]
\label{lemma:PohozaevGreen_DjitteSeur}Let $ N\geq 2$, $s\in (0,1)$, 
and let $\Omega$ be a bounded open set of $\mathbb R^N$ of class $C^{1,1}$. 
Then,
\begin{equation*}
\nabla\mathcal R(x_0)=\Gamma^2(1+s)\int_{\partial\Omega}\left( \frac{\mathcal G(x_0, \cdot)}{d^s(\cdot)}\right)^2 \nu  \, d\sigma\qquad \text{for every $x_0\in\Omega$}.
\end{equation*}
\end{lemma}

\section{Properties of positive solutions}

In this section we show some properties that hold for general families of solutions to \eqref{problema}.

\begin{lemma}\label{lemma:boundary estimate}
Let $\Omega$ be a smooth open bounded set and let $u_\eps$ be a family of solutions to \eqref{problema}. Then there exist $\delta>0$ and $C>0$ such that 
\begin{equation}\label{claim:interiorBlow-up}
\sup_{O(\delta)}u_\eps\leq C\qquad \text{ for all }\,\eps>0,
\end{equation}
where $O(\delta):=\{z\in\Omega\,:\,d(z,\partial\Omega)<\delta\}$.  As a consequence, if $x_0$ is a blow-up point for $u_\eps$, then $x_0\in\Omega$.
\end{lemma}
\begin{proof}
We first establish a priori bounds in a neighborhood of the boundary via the moving planes method. In the local setting this idea was introduced by de Figueiredo, Lions, and Nussbaum in \cite{dFLN}, and in particular  applied to the slightly subcritical standard Lane-Emden  problem in \cite[Section 2]{H91}. We follow a similar strategy, using the  nonlocal version of the moving plane introduced in \cite{JW14}.

Let $\lambda_1>0$ be the first eigenvalue of $(-\Delta)^s$ in $\Omega$ with exterior zero conditions and let $\phi_1\geq 0$ be the corresponding eigenfunction. Integrating by parts,
\begin{align*}
\lambda_1\int_{\Omega}\phi_1u_\eps\ dx=\int_{\Omega}(-\Delta)^s\phi_1u_\eps\ dx=\int_{\mathbb R^N}(-\Delta)^s\phi_1u_\eps\ dx=\int_{\mathbb R^N}\phi_1(-\Delta)^su_\eps\ dx=\int_{\Omega}\phi_1(-\Delta)^su_\eps\ dx=\int_{\Omega}\phi_1u_\eps^{p_\eps}\ dx.
\end{align*}
Let $\vartheta>0$ be a constant such that $\vartheta^{p_\eps-1}\geq 2\lambda_1$, then
\begin{align*}
\lambda_1\int_{\Omega}\phi_1u_\eps\ dx&=\int_{\Omega}\phi_1u_\eps^{p_\eps}\ dx
=\int_{\{u_\eps\geq\vartheta\}}\phi_1u_\eps^{p_\eps}\ dx+\int_{\{u_\eps<\vartheta\}}\phi_1u_\eps^{p_\eps}\ dx
\geq \int_{\{u_\eps\geq\vartheta\}}\phi_1u_\eps^{p_\eps}\ dx
\geq
\vartheta^{p_\eps-1}\int_{\{u_\eps\geq\vartheta\}}\phi_1u_\eps\ dx\\
&=\vartheta^{p_\eps-1}\int_{\Omega}\phi_1u_\eps\ dx
-\vartheta^{p_\eps-1}\int_{\{u_\eps<\vartheta\}}\phi_1u_\eps\ dx
\geq 2\lambda_1 \int_{\Omega}\phi_1u_\eps\ dx-C_1,
\end{align*}
where we used that $\vartheta^{p_\eps-1}\int_{\{u_\eps<\vartheta\}}\phi_1u_\eps\, dx
\leq \vartheta^{p_\eps}\int_{\Omega}\phi_1\, dx\leq C\int_{\Omega}\phi_1\, dx=:C_1$.
Therefore, 
\begin{align*}
\int_{\Omega}\phi_1u_\eps\ dx\leq \frac{C_1}{\lambda_1}.
\end{align*}
In any $\Omega'\subset\subset\Omega$, since $\phi_1\geq \widetilde C(\Omega')>0$, we obtain
\begin{equation}\label{localL1Bound}\int_{\Omega'}u_\eps\ dx\leq  C(\Omega')\qquad \text{ for all } \Omega'\subset\subset\Omega.
\end{equation}

If $\Omega$ is strictly convex, by applying the moving planes method in the fractional setting (see \cite{JW14}), it follows that each solution $u_\eps$ is non-decreasing along an arbitrary straight line emanating from a point on $\partial\Omega$ towards the inside of $\Omega$. Hence, we can find $\gamma,\delta>0$ such that for any $x\in O(\delta)$, there exists a measurable set $\Gamma_x$ such that
\begin{itemize}
\item[(i)] $meas (\Gamma_x)\geq \gamma$,
\item[(ii)] $\Gamma_x\subset \{z\in\Omega \ :\ dist (z,\partial\Omega)>\frac{\delta}{2}\}$,
\item[(iii)] $u_\eps(y)\geq u_\eps(x)$
for any $y\in\Gamma_x$.
\end{itemize}
In particular, $\Gamma_x$ can be taken as a piece of a cone with vertex at $x$. Let us choose $\Omega':=\{z\in\Omega \ :\ dist (z,\partial\Omega)>\frac{\delta}{2}\}\ (\subset\subset\Omega)$, then \eqref{localL1Bound} implies that, for $x\in O(\delta),$
\begin{align*}
u_\eps(x)\overset{(iii)}{\leq} \frac{1}{meas (\Gamma_x)}\int_{\Gamma_x}u_\eps(y)\ dy\overset{(ii)}{\leq} \frac{1}{meas (\Gamma_x)}\int_{\Omega'}u_\eps(y)\ dy\overset{(i)}{\leq} \frac{1}{\gamma}\int_{\Omega'}u_\eps(y)\ dy\overset{\eqref{localL1Bound}}{\leq}\frac{C(\Omega')}{\gamma},
\end{align*}
which proves \eqref{claim:interiorBlow-up} when $\Omega$ is strictly convex. For a general smooth domain $\Omega$, not necessarily convex, following the approach in \cite{dFLN}, we  apply the moving planes method on the Kelvin transform of $u_\eps$ near each boundary point in order to avoid the difficulty of an empty cap. 
Indeed, in this situation one can use the nonlocal Kelvin transform (see Section~\ref{subsec:Kelvin}) and \eqref{general:Kelvin} for some $v\in \rn$ and some $c>0$, to obtain a transformed problem where the nonlocal moving plane method in \cite{JW14} can be applied locally. We thus get the monotonicity in a neighborhood of the boundary point.  The conclusion then follows from an interior integral estimate, similarly as in the previous case. We omit the details, since they are the same as in the local case, and refer the reader to \cite{dFLN, H91} for them.
\end{proof}

\begin{lemma}\label{blowup:positive:lem}
Let $\Omega\subset \rn$ be a smooth bounded starshaped domain and let $u_\eps$ be a family of solutions to \eqref{problema}. Then, 
\begin{align}\label{blowupueps}
\|u_\eps\|_\infty\to \infty\qquad \text{ as $\eps\to 0$.}     
\end{align}
\end{lemma}
\begin{proof}
We show first \eqref{blowupueps}. By contradiction, assume that 
 \begin{align}\label{cont:bd:1}
 \|u_\eps\|_\infty \le C\qquad \text{ for some $C>0$.}
 \end{align}
 By Hölder's inequality, $\|u_\eps\|^2_{p_\eps+1}\leq
|\Omega|^\frac{2(2^*_s-p_\eps-1)}{(p_\eps+1)2^*_s} \|u_\eps\|^2_{2^*_s}$. Then, using that $\|u_\eps\|_{\mathcal D^s}^2=\|u_\eps\|_{p_\eps+1}^{p_\eps+1}$, 
 \begin{align}\label{cont:bd2}
S_{N,s} \le 
\frac{\|u_\eps\|_{\mathcal D^s}^2}
{
\|u_\eps\|^{2}_{2^*_s}
}
\le |\Omega|^\frac{2(2^*_s-p_\eps-1)}{(p_\eps+1)2^*_s}
\frac{\|u_\eps\|_{\mathcal D^s}^2}
{\|u_\eps\|^2_{p_\eps+1}}
= |\Omega|^\frac{2(2^*_s-p_\eps-1)}{(p_\eps+1)2^*_s} 
\|u_\eps\|^{p_\eps-1}_{p_\eps+1}
\le \|u_\eps\|_\infty^{p_\eps-1} |\Omega|^{\frac{2(2^*_s-p_\eps-1)}{(p_\eps+1)2^*_s}+\frac{p_\eps-1}{p_\eps+1}}.
 \end{align}

In the following, we use $C$ to denote possibly different positive constants independent of $\eps$.  By \eqref{cont:bd:1} and \eqref{cont:bd2}, passing to a subsequence, we have that $\|u_\eps\|_\infty \to c_0$ for some $c_0>0$. Since $u_\eps$ is uniformly bounded in $\Omega$, one can use regularity theory \cite{RS14} (see also \cite[Lemma 2.2]{DIS23}) to deduce that $\|u_{\eps_n}\|_{C^s(\overline{\Omega})}<C$. By interior regularity, we then have that $\|u_{\eps_n}\|_{C^{2s+\alpha}(O)}<C(O)$ for every open set $O\subset \subset \Omega$ and some $\alpha\in(0,1)$, see \cite{RS14} or \cite[Lemma 2.4]{DIS23}.  Furthermore, using higher-order Hölder regularity \cite[Theorem 1.4]{fall2020regularity}, one can have that $\|u_{\eps_n}\|_{C^{2+\alpha}(O)}<C(O)$ for every open set $O\subset \subset \Omega$ (in fact, iterating these estimates, one has that $u_\eps$ is smooth in $\Omega$). Now, using the compact embedding $C^{2+\alpha}(O) \hookrightarrow C^{2}(O),$ there is $u_0$ such that, up to a subsequence, 
\begin{align}\label{reg1}
u_\eps \rightarrow u_0\qquad \text{ in $C_{loc}^{2}(\Omega)\cap C_{loc}^{s-\delta}(\rn)$ for any $\delta>0$.}    
\end{align}
 Using \eqref{cont:bd:1} and \eqref{reg1}, for every $x\in \Omega$, we have that $(-\Delta)^s u_\eps(x) \to (-\Delta)^s u_0(x)$ as $\eps\to 0$.  As a consequence, $u_0$ is a nonnegative nontrivial classical solution of
\begin{align}\label{critprob}
(-\Delta)^s u_0 = u_0^{2^*_s-1} \quad  \text{in } \Omega, \qquad 
u_0 = 0 \quad \text{ in } \mathbb R^N\setminus\Omega,\qquad 
u_0\geq 0\quad \text{ in }\rn.
\end{align}
Since $\Omega$ is starshaped with respect to the origin, the fractional Pohozaev identity implies that \eqref{critprob} only has the trivial solution (see \cite[Corollary 1.3]{Ros-OtonSerra2014}), and we reach a contradiction and \eqref{blowupueps} follows. 
\end{proof}

\section{Asymptotic analysis for low-energy solutions in general sets}\label{section:asymptotic}

Let $s \in(0,1)$ and let $\Omega$ be a smooth bounded open set in $\mathbb R^N$, $N\geq 2$, $p_\eps=2^{\ast}_s-1-\eps=\frac{N+2s}{N-2s}-\eps$, $\eps>0$, and let us consider the minimization problem
\begin{equation}\label{SNseps}
S_{N,s,\eps}
:=\inf_{w\in \mathcal H_0^s(\Omega)\setminus\{0\}}\frac{\|w\|_{\mathcal D^{s}}^2}{\|w\|_{2^*_s-\eps}^2}
=
\inf_{\substack{
w\in \mathcal H_0^s(\Omega)\\
\|w\|_{2_s^*-\varepsilon}=1
}}
\|w\|_{\mathcal D^s}^2.
\end{equation}
By the compactness of the embedding $\mathcal H_0^s(\Omega)\hookrightarrow L^{2^*_s-\eps}(\Omega)$,  for any $\eps\in(0,2^*_s-2)$, it follows that  \eqref{SNseps} admits a nontrivial minimizer $\tilde u_\eps\in \mathcal H_0^s(\Omega)$ such that $\|\tilde u_\eps\|_{2^*_s-\eps}=1$. As a consequence, the dilation $u_\eps:=(S_{N,s,\eps})^{\frac{1}{2^*_s-2-\eps}}\tilde u_\eps$ is still a minimizer of \eqref{SNseps} and solves
\begin{equation}\label{problema}
\left\{
\begin{array}{ll}
(-\Delta)^s u=u^{2^*_s-1-\eps }& \mbox{in }\Omega,\\
\hspace{1cm}u>0 & \mbox{in }\Omega,\\
\hspace{1cm}u=0 & \mbox{in }\mathbb R^N\setminus\Omega.
\end{array}
\right.\end{equation}
We refer to solutions to \eqref{problema} which are also minimizers of \eqref{SNseps} as \emph{least energy solutions} of \eqref{problema}. Least energy solutions can be also obtained by minimizing the associated energy functional on the Nehari manifold.   
\begin{remark}\label{rem:least_are_low}
In \cite[Theorem 1.6]{PP2014} it has been proved that
\begin{equation}\label{SnsepstoSNs}
S_{N,s,\eps}\to S_{N,s}\quad\text{ as $\eps\to0$}.
\end{equation}
\eqref{SnsepstoSNs} implies that any family $u_\eps$ of least energy solutions to \eqref{problema} satisfies
\begin{equation}\label{limSn}
\frac{\|u_\eps\|_{\mathcal D^s}^2}{\left(\int_{\Omega}\left|u_\eps\right|^{2^*_s-\eps}\right)^{\frac{2}{2^*_s-\eps}}} \to  S_{N,s}\qquad \text{ as }\eps\to 0.
\end{equation}  
\end{remark}

\begin{definition}\label{def:low_energy_solutions}
We say that a family of solutions to \eqref{problema} satisfying \eqref{limSn} is a family of {low energy solutions}.
\end{definition}

Since any solution $u_\eps$ of \eqref{problema} verifies
\[
\|u_\eps\|_{\mathcal D^s}^2=\int_{\Omega}u_\eps^{p_\eps+1}\,dx,
\]
by \eqref{limSn} one has that, for any family of low energy solutions,
\begin{equation}\label{H_s_bound}\lim_{\eps\rightarrow 0}\int_{\Omega}u_\eps^{p_\eps+1}\,dx=\lim_{\eps\rightarrow 0}\|u_\eps\|_{\mathcal D^s}^2 =\left(S_{N,s}\right)^{\frac{N}{2s}}.\end{equation}
Next we show that low-energy solutions blow-up as $\eps\rightarrow 0$.
\begin{lemma}\label{lemma:blow-up}
    Let $\Omega$ be a smooth  bounded open set and let $u_\eps$ be a family of low energy solutions of \eqref{problema}. 
Let $x_\eps\in\Omega$ be such that
\begin{align}\label{xepsdef}
u_\eps(x_\eps)=\|u_\eps\|_{\infty}.  
\end{align}    
 Then,  up to a subsequence,  $x_\eps$ converges, and let 
\begin{equation}\label{x0}
x_0:=\lim_{\eps\rightarrow 0} x_\eps.
\end{equation}   
Then $x_0$ is a blow-up point with blowing-up sequence $x_\eps$, namely
\begin{equation}\label{blowupoccursLOW}
\|u_\eps\|_{\infty}\rightarrow +\infty\quad\mbox{ as }\eps\rightarrow 0.
\end{equation}
Furthermore, $x_0$ is an interior point of $\Omega$.
\end{lemma}
\begin{proof}
In order to get \eqref{blowupoccursLOW} we adapt some ideas in \cite{H91}. Let us assume by contradiction that 
\begin{equation}\label{ContradictionAssumption:NoBlowup}
\|u_{\eps_n}\|_{\infty}\leq C\end{equation} for a sequence $\eps_n\rightarrow 0$ as $n\rightarrow +\infty$.
By \eqref{H_s_bound}, $u_{\eps_n}$ is also bounded in $\mathcal H^s_0(\Omega)$. Hence, up to a subsequence, $u_{\eps_n}$ converges to a function $v$ weakly in $\mathcal H^s_0(\Omega)$, strongly in $L^2(\Omega)$, and pointwise a.e. in $\Omega$, with
\begin{equation}\label{wlsemicontinuity}
\|v\|_{\mathcal D^s}^2\leq \liminf_{n\to\infty} \|u_{\eps_n}\|_{\mathcal D^s}^2.\end{equation}
Moreover, by \eqref{ContradictionAssumption:NoBlowup} and dominated convergence,
\begin{equation}\label{convergLp}\int_{\Omega}|u_{\eps_n}|^{p_{\eps_n}+1}\,dx\rightarrow \int_{\Omega}|v|^{2^*_s}\,dx\quad\mbox{ as }n\rightarrow\infty.\end{equation}
Hence, by \eqref{H_s_bound}, $\int_{\Omega}|v|^{2^*_s}\,dx=\left(S_{N,s}\right)^{\frac{N}{2s}}$; in particular, $v\not\equiv 0$. Furthermore, by \eqref{limSn},
\[ \|u_{\eps_n}\|_{\mathcal D^s}^2=\left(S_{N,s}+o(1)\right)\left(\int_{\Omega}\left|u_{\eps_n}\right|^{p_{\eps_n}+1}\right)^{\frac{2}{p_{\eps_n}+1}} \qquad \text{ as }n\to +\infty.\]
Passing to the limit as $n\rightarrow \infty$ and using \eqref{wlsemicontinuity}, \eqref{convergLp} and \eqref{SNS}
\[ \frac{\|v\|_{\mathcal D^s}^2}{\|v\|_{2^*_s}^{2} }= S_{N,s},\]
a contradiction with the fact that the best Sobolev constant $S_{N,s}$ is never achieved on a bounded open set (see Section~\ref{section:SobolevConstant_and_bubbles}). Thus \eqref{blowupoccursLOW} is proved.

As a consequence of \eqref{blowupoccursLOW} and of the general uniform boundary bound in Lemma~\ref{lemma:boundary estimate} we also deduce that $x_0\in \Omega$.
%
 %
%
\end{proof}

\begin{remark}
Under the additional assumption that $\Omega$ is star-shaped, one can prove that \emph{any} family of solutions to \eqref{problema} blows-up, namely that the low energy assumption in Lemma~\ref{lemma:blow-up} can be removed, see Lemma~\ref{blowup:positive:lem}. 
\end{remark}

In order to characterize the blow-up behavior of $u_\eps$, we scale it around the blow-up sequence $x_\eps$ defined in \eqref{xepsdef}, defining
\begin{equation}\label{defvepsilon}
v_\eps(x):=
\frac{1}{\|u_\eps\|_{\infty}}u_\eps \left( x\|u_\eps\|_{\infty}^{-\frac{p_\eps-1}{2s}} +x_\eps\right)
,\qquad x\in \rn,
\end{equation}
\begin{equation}\label{gamma:def}
{
\|u_\eps\|_{\infty}^{-\frac{p_\eps-1}{2s}}}\to 0\qquad \text{as }\eps\to 0,
\end{equation}
(where the convergence follows from Lemma~\ref{lemma:blow-up}). 

Then $\|v_\eps\|_{\infty}=v_\eps(0)=1$ and $v_\eps$ satisfies
\begin{equation}
    \label{eq:v_eps}
\left\{
\begin{array}{ll}
(-\Delta)^s v_\eps=v_\eps^{p_\eps}& \mbox{in }\Omega_\eps,\\
\hspace{1cm}v_\eps>0 & \mbox{in }\Omega_\eps,\\
\hspace{1cm}v_\eps=0 & \mbox{in }\mathbb R^N\setminus\Omega_\eps,
\end{array}
\right.\end{equation}
where $\Omega_\eps:= \|u_\eps\|_{\infty}^{\frac{p_\eps-1}{2s}}\left(\Omega-x_\eps\right)
$.

\begin{lemma}\label{lemma:wetoW}
Let $u_\eps$ be a family of low energy solutions to \eqref{problema}, let $v_\eps$ be as in \eqref{defvepsilon} and let $U$ be as in \eqref{Bubble}. Then,
    \begin{align}\label{convv}
    v_\eps\rightarrow U\qquad \text{  in  $C^{2}_{loc}(\mathbb R^N)$  as } \eps\rightarrow 0.
    \end{align}   
\end{lemma}
\begin{proof}
Let us observe that, by Lemma~\ref{lemma:blow-up},
\[
\|u_\eps\|_{\infty}^\eps\geq 1.
\]
Hence, using the definition of $v_\eps$ in \eqref{defvepsilon} and \eqref{H_s_bound} we get a uniform $\mathcal D^{s}(\mathbb R^N)$ bound, namely,
\begin{equation}\label{bounds_w_eps}\|v_\eps\|^2_{\mathcal D^s}=\int_{\Omega_\eps}v_\eps^{p_\eps+1}
=
\|u_\eps\|_{\infty}^{-\eps\frac{N-2s}{2s}}
\int_{\Omega}u_\eps^{p_\eps+1}
\leq \int_{\Omega}u_\eps^{p_\eps+1}
= \left(S_{N,s}\right)^{\frac{N}{2s}} +o(1).\end{equation}
As a consequence, up to a subsequence, 
\begin{align}\label{vweakto}
v_\eps\rightharpoonup v\qquad \text{ weakly in $\mathcal D^{s}(\mathbb R^N)$. }    
\end{align}
Furthermore, by Lemma~\ref{lemma:blow-up}, it also follows that   $\Omega_\eps\rightarrow \mathbb R^N$ (which means that for any bounded $K\subset \mathbb R^N$, we have $\Omega_\eps\supset K$ for $\eps$ small). Let $\varphi\in C^{\infty}_c(\mathbb R^N)$, then $\varphi\in C^{\infty}_c(\Omega_\eps)$ for $\eps$ small enough. Hence, by \eqref{eq:v_eps} and the definition of weak solution
\begin{align}\label{weakform}
[v_\eps,\varphi]_s=\int_{\Omega_\eps}v_\eps^{p_\eps}\varphi\,dx.   
\end{align}
Note that $|v_{\eps_n}^{p_\eps}\varphi|\leq |\varphi|$, since $|v_\eps|\leq 1$, by definition. Using \eqref{vweakto}, \eqref{weakform}, and dominated convergence,
\[
[v,\varphi]_s=\int_{\mathbb R^N}v^{\frac{N+2s}{N-2s}}\varphi\,dx \qquad \text{ for all } \varphi \in C^{\infty}_c(\mathbb R^N).
\]
Moreover, since $\|v_\eps\|_{\infty}=1$, by definition, one can use a priori regularity estimates \cite{RS14,fall2020regularity} (see the proof of Lemma~\ref{blowup:positive:lem}) to deduce that $v_{\eps_n}\rightarrow v$ in $C^{2}_{loc}(\mathbb R^N)$. As a consequence, $v\geq 0$ and,  moreover, as $\max_{\Omega_\eps} v_\eps=v_\eps(0)=1$, we have that $v(0)=1$ and $v\leq 1$. Hence,  $v$ is a nontrivial solution of the limit critical  equation \eqref{equationBubble}.  In particular, by \cite[Theorem 1]{CLO}, it follows that $v=U$, where $U$ is the function in \eqref{Bubble}.

\end{proof}

\begin{corollary}\label{cor:Corollary1_HAN}
Let $u_\eps$ be a family of low energy solutions to \eqref{problema}, then there exists $C>0$ such that
    \begin{equation}\label{2bounds}
    1\leq \|u_\eps\|_{\infty}^\eps\leq C \qquad \text{as $\eps\to0$}.
    \end{equation}
\end{corollary}

\begin{proof}[Proof of Corollary~\ref{cor:Corollary1_HAN}]
  The lower bound  in \eqref{2bounds} is a direct consequence of  Lemma~\ref{lemma:blow-up}, hence it is enough to prove the upper bound. Let $v_\eps$ be the rescaled function defined in \eqref{defvepsilon}.
    Lemma~\ref{lemma:wetoW} guarantees that, up to a subsequence, 
    \begin{align}\label{vconvU}
    \text{$v_\eps$ converges uniformly to $U$ in $B_1(0)$, }    
    \end{align}
    so on the one hand, there exists $C>0$, depending only on $U$, such that
    \begin{equation}\label{otoh}
        \int_{B_1(0)}v_\eps^{p_\eps+1}\geq C\qquad \text{as $\eps\to0$}.
    \end{equation}
    On the other hand,
    \begin{eqnarray*} \int_{B_1(0)}v_\eps^{p_\eps+1}&=&\frac{1}{\|u_\eps\|_\infty^{p_\eps+1}}\int_{B_1(0)}u_\eps^{p_\eps+1} \left( x \, \|u_\eps\|_{\infty}^{-\frac{p_\eps-1}{2s}} +x_\eps\right)dx\nonumber\\
    &=&\|u_\eps\|_\infty^{-\eps\frac{N-2s}{2s}}\int_{B_{r_\eps}(x_\eps)}u_\eps^{p_\eps+1}(y)dy\leq\|u_\eps\|_\infty^{-\eps\frac{N-2s}{2s}}\int_{\Omega}u_\eps^{p_\eps+1}(y)dy,
    \end{eqnarray*}
    where $r_\eps:=\|u_\eps\|_{\infty}^{-\frac{p_\eps-1}{2s}}$($\rightarrow 0$ by  \eqref{gamma:def}), hence 
\begin{equation}\label{otoh2}
\int_{B_1(0)}v_\eps^{p_\eps+1}\leq C\|u_\eps\|_\infty^{-\eps\frac{N-2s}{2s}} 
\end{equation}
 by \eqref{H_s_bound}.   
    Combining \eqref{otoh} and \eqref{otoh2}, we get that  $\|u_\eps\|_\infty^{\eps\frac{N-2s}{2s}}$ is bounded from above, by a constant independent of the subsequence, and hence the upper bound in \eqref{2bounds} follows.   
        \end{proof}

\begin{remark}\label{rem:finite energy}
If $u_\eps$ is a family of solutions to \eqref{problema} which is uniformly bounded in $\mathcal H^s_0(\Omega)$ and such that 
\[
\|u_\eps\|_{\infty}\to+\infty,
\]
then, as in Lemma~\ref{lemma:blow-up}, $x_0\in\Omega$ and also the conclusion in Corollary~\ref{cor:Corollary1_HAN} still holds.
\end{remark}

\begin{remark}\label{Remark:MissiongYanYanLi}
In the following, the assumption that $u_\eps$ is a family of low energy solutions is crucial in order to derive pointwise estimates. 
In general, for families of bounded energy nonnegative solutions it does not seem direct to apply a fractional version of Struwe's concentration-compactness result (see \cite{Struwe}, \cite{FangGonzalez}), because it is not clear whether such families of solutions are $(PS)$-sequences for the limiting functional. 

We emphasize that, in some papers in the literature, this has been assumed without a proper argument, see for example \cite[proof of Lemma 2.2]{choikim}, where a sequence of solutions with bounded energy is implicitly assumed to be a Palais-Smale sequence for the critical problem. However, this does not seem to be direct, even in the local case $s=1$.
\end{remark}

\begin{lemma}\label{lemma:Lemma3_HAN} 
Let $u_\eps$ be a family of low energy solutions to \eqref{problema} and let $v_\eps$ be given by \eqref{defvepsilon}.
There exists $\Lambda>0$ independent of $\eps$ such that
\begin{equation}\label{pointwiseBound}v_\eps(x)\leq \Lambda U(x)\qquad \mbox{ as $\eps\to0$ for all }x\in\Omega_\eps,\end{equation}
where $U$ is the function in \eqref{Bubble}.
\end{lemma}
\begin{remark}\label{rmk:bound}
    Using the definition of $v_\eps$ in \eqref{defvepsilon} one easily sees that \eqref{pointwiseBound} is equivalent to
    \begin{equation}\label{pointwiseBoundU_eps}
    u_\eps(y)\leq \Lambda \|u_\eps\|_{\infty}U\left((y-x_\eps)\|u_\eps\|_{\infty}^{\frac{p_\eps-1}{2s}} \right)\qquad \mbox{ for all }y\in\Omega.
    \end{equation}
\end{remark}
\begin{proof}[Proof of Lemma~\ref{lemma:Lemma3_HAN}]
Since $v_\eps\leq 1$ in $\rn$, the bound \eqref{pointwiseBound} trivially holds true locally. Indeed  for any fixed region around the origin there exists a constant $\Lambda>0$ such that $\Lambda U(x)>1$ for any $x$ in this region.

Next we show that  \eqref{pointwiseBound} also holds true near infinity. By the definition of $U$ (see  \eqref{Bubble}), this is equivalent to proving that
there exists a constant $C>0$, independent of $\eps$, such that 
\begin{equation}\label{LimnftyBoundKelvin}\hat v_\eps(x)\leq C \qquad \mbox{ for $x$ near the origin $O$,}\end{equation}
where 
\[\hat v_\eps(x):=\frac{1}{|x|^{N-2s}}v_\eps\left(\frac{x}{|x|^2}\right), \qquad  x\in\hat\Omega_\eps,\]
is the Kelvin transform of $v_\eps$ and $\hat\Omega_\eps$ is the image of $\Omega_\eps$ under the  Kelvin transform, which is the whole $\mathbb R^N$ except a small region near the origin $O$ (see Section~\ref{subsec:Kelvin} for the definition).

 By \eqref{sharm}, $\hat v_\eps$ solves 
\begin{align}\label{vepsP}
\left\{ 
\begin{array}{ll}
(-\Delta)^s \hat v_\eps=\frac{1}{|x|^{\eps(N-2s)}}{\hat v_\eps}^{p_\eps}& \mbox{in }\hat\Omega_\eps,\\
\hspace{1cm}\hat v_\eps>0 & \mbox{in }\hat\Omega_\eps,\\
\hspace{1cm}\hat v_\eps=0 & \mbox{in }\mathbb R^N\setminus\hat\Omega_\eps.
\end{array}\right.
\end{align}

We argue similarly as in \cite[Proof of Lemma 3]{H91}.  Let $b_\eps(x):=|x|^{-\eps(N-2s)}$. We show  that $\|b_\eps\|_{L^\infty(\hat\Omega_\eps)}$ is bounded independently of $\eps$, indeed for $x\in \hat\Omega_\eps$ there exists $y_{x,\eps}\in\Omega$ such that
\begin{equation}
\label{stima|x|}
|x|=\frac{\|u_\eps\|_{\infty}^{-\frac{p_\eps-1}{2s}}}{|y_{x,\eps}-x_\eps|},\end{equation}
and so, by Corollary~\ref{cor:Corollary1_HAN},
\begin{equation}\label{Uniform_bound_b}
b_\eps(x)=\|u_\eps\|_{\infty}^{\eps\frac{p_\eps-1}{2s}(N-2s)}|y_{x,\eps}-x_\eps|^{\eps (N-2s)}\leq C.\end{equation}

Let $\Theta_0,\,r_0>0$ be given by Lemma~\ref{lemma:Lemma_6_Han} applied to $\hat v_\eps$ in $\hat\Omega_\eps$ with $b=b_\eps$ and $q=p_\eps+1=2^*_s-\eps$, and let us observe that $\Theta_0$ and $r_0$ are independent of $\eps$.  By \eqref{ExactValue2}, there exists $r\in(0,r_0)$ such that 
\begin{equation}\label{estimate_delta0}
\int_{B_{\frac{1}{2r}(0)}}U^{2^*_s}(y)dy\geq(S_{N,s})^{\frac{N}{2s}}-\frac{\Theta_0}{2}.
\end{equation}
Furthermore, using the definition of $v_\eps$ in \eqref{defvepsilon} and \eqref{H_s_bound} we get a uniform $\mathcal D^{s}(\mathbb R^N)$ bound, namely,
\begin{equation}\label{bounds_w_eps2}\|v_\eps\|^2_{\mathcal D^s}=\int_{\Omega_\eps}v_\eps^{p_\eps+1}
=
\|u_\eps\|_{\infty}^{-\eps\frac{N-2s}{2s}}
\int_{\Omega}u_\eps^{p_\eps+1}
\leq \int_{\Omega}u_\eps^{p_\eps+1}
= \left(S_{N,s}\right)^{\frac{N}{2s}} +o(1).\end{equation}
Hence, by \eqref{estimate_delta0} and \eqref{bounds_w_eps2}, and using also Lemma \ref{lemma:wetoW}, we deduce
\begin{eqnarray}\label{hpLemVerified}\nonumber
\int_{\hat\Omega_\eps\cap B_{2r}(0)}\hat v_\eps^{2^*_s-\eps}(x)\,dx&=&\int_{\hat\Omega_\eps\cap B_{2r}(0)}\frac{1}{|x|^{2N-\eps(N-2s)}}v_\eps^{2^*_\eps-\eps}\!\left(\frac{x}{|x|^2}\right)dx\leq(2r)^{\eps(N-2s)}\int_{\hat\Omega_\eps\cap B_{2r}(0)}\frac{1}{|x|^{2N}}v_\eps^{2^*_s-\eps}\!\left(\frac{x}{|x|^2}\right)dx\\\nonumber
&=&(2r)^{\eps(N-2s)}\int_{\Omega_\eps\setminus B_{\frac1{2r}}(0)}v_\eps^{2^*_s-\eps}(y)\,dy=(1+o(1))\left[\int_{\Omega_\eps}v_\eps^{2^*_s-\eps}(y)\,dy-\int_{B_{\frac1{2r}}(0)}v_\eps^{2^*_s-\eps}(y)\,dy\right]\\\nonumber
&\overset{\eqref{bounds_w_eps2}}{\underset{\eqref{convv}}\leq}&(1+o(1))\left[(S_{N,s})^{\frac{N}{2s}}-\int_{B_{\frac1{2r}}(0)}U^{2^*_s}(y)\,dy+o(1)\right]\\
&\overset{\eqref{estimate_delta0}}{\leq}&\frac{\Theta_0}{2}(1+o(1))\leq\Theta_0,\qquad\text{for any $\eps$ sufficiently small.}
\end{eqnarray}
Thanks to \eqref{hpLemVerified}, we may then apply Lemma~\ref{lemma:Lemma_6_Han} to $\hat v_\eps$ in $\hat\Omega_\eps$ (with $b=b_\eps$, $q=2^*_s-\eps$ and $\beta=2^*_s$) and deduce that there exists $C=C(N,s,\beta)>0$ such that
\begin{equation}\label{int_vhat_1}
\int_{\hat\Omega_\eps\cap B_r(0)}\hat v_\eps^{\frac{(2_s^*)^2}{2}}(x)\,dx\leq 
\left( C r^{-\frac{2s}{2^*_s}}\left(\int_{\hat\Omega_\eps}\hat v_\eps^{2^*_s}(x)\,dx\right)^{\frac1{2^*_s}}
       \right)^{\frac{(2_s^*)^2}{2}}\!\!\!\!\!\!\!\!\!\!\!\!.
\end{equation}
Moreover, since $\|v_\eps\|_\infty=1$, and using \eqref{bounds_w_eps2},
\begin{equation}\label{int_vhat_2}
\int_{\hat\Omega_\eps}\hat v_\eps^{2_s^*}(x)\,dx\overset{\eqref{normInvarianceKelvin}}{=}
\int_{\Omega_\eps} v_\eps^{2_s^*}(y)\,dy\leq\|v_\eps\|_{\infty}^\eps
\int_{\Omega_\eps} v_\eps^{2_s^*-\eps}(y)\,dy=\int_{\Omega_\eps} v_\eps^{2_s^*-\eps}(y)\,dy\leq(S_{N,s})^{\frac{N}{2s}}+o(1),
\end{equation}
therefore, combining \eqref{int_vhat_1} and \eqref{int_vhat_2} we have that
\begin{equation}\label{int_vhat_3}
\int_{\hat\Omega_\eps\cap B_r(0)}\hat v_\eps^{\frac{(2_s^*)^2}{2}}(x)\,dx\leq C_1,
\end{equation}
for some $C_1>0$ independent of $\eps$.  We now write the equation in \eqref{vepsP} as
$$
(-\Delta)^s \hat v_\eps=a_\eps(x) \hat v_\eps\quad \mbox{ in }\hat\Omega_\eps
$$
with $a_\eps(x):=b_\eps(x)\hat v_\eps^{p_\eps-1}.$  Note that 
\begin{align*}
\theta_\eps
:=\frac{(2^*_s)^2}{2} \frac{1}{p_\eps-1}
>\frac{(2N)^2}{2(N-2s)^2} \frac{N-2s}{4s}
=\frac{N}{N-2s} \frac{N}{2s}>\frac{N}{2s}.    
\end{align*}
And, by \eqref{Uniform_bound_b},
\[|a_\eps(x)|^{\theta_\eps}=b_\eps(x)^{\theta_\eps}\hat v_\eps^{(p_\eps-1)\theta_\eps}\leq C\hat v_\eps^{(p_\eps-1)\theta_\eps} \quad\mbox{ for }x\in\hat\Omega_\eps.\]
Hence, 
\begin{align*}
\int_{B_r(0)}|a_\eps(x)|^{\theta_\eps}\, dx =\int_{B_r(0)\cap \hat\Omega_\eps}|a_\eps(x)|^{\theta_\eps}\, dx 
\leq C \int_{B_r(0)\cap \hat\Omega_\eps}(\hat v_\eps^{p_\eps-1}(x))^{\frac{(2^*_s)^2}{2} \frac{1}{p_\eps-1}}\,dx
= C \int_{B_r(0)\cap \hat\Omega_\eps}\hat v_\eps^{\frac{(2^*_s)^2}{2}}\!(x)\,dx\leq C,
\end{align*}
where the last inequality follows from  \eqref{int_vhat_3}.

Now, let $Q_r:=B_r(0)\times (0,r)\subset \R^{N+1}_+$ and let $W_\eps$ be the extension of $\hat v_\eps$ (see Subsection~\ref{subsec:ext}), then,
\begin{align*}
\operatorname{div}\left(t^{1-2s} \nabla W_\eps(x,t)\right)=0\ \text { for } (x,t)\in Q_r,
\quad -\lim _{t \rightarrow 0^{+}} t^{1-2s} \partial_t W_\eps(x, t)=k_s a_\eps(x) W_\eps(x, 0) { \text{ for } x\in B_r(0)\cap \hat \Omega_\eps.}
\end{align*}
By Proposition~\ref{prop:yyLiProp2.6}, Lemma~\ref{est:lem}, and \eqref{int_vhat_2}, 
\begin{align}
\|\hat v_\eps\|_{L^{\infty}(B_{r/2})}
\leq C_2\left(\int_{Q_r}t^{1-2s}(W^+_\eps)^{2^*_s}dx dt\right)^{\frac1{2^*_s}}
\leq C_3 \left(\int_{\R^N} \hat v_\eps^{2^*_s}(x)\, dx\right)^{\frac1{2^*_s}}
\leq C_4
\end{align}
for some positive constants $C_2,C_3,$ and $C_4$ independent of $\eps$.  This establishes \eqref{LimnftyBoundKelvin} and the proof is finished.
\end{proof}

\[\]

\begin{proposition}\label{prop:x_0isSimpleLowEnergy}
Let $\Omega$ be an open bounded smooth set and let $u_\eps$ be a family of low energy solutions of \eqref{problema}. Let  $x_0\in\Omega$ be the blow-up point in Lemma~\ref{lemma:blow-up}. Then $x_0$ is isolated and simple.
\end{proposition}
\begin{proof}
By Proposition~\ref{prop:Li4.16} it is enough to prove that $x_0$ is  isolated. Let $x_\eps$ be the sequence defined in \eqref{xepsdef}.  From  Lemma~\ref{lemma:blow-up} we know that $x_\eps$ is a blowing-up sequence of maximum points for $u_\eps$ converging to the blow-up point $x_0$. Next we show the existence of $C>0$, $r<1$ such that
\begin{align}\label{claimIsolatedBlowUp}
    u_\eps(x)\leq \frac{C}{|x-x_\eps|^{\frac{2s}{p_\eps-1}}}, \mbox{ for } x\in B_r(x_\eps).
\end{align}
As a consequence, it follows that $x_0$ is isolated. 

In order to prove \eqref{claimIsolatedBlowUp}, let $r<1$ be such that $B_{2r}(x_0)\subset\Omega$, then for any $\eps$ small $B_r(x_\eps)\subset\Omega$. Furthermore, let us define
\[r_\eps:= \|u_\eps\|_{\infty}^{-\frac{p_\eps-1}{2s}}\quad \rightarrow 0 \ \mbox{ as }\eps\rightarrow 0,\]
hence $B_{r_\eps}(x_\eps)\subset B_r(x_\eps)$. We divide the proof of \eqref{claimIsolatedBlowUp} into $2$ cases.
\\\\
{\bf Case 1.} If $x\in B_{r_\eps}(x_\eps)$.  The definition of $r_\eps$ immediately implies that
\[u_\eps (x)\leq \|u_\eps\|_{\infty}=\frac{1}{r_\eps^{\frac{2s}{p_\eps-1}}}\leq \frac{1}{|x-x_\eps|^{\frac{2s}{p_\eps-1}}}.\]
{\bf Case 2.} If $x\in B_{r}(x_\eps)\setminus B_{r_\eps}(x_\eps)$. The proof is a consequence of Lemma \ref{lemma:Lemma3_HAN}, indeed  by \eqref{pointwiseBoundU_eps}
\[
    u_\eps(x)\leq \Lambda \|u_\eps\|_{\infty}U\left((x-x_\eps)\|u_\eps\|_{\infty}^{\frac{p_\eps-1}{2s}} \right),
    \]
    namely, recalling the definition of $U$ in \eqref{Bubble}
\begin{equation}\label{esti1}
    u_\eps(x)\leq  \frac{\widetilde\Lambda}{|x-x_\eps|^{N-2s}\|u_\eps\|_{\infty}^{\frac{p_\eps-1}{2s}(N-2s)-1}}=\frac{\widetilde\Lambda}{|x-x_\eps|^{\frac{2s}{p_\eps-1}}}A_\eps,
    \end{equation}
    for a constant $\widetilde\Lambda>0$, where \[A_\eps:=\frac{1}{|x-x_\eps|^{N-2s-\frac{2s}{p_\eps-1}} \|u_\eps\|_{\infty}^{\frac{p_\eps-1}{2s}(N-2s)-1}}.\]
    The definition of $r_\eps$ implies that
\[|x-x_\eps|^{N-2s-\frac{2s}{p_\eps-1}}\geq r_\eps^{N-2s-\frac{2s}{p_\eps-1}}=
   \|u_\eps\|_{\infty}^{-\frac{p_\eps-1}{2s}(N-2s) +1},\]
as a consequence
  \[A_\eps\leq 1,\]
 hence, from \eqref{esti1} we deduce that
\[u_\eps(x)\leq  \frac{\widetilde\Lambda}{|x-x_\eps|^{\frac{2s}{p_\eps-1}}}.
    \]
\end{proof}

\begin{corollary}
\label{cor:OnceMissing}
Let $u_\eps$ be a family of low energy solutions to \eqref{problema}. Then,
\begin{align*}
\lim_{\eps\to 0}\|u_\eps\|^\eps_\infty=1.
\end{align*}
\end{corollary}

\begin{proof}
It is a consequence of Proposition~\ref{prop:x_0isSimpleLowEnergy} and the general result stated in Lemma~\ref{lemma:Li4.8}.
\end{proof}

\begin{proposition}\label{Analog_Lemma3.7_BMS}
Let $u_\eps$ be a family of low energy solutions to \eqref{problema}
. Let us define, up to a subsequence,
\begin{align}\label{gamma0:Def}
\gamma_{N,s}:={\lim_{\eps\rightarrow 0} }\|u_\eps\|_{\infty}\int_{\Omega} u_\eps^{p_\eps}(y)\,dy. \end{align}
Then 
\begin{equation}\label{ExactVAlueGamma0}
 \gamma_{N,s}=
 \frac{\sigma_N}{2} (\tau_{N,s})^{-\frac{N}{2}} B\Big(\frac{N}{2},s \Big)
 =2^N \pi ^{N/2} \Gamma (s)\frac{\Gamma
   \left(\frac{N}{2}+s\right)^{\frac{N}{2 s}-1}}{ \Gamma \left(\frac{N}{2}-s\right)^{\frac{N}{2 s}}},
 \end{equation}
 where the constant $\tau_{N,s}$ is defined in \eqref{def:constanTau}, $\sigma_N=|\partial B_1|=2\frac{\pi^{\frac{N}{2}}}{\Gamma (\frac{N}{2})}$, and $B(a,b)=\frac{\Gamma(a ) \Gamma(b)}{ \Gamma(a + b)}$ is the Beta function, so in particular $\gamma_{N,s}$ does not depend on the subsequence. Furthermore,
\begin{eqnarray}&
\displaystyle \lim_{\eps\rightarrow 0} \|u_\eps\|_{\infty}u_\eps(x)=  
 \gamma_{N,s} \mathcal G(x,x_0)
 \qquad\mbox{ in  }\ C_{loc}(\overline\Omega\setminus \{x_0\}),\label{laPrima}\\
&\displaystyle  \lim_{\eps\rightarrow 0} \frac{\|u_\eps\|_{\infty}u_\eps(x)}{d^s(x)}=  
 \gamma_{N,s}\frac{\mathcal G(x,x_0)}{d^s(x)}
 \qquad\mbox{ in  }\ C_{loc}(\overline\Omega\setminus \{x_0\}).\label{laSeconda}\end{eqnarray}
 where $x_0$ is the blow-up point in Lemma~\ref{lemma:blow-up}.
\end{proposition}
\begin{proof} 
We show first that $\gamma_{N,s}\in(0,+\infty)$ and that \eqref{ExactVAlueGamma0} holds. By a change of variables, recalling the definition of
$v_\eps$ in \eqref{defvepsilon}
\begin{align*}
\|u_\eps\|_\infty\int_\Omega u_\eps^{p_\eps}(y)\,dy
&=
\|u_\eps\|_\infty^{p_\eps+1-\frac{N}{2s}(p_\eps-1)}
\int_{\Omega_\eps}v_\eps^{p_\eps}(x)\,dx
=
\left(\|u_\eps\|_\infty^\eps\right)^{\frac{N-2s}{2s}}
\int_{\Omega_\eps}v_\eps^{p_\eps}(x)\,dx.
\end{align*}
By Corollary~\ref{cor:OnceMissing}, we have $\|u_\eps\|_\infty^\eps\to1$ as $\eps\to0$.  By Lemmas~\ref{lemma:Lemma3_HAN} and \ref{lemma:wetoW} and dominated convergence, $\int_{\Omega_\eps}v_\eps^{p_\eps}(x)\,dx
\longrightarrow
\int_{\R^N}U^{2_s^*-1}(x)\,dx$  as $\eps\to 0.$ Hence, 
\begin{align*}
\lim_{\eps\to0}
\|u_\eps\|_\infty\int_\Omega u_\eps^{p_\eps}(y)\,dy
&=
\int_{\R^N}U^{2_s^*-1}(x)\,dx=
\sigma_N\int_0^\infty
\frac{\rho^{N-1}}
     {(1+\tau_{N,s}\rho^2)^{\frac{N+2s}{2}}}\,d\rho
     =
\frac{\sigma_N}{2}
(\tau_{N,s})^{-\frac N2}
B\left(\frac N2,s\right).
\end{align*}
Using the definition of $\tau_{N,s}$ in
\eqref{def:constanTau} and the identities $\sigma_N=\frac{2\pi^{N/2}}{\Gamma(N/2)}$ and $B\left(\frac N2,s\right)
=
\frac{\Gamma(N/2)\Gamma(s)}
     {\Gamma(N/2+s)},$ we obtain
\begin{align*}
\gamma_{N,s}
&=
\frac{\sigma_N}{2}
(\tau_{N,s})^{-\frac N2}
B\left(\frac N2,s\right)=
2^N\pi^{\frac N2}\Gamma(s)
\frac{
\Gamma\left(\frac N2+s\right)^{\frac{N}{2s}-1}
}{
\Gamma\left(\frac N2-s\right)^{\frac{N}{2s}}
}
\in(0,+\infty).
\end{align*}
We are now in the position to prove \eqref{laSeconda}. Notice that  \eqref{laPrima} can be obtained similarly. \\
Let us define
\begin{equation}\label{defw_eps}
w_\eps(x):=\|u_\eps\|_{\infty}u_\eps(x),\qquad x\in \rn, 
\end{equation}
which solves
\begin{equation}\label{eq:w_eps}
(-\Delta)^s w_\eps=\|u_\eps\|_{\infty}u_\eps^{p_\eps}= \frac{1}{\|u_\eps\|_{\infty}^{p_\eps-1}}w_\eps^{p_\eps}
\quad \mbox{in }\Omega,\qquad 
w_\eps>0\quad \mbox{in }\Omega,\qquad
w_\eps=0 \quad \mbox{in }\mathbb R^N\setminus\Omega.
\end{equation}
Let $r>0$, we show that
\begin{align}\label{lemma:boundednessOutisex_0}
\|w_\eps\|_{L^{\infty}(\Omega\setminus B_r(x_0))}\leq C(r)\qquad \text{for all $\eps>0$ sufficiently small.}  
\end{align}
Indeed, if $x\in \Omega\setminus B_r(x_0)$ then $|x-x_\eps|\geq \frac{r}{2}$ for all $\eps$ small enough and, by \eqref{pointwiseBoundU_eps}, \eqref{UIntro}, and Corollary~\ref{cor:OnceMissing},
\begin{align*}
w_\eps(x)
&=\|u_\eps\|_{\infty}u_\eps(x)
\leq \Lambda \|u_\eps\|^2_{\infty} U\left((x-x_\eps)\|u_\eps\|_{\infty}^{\frac{p_\eps-1}{2s}} \right)
\leq \Lambda \|u_\eps\|^2_{\infty} \left(1+\tau_{N,s} \|u_\eps\|_{\infty}^{\frac{p_\eps-1}{s}} |x-x_\eps|^2\right)^{-\frac{N-2s}{2}}\\
&\leq \Lambda \left(\tau_{N,s} \frac{r^2}{4}\right)^{-\frac{N-2s}{2}}
 \|u_\eps\|_{\infty}^{2-\frac{N-2s}{2}\frac{p_\eps-1}{s}}
=\Lambda \left(\tau_{N,s} \frac{r^2}{4}\right)^{-\frac{N-2s}{2}}
 (\|u_\eps\|_{\infty}^\eps)^{\frac{N-2s}{2s}}\leq C(r)
\end{align*}
for some constant $C(r)>0$ independent of $\eps$, which proves  \eqref{lemma:boundednessOutisex_0}.
Furthermore, using the definition of $w_\eps$ in \eqref{defw_eps} and Lemma~\ref{lemma:blow-up}, one easily sees that \eqref{lemma:boundednessOutisex_0} implies
\begin{equation}\label{uGoesto0Outisidex_0}
    \|u_\eps\|_{L^{\infty}(\Omega\setminus B_r(x_0))}=o_\eps(1)\qquad \text{ as }\eps\to 0.
\end{equation}
Let $\Omega'$ be any neighborhood of $\partial\Omega$, such that $\Omega'$ does not contain $x_0$ and let us fix $\sigma>0$ such that $\overline{\Omega'}\cap \overline{B_{\sigma}(x_0)}=\emptyset$. We claim that
\begin{equation}\label{claim_oct_0}
\frac{w_\eps(x)}{d^s(x)}\to \gamma_{N,s}\frac{\mathcal G(x,x_0)}{d^s(x)}\qquad \text{uniformly in $\Omega'$,}
\end{equation}
namely that, for any $\eta>0$,
\begin{equation}\label{claim_oct}
\sup_{x\in\Omega'}\left|\frac{w_\eps(x)}{d^s(x)}-\gamma_{N,s}\frac{\mathcal G(x,x_0)}{d^s(x)}\right|<\eta\qquad \text{for any $\eps$ sufficiently small.}
\end{equation}
From the uniform continuity of $(x,y)\mapsto \frac{\mathcal G(x,y)}{d^s(x)}$ in $\overline{\Omega'}\times\overline{B_{\sigma}(x_0)}$, we deduce that there exists $r_\eta\in(0,\sigma)$ such that 
\begin{equation}\label{Gxx0}
\left|\frac{\mathcal G(x,y)}{d^s(x)}-\frac{\mathcal G(x,x_0)}{d^s(x)}  \right|\leq \frac{\eta}{2\gamma_{N,s}}
\qquad \text{for any $x\in\Omega'$ and for any $y\in B_{2r_\eta}(x_0)$}.    
\end{equation}
In order to prove \eqref{claim_oct} we use the equation \eqref{eq:w_eps} and the Green's representation formula of $w_\eps$, obtaining
\begin{align}
\sup_{x\in\Omega'}&\left|\frac{w_\eps(x)}{d^s(x)}-\gamma_{N,s}\frac{\mathcal G(x,x_0)}{d^s(x)}\right|=
\sup_{x\in\Omega'}\left|\int_\Omega\frac{\mathcal G(x,y)\|u_\eps\|_\infty u_\eps^{p_\eps}(y)}{d^s(x)}dy-\gamma_{N,s}\frac{\mathcal G(x,x_
0)}{d^s(x)}\right|\nonumber\\
&\leq\sup_{x\in\Omega'}\left|\int_{\Omega\setminus B_{r_\eta}(x_0)}\!\!\!\frac{\mathcal G(x,y) w_\eps(y) u_\eps^{p_\eps-1}(y)}{d^s(x)}dy\right|+
\sup_{x\in\Omega'}\left|\int_{B_{r_\eta}(x_0)}\!\!\frac{\mathcal G(x,y) \|u_\eps\|_\infty u_\eps^{p_\eps}(y)}{d^s(x)}dy-\gamma_{N,s}\frac{\mathcal G(x,x_0)}{d^s(x)}\right|.
\label{twoterms_oct}
\end{align}
By \eqref{Green1}, \eqref{lemma:boundednessOutisex_0}, \eqref{uGoesto0Outisidex_0} and the integrability of $|x-y|^{s-N}$ in $B_{2\diam(\Omega)}$,
\begin{equation}\label{firstterm_oct}
\sup_{x\in\Omega'}\left|\int_{\Omega\setminus B_{r_\eta}(x_0)}\!\!\!\frac{\mathcal G(x,y)w_\eps(y)u_\eps^{p_\eps-1}(y)}{d^s(x)}dy\right|\leq o_\eps(1)\\\sup_{x\in\Omega'}\left|\int_{\Omega\setminus B_{r_\eta}(x_0)}\frac{C}{|x-y|^{N-s}}dy\right|=o_\eps(1).
\end{equation}
Moreover, 
\begin{align}
\sup_{x\in\Omega'}&\left|\int_{B_{r_\eta}(x_0)}\!\!\frac{\mathcal G(x,y) \|u_\eps\|_\infty u_\eps^{p_\eps}(y)}{d^s(x)}dy-\gamma_{N,s}\frac{\mathcal G(x,x_0)}{d^s(x)}\right|\notag\\
&\overset{\eqref{Gxx0}}{\leq}
\frac{\eta}{2\gamma_{N,s}}\left|\int_{B_{r_\eta}(x_0)}  \|u_\eps\|_\infty u_\eps^{p_\eps}(y)\,dy\right|+\sup_{x\in\Omega'}\left|
\frac{\mathcal G(x,x_0)}{d^s(x)}
\left[
\int_{\Omega} \|u_\eps\|_\infty u_\eps^{p_\eps}(y)\,dy-\gamma_{N,s}-\int_{\Omega\setminus B_{r_\eta}(x_0)} \|u_\eps\|_\infty u_\eps^{p_\eps}(y)\,dy
\right]
\right|\nonumber\\
&{\leq}\frac{\eta}{2}(1+o_\eps(1))+o_\eps(1)
\leq \frac\eta2+o_\eps(1),
\label{secondterm_oct}
\end{align}
where in the last estimate we have used \eqref{Green2}, the definition of $\gamma_{N,s}$, \eqref{lemma:boundednessOutisex_0}, and \eqref{uGoesto0Outisidex_0}.  Combining \eqref{twoterms_oct}, \eqref{firstterm_oct}, and \eqref{secondterm_oct} we conclude that \eqref{claim_oct} holds true for any $\eps$ sufficiently small. \eqref{laSeconda} follows from it, recalling the
definition of $w_\eps$ in \eqref{defw_eps}.
\end{proof}

\begin{proposition} \label{Prop:CanPassTotheLimit}
Let $u_\eps$ be a family of low energy solutions to \eqref{problema}. Then, up to a subsequence, 
\[\lim_{\eps \rightarrow 0} \eps \|u_\eps\|_{\infty}^2=\gamma_{N,s}^2 \frac{2N }{(N-2s)\left(S_{N,s}\right)^{\frac{N}{2s}}}\mathcal R(x_0) \ ( >0), \]
where  the constant  $\gamma_{N,s}$ is the same as in Proposition~\ref{Analog_Lemma3.7_BMS}, $S_{N,s}$ is as in \eqref{constant:SnsVALUE} and   $\mathcal R$ is the fractional Robin function, and $x_0\in\Omega$ is the blow-up point in Lemma~\ref{lemma:blow-up}.
\end{proposition}
\begin{proof}From the fractional Pohozaev identity \eqref{Pohozaev}, choosing $y_0=x_0$ we get
\begin{equation}\label{PohozaevExplicit}\eps\frac{(N-2s)^2}{2N-\eps (N-2s)}\int_{\Omega}u_\eps^{p_\eps+1}(x)\,dx=\Gamma^2 (1+s)\int_{\partial\Omega}\left(\frac{u_\eps(x)}{d^s(x)}\right)^2\langle x-x_0,\nu (x)\rangle d\sigma_x.\end{equation}
Multiplying both sides by $\|u_\eps\|_\infty$, \eqref{PohozaevExplicit} becomes
\begin{equation}\label{PohozaevExplicitW}\eps\frac{(N-2s)^2}{2N-\eps (N-2s)}\|u_\eps\|_{\infty}^2\int_{\Omega}u_\eps^{p_\eps+1}(x)\,dx=\Gamma^2 (1+s)\int_{\partial\Omega}\left(\frac{\|u_\eps\|_{\infty}u_\eps(x)}{d^s(x)}\right)^2\langle x-x_0,\nu (x)\rangle d\sigma_x.\end{equation}
For the right-hand side of \eqref{PohozaevExplicitW}, Proposition~\ref{Analog_Lemma3.7_BMS} and the dominated convergence theorem imply that, up to a subsequence,
\begin{equation}\label{RHPoho}\lim_{\eps\rightarrow 0}
\Gamma^2 (1+s)\int_{\partial\Omega}\left(\frac{\|u_\eps\|_{\infty}u_\eps(x)}{d^s}\right)^2\langle x-x_0,\nu (x)\rangle d\sigma_x=
\gamma_{N,s}^2 \Gamma^2 (1+s)\int_{\partial\Omega}\left(\frac{\mathcal G(x,x_0)}{d^s(x)}\right)^2\langle x-x_0,\nu (x)\rangle d\sigma_x=
\gamma_{N,s}^2
\mathcal R(x_0)(N-2s),
\end{equation}
where the last equality follows from the Pohozaev identity for the fractional Green's function in Lemma~\ref{lemma:PohozaevGreen_DiebIanni}.
 For the left-hand side of \eqref{PohozaevExplicitW}, recalling the definition of 
$v_\eps$ (see \eqref{defvepsilon}), one has
\begin{eqnarray}
\label{LHPoho}\eps\frac{(N-2s)^2}{2N-\eps (N-2s)}\|u_\eps\|_{\infty}^2\int_{\Omega}u_\eps^{p_\eps+1}(x)\,dx
&=&
\eps\frac{(N-2s)^2}{2N-\eps (N-2s)}\|u_\eps\|_{\infty}^{2+\eps\frac{N-2s}{2s}}\int_{\Omega_\eps}v_\eps^{p_\eps+1}(x)\,dx\nonumber\\
&=& \eps \frac{(N-2s)^2}{2N}\|u_\eps\|_{\infty}^{2} \left(\int_{\mathbb R^N}U^{\frac{2N}{N-2s}}\,dx +o(1)\right),
\end{eqnarray}
where we used Corollary~\ref{cor:OnceMissing}, that $v_\eps\rightarrow U$ (Lemma~\ref{lemma:wetoW}), and the dominated convergence theorem, thanks to Lemma~\ref{lemma:Lemma3_HAN}.
Therefore, substituting \eqref{RHPoho} and \eqref{LHPoho} into \eqref{PohozaevExplicitW},
\[
\lim_{\eps\rightarrow 0}\eps \|u_\eps\|_{\infty}^2=\frac{\gamma_{N,s}^2 \mathcal R(x_0)}{\frac{(N-2s)}{2N}\int_{\mathbb R^N}U^{\frac{2N}{N-2s}}\,dx}.
\]
The conclusion follows recalling that $\int_{\mathbb R^N} U^{\frac{2N}{N-2s}}\,dy=\left(S_{N,s}\right)^{\frac{N}{2s}}$ (see \eqref{ExactValue2}). 
\end{proof}

\begin{remark}
  Proposition \ref{Prop:CanPassTotheLimit}
may be rephrased as $u_\eps(x_\eps)\sim C\eps^{-1/2}$, which is sharper than the blowup rate estimate 
for  isolated and simple blowup points that one can deduce
from the general theory developed 
in \cite{JLX14} (see Lemma \ref{lemma:Li4.8}).
\end{remark}

\begin{lemma}\label{lemmaCriticalRobinx_0}
Let $u_\eps$ be a family of low energy solutions to \eqref{problema} and let $x_0\in\Omega$ be the blow-up point in  Lemma~\ref{lemma:blow-up}. Then
\begin{equation}\label{CriticalPointRobin}\nabla \mathcal R(x_0)=0,\end{equation}
 where $\mathcal R$ is the fractional Robin function (see \eqref{def:Robin}).
\end{lemma}
\begin{proof}
We consider the Pohozaev identity \eqref{Pohozaev} and multiply both sides by $\|u_\eps\|_{\infty}$,
\begin{equation}\label{PohozaevTemp}\eps\|u_\eps\|_{\infty}^2\frac{(N-2s)^2}{2N-\eps (N-2s)}\int_{\Omega}u_\eps^{p_\eps+1}(x)\,dx=\Gamma^2 (1+s)\int_{\partial\Omega}\left(\frac{\|u_\eps\|_{\infty}u_\eps(x)}{d^s(x)}\right)^2\langle x-y_0,\nu (x)\rangle d\sigma_x\end{equation}
for any $y_0\in\mathbb R^N$. 
Using Proposition~\ref{Prop:CanPassTotheLimit} and \eqref{H_s_bound}, we can pass to the limit in the left hand side of \eqref{PohozaevTemp} as $\eps\rightarrow 0$, getting
\begin{equation}\label{PohozaevExplicitLimitB}\lim_{\eps\rightarrow 0}\eps\|u_\eps\|_{\infty}^2\frac{(N-2s)^2}{2N-\eps (N-2s)}\int_{\Omega}u_\eps^{p_\eps+1}(x)\,dx=\gamma_{N,s}^2 (N-2s)\mathcal R(x_0).\end{equation}
Furthermore, using Proposition~\ref{Analog_Lemma3.7_BMS} and the dominated convergence theorem, we can also pass to the limit in the right hand side of \eqref{PohozaevTemp} as $\eps\rightarrow 0$, which
gives
\begin{equation}\label{PohozaevExplicitLimitA}
\lim_{\eps\rightarrow 0}
\Gamma^2 (1+s)\int_{\partial\Omega}\left(\frac{\|u_\eps\|_{\infty}u_\eps(x)}{d^s(x)}\right)^2\langle x-y_0,\nu (x)\rangle d\sigma_x=
\Gamma^2 (1+s)\gamma_{N,s}^2\int_{\partial\Omega}\left(\frac{\mathcal G(x,x_0)}{d^s(x)}\right)^2\langle x-y_0,\nu (x)\rangle d\sigma_x
\end{equation}
for any $y_0\in\mathbb R^N$. Substituting   \eqref{PohozaevExplicitLimitA} and \eqref{PohozaevExplicitLimitB} into \eqref{PohozaevTemp},
\[
\int_{\partial\Omega}\left(\frac{\mathcal G(x,x_0)}{d^s(x)}\right)^2\langle x-y_0,\nu (x)\rangle d\sigma_x= \frac{(N-2s)\mathcal R(x_0)}{\Gamma^2(1+s)}
,\]
namely, 
\begin{equation}\label{PohozaevExplicitLimit}
\int_{\partial\Omega}\left(\frac{\mathcal G(x,x_0)}{d^s(x)}\right)^2\langle x,\nu (x)\rangle d\sigma_x- \frac{(N-2s)\mathcal R(x_0)}{\Gamma^2(1+s)} =\int_{\partial\Omega}\left(\frac{\mathcal G(x,x_0)}{d^s(x)}\right)^2\langle y_0,\nu (x)\rangle d\sigma_x
,\end{equation} for any $y_0\in\mathbb R^N$.
As a consequence, from the arbitrariness of $y_0$, one can take $y_0=e_i,-e_i$, $i=1,\ldots, N$,  where $\{e_i: \ i=1,\ldots, N\}$ is the canonical orthonormal basis of $\mathbb R^N$, getting 
\[\int_{\partial\Omega}\left(\frac{\mathcal G(x,x_0)}{d^s(x)}\right)^2\langle e_i,\nu (x)\rangle d\sigma_x= \int_{\partial\Omega}\left(\frac{\mathcal G(x,x_0)}{d^s(x)}\right)^2\langle -e_i,\nu (x)\rangle d\sigma_x=-\int_{\partial\Omega}\left(\frac{\mathcal G(x,x_0)}{d^s(x)}\right)^2\langle e_i,\nu (x)\rangle d\sigma_x,\]
hence
\[\int_{\partial\Omega}\left(\frac{\mathcal G(x,x_0)}{d^s(x)}\right)^2\langle e_i,\nu (x)\rangle d\sigma_x=0, \quad i=1,\ldots, N,\]
namely,
\[\int_{\partial\Omega}\left(\frac{\mathcal G(x,x_0)}{d^s(x)}\right)^2\nu (x) d\sigma_x=0.\]
Equality \eqref{CriticalPointRobin} finally follows by the representation formula for the gradient of the fractional Robin function obtained by \cite{DjitteSueur} (see Lemma~\ref{lemma:PohozaevGreen_DjitteSeur}): 
\begin{equation}\label{fractionaGreenPohozaev1}\nabla \mathcal R(x_0)=\Gamma^2(1+s)\int_{\partial\Omega} \left(\frac{\mathcal G(x,x_0)}{d^s(x)}\right)^2\nu (x) d\sigma_x.\end{equation}
\end{proof}

\begin{remark}
The localization of $x_0$ as a critical point of the Robin function has been deduced here using the Pohozaev type identity \eqref{fractionaGreenPohozaev1} for the fractional Green's function $\mathcal G$
recently proved in \cite[Corollary 1.7]{DjitteSueur}, which extends to the  fractional Laplacian the analogous result for the classical Laplacian proved in \cite[Theorem 4.4]{BrezisPeletier}. 
\end{remark}

We are ready to prove Theorem ~\ref{fractional:Han}, which follows from the various partial results obtained throughout the section. 

\begin{proof}[Proof of Theorem~\ref{fractional:Han}] 
Claim \eqref{blow-upintro} is proved in  Lemma~\ref{lemma:blow-up}; moreover, from Lemma~\ref{lemma:blow-up} and Lemma~\ref{lemmaCriticalRobinx_0}, it follows that  $x_\eps\rightarrow x_0$, where $x_0$ is an interior point of $\Omega$, and satisfies $\nabla\mathcal R(x_0)=0$, respectively.  Claim \eqref{eps_normueps2} is a consequence of Proposition~\ref{Prop:CanPassTotheLimit}.  Properties \eqref{gamma0_Ga} and 
\eqref{gamma0_G_delta} follow from Proposition~\ref{Analog_Lemma3.7_BMS}. 
The estimate \eqref{BuonduIntro} follows from Lemma~\ref{lemma:Lemma3_HAN} (see Remark~\ref{rmk:bound}). 
Finally, \eqref{convvIntro} is contained  in Lemma~\ref{lemma:wetoW}.    
\end{proof}

\subsection{The case of a unit ball}\label{ball:sec}
In this Section we show that if $\Omega=B_1(0)$ is the unit ball centered at zero then the estimates in Theorem~\ref{fractional:Han} are significantly more explicit.

\begin{proof}[Proof of Corollary~\ref{ball:thm}]
By a moving planes argument \cite{JW14}, we know that $u_\eps$ is radially symmetric and decreasing in the radial variable. Therefore $x_0=0$. The Robin function $\mathcal R(x_0)=\mathcal R(0)$ can be computed as follows. If $H$ denotes the regular part of the Green's function, then we have that $H(x,0)$ solves
\begin{align*}
    (-\Delta)^s H(x,0)=0\quad \text{ in }B_1,\qquad 
    H(x,0)=\Phi(|x|)\quad \text{ for }x\in \rn\backslash B_1,
\end{align*}
where $\Phi(|x|)=a_{N,s}|x|^{2s-N}$ is the fundamental solution for the fractional Laplacian and $a_{N,s}=\frac{\Gamma(\frac{N}{2}-s)}{4^s\pi^\frac{N}{2}\Gamma(s)}$.  Then, $H(\cdot,0)$ can be obtained with the fractional Poisson kernel $\Gamma_s$ for the ball (see, for instance, \cite{N15,B99,AJS18}), namely,
\begin{align*}
H(x,0)
=\int_{\rn\backslash B_1}\Gamma_s(x,y)\Phi(|y|)\, dy
=\widetilde\gamma_{N,s}\int_{\rn\backslash B_1}\frac{(1-|x|^2)^s}{|x-y|^N(|y|^2-1)^s}a_{N,s}|y|^{2s-N}\, dy
\end{align*}
for $x\in B_1$ with $\widetilde\gamma_{N,s}=\frac{2}{\Gamma(s)\Gamma(1-s)\sigma_N}$ and $\sigma_N=|\partial B_1|=\frac{2\pi^\frac{N}{2}}{\Gamma(\frac{N}{2})}$. Hence,
\begin{align}
\mathcal R(0)
&=H(0,0)
=\widetilde\gamma_{N,s}a_{N,s}\int_{\rn\backslash B_1}\frac{|y|^{2s-2N}}{(|y|^2-1)^s}\, dy
=\widetilde\gamma_{N,s}a_{N,s}\sigma_N\int_{1}^\infty\frac{\rho^{2s-N-1}}{(\rho^2-1)^s}\, d\rho\notag\\
&=\frac{2}{\Gamma(s)\Gamma(1-s)}
\frac{\Gamma(\frac{N}{2}-s)}{4^s\pi^\frac{N}{2}\Gamma(s)}
\frac{\Gamma \left(\frac{N}{2}\right) \Gamma (1-s)}{2 \Gamma
   \left(\frac{N}{2}-s+1\right)}
   =\frac{2^{1-2 s} \Gamma \left(\frac{N}{2}\right)}{(N-2 s) \Gamma (s)^2\pi ^{\frac{N}{2}}}.\label{R0}
\end{align}
Now, \eqref{1ball} follows by Theorem~\ref{fractional:Han} and \eqref{ExactVAlueGamma0}, since
\begin{align}
    \lim_{n \rightarrow +\infty} \eps_n \|u_{\eps_n}\|_{\infty}^2
    &=\gamma_{N,s}^2 \frac{2N }{(N-2s)\left(S_{N,s}\right)^{\frac{N}{2s}}}\mathcal R(x_0)
   =\frac{2^{N-2 s+2}N! }{(N-2 s)^2} \frac{\Gamma
   \left(\frac{N}{2}+s\right)^{\frac{N}{2 s}-2}}{\Gamma \left(\frac{N}{2}-s\right)^{\frac{N}{2 s}}}=:C(N,s).\label{Cdef}
\end{align}

Next, recall that the \emph{Martin kernel} of $(-\Delta)^s$ in $B_1$ is given by
\[
M_s(x,z)
:= \lim_{B \ni y \to z} \frac{\mathcal G(x,y)}{(1-|y|^2)^s}
= \frac{2^{1-2s}}{s\Gamma(s)^2\sigma_N}\,\frac{(1-|x|^2)_+^s}{|x-z|^N},
\qquad
x \in \mathbb{R}^N,\ z \in \partial B_1,\ x \neq z,
\]
see, for example, \cite{N15,AJS18,AJS182,B99}. Note also that $(1-|y|^2)^s=(1-|y|)^s(1+|y|)^s=d^s(y)(1+|y|)^s.$ Then, by Theorem~\ref{fractional:Han},
\begin{align*}
\lim_{n \rightarrow +\infty}\|u_{\eps_n}\|_{\infty} \frac{u_{\eps_n}}{d^s}(e_1)
 &=  
 2^s\gamma_{N,s}M_s(0,e_1)
 =2^s\left(2^N \pi ^{N/2} \Gamma (s)\frac{\Gamma
   \left(\frac{N}{2}+s\right)^{\frac{N}{2 s}-1}}{ \Gamma \left(\frac{N}{2}-s\right)^{\frac{N}{2 s}}}\right)
   \frac{2^{1-2s}}{s\Gamma(s)^2\sigma_N}\\
  &=\frac{s 2^{N-s} B\left(\frac{N}{2},s\right)}{\Gamma (s+1)^2}
  \left(\frac{\Gamma\left(\frac{N}{2}+s\right)}{\Gamma \left(\frac{N}{2}-s\right)}\right)^{\frac{N}{2 s}}=:D(N,s).
\end{align*}
Namely, using \eqref{Cdef},
\begin{align*}
\lim_{n \rightarrow +\infty}\eps^{-\frac{1}{2}} \frac{u_{\eps_n}}{d^s}(e_1)
 &=  \frac{D(N,s)}{\sqrt{C(N,s)}}
 =\frac{2^{\frac{N}{2}-1} (N-2 s) \Gamma \left(\frac{N}{2}\right) 
  }{ \Gamma(s+1) \sqrt{N!}}\left(\frac{\Gamma \left(\frac{N}{2}+s\right)}{\Gamma\left(\frac{N}{2}-s\right)}\right)^{\frac{N}{4 s}}.
\end{align*}
This yields \eqref{2ball} and ends the proof.
\end{proof}

\begin{remark}\label{ball:rmk}
We emphasize that the formulas in Corollary~\ref{ball:thm} also make sense for $s=1$ and, after some simplifications, these formulas recover the well-known explicit constants in the local case (see, for instance, \cite{AP87,H91} and \cite[eqs. (1.6) and (1.7)]{GST20})
\begin{align*}
\lim_{\eps \to 0} \eps \|u_\eps\|_{\infty}^2
 &=\frac{4 ((N-2) N)^{\frac{N-2}{2}} \Gamma (N)}{(N-2) \Gamma \left(\frac{N}{2}\right)^2}
 ;\\
 \lim_{\eps\to 0}\eps^{-\frac{1}{2}}|\nabla u_\eps(z)|
 &=
 \frac{1}{2} (N-2) ((N-2) N)^{\frac{N-2}{4}} \sqrt{\frac{(N-2) \Gamma \left(\frac{N}{2}\right)^2}{\Gamma(N)}}\quad  \text{ for all }z\in \partial B_1.
\end{align*}
\end{remark}

\section{Uniqueness and nondegeneracy in GNN type domains}
\subsection{Positive solutions are low energy solutions in symmetric domains}\label{positive:section}

\subsubsection{The origin is an isolated blow-up point in a symmetric domain}

In this section we assume that $\Omega$ is a smooth bounded  open set of $\mathbb R^N$, $N\geq 2$ which satisfies  \emph{(H1)} and \emph{(H2)}.

\begin{remark}\label{Omega_starshaped}
    Any $\Omega$ satisfying \emph{(H1)} and \emph{(H2)} is star-shaped with respect to the origin.
\end{remark}
\begin{remark}\label{rem:symmetry}
    If $\Omega$ satisfies \emph{(H1)} and \emph{(H2)}, any nonnegative solution $u$ to \eqref{equn} has reflectional symmetry with respect to the hyperplanes $\{x_i=0\}$ and it is strictly monotone decreasing in $|x_i|$ (where $x=(x_1,\ldots,x_N)$) for any $i\in\{1,\ldots,N\}$. This follows from a moving plane argument, see \cite[Corollary 1.2]{JW14}.   In particular, $\|u\|_{\infty}=u(0)$, and $x=0$ is the only point in $\Omega$ such that $\nabla u(x)=0$.
\end{remark}

The following lemma shows that the origin is a blow-up point.
\begin{lemma}\label{lem:yyl1}
    Let $s\in (0,1)$ and let $\Omega$ be a smooth bounded  open set of $\mathbb R^N$, $N\geq 2$ which satisfies  (H1) and (H2).
    Let $u_\eps$ be a solution of \eqref{equn}-\eqref{defp_eps}. Then there exists a sequence $\eps_n\rightarrow 0$ as $n\rightarrow +\infty$,  such that
\begin{equation}\label{2.4}
\|u_{\eps_n}\|_\infty=u_{\eps_n}(0)\rightarrow +\infty\qquad \text{as $n\rightarrow +\infty$.}
\end{equation}
\end{lemma}
\begin{proof}
By Remark~\ref{rem:symmetry} we know that all solutions take the maximum at the origin, namely $\|u_\eps\|_\infty=u_\eps(0)$. The claim now follows from  Lemma~\ref{blowup:positive:lem}.
\end{proof}

Next we show that the origin $\bar x=0$ is an isolated blow-up point, with blowing-up sequence $\bar x_\eps\equiv 0$.
\begin{proposition}\label{prop:0isIsolatedBU}
   Let $s\in (0,1)$ and let $\Omega$ be a smooth bounded domain of $\R^N$, $N\geq 2$, satisfying (H1) and (H2).
Let $u_\eps$ be a solution of \eqref{equn}-\eqref{defp_eps}, then there exists a constant $C_0=C_0(N,s,\Omega)>0$  such that
\begin{equation}
    u_\eps(x)|x|^{\frac{2s}{p_\eps-1}}\leq C_0\quad\mbox{for }x\in \Omega
\end{equation}
\end{proposition}
\begin{proof}
    By contradiction, suppose that there exist a sequence $\eps_n\rightarrow 0$ as $n\rightarrow +\infty$ such that
    \[M_n:=\sup_{\Omega} |x|^{\frac{2s}{p_{n}-1}}u_{n}(x)\rightarrow +\infty, \mbox{ as }n\rightarrow +\infty\]
    where $p_n:=p_{\eps_n}$, and $u_n:=u_{\eps_n}$.
    Let $Q_n\in \Omega$ be a sequence of points such that
\begin{equation}\label{Mnblowsup}M_n=|Q_n|^{\frac{2s}{p_{n}-1}}u_{n}(Q_n)\rightarrow +\infty.\end{equation}
    In particular, $u_{n}(Q_n)\rightarrow +\infty$, and thus, since by Lemma~\ref{lemma:boundary estimate} $u_n$ is bounded (uniformly on $n$) in a neighborhood of the boundary of $\Omega$, it follows  that  $Q_n$ do not accumulate at the boundary of $\Omega$. Therefore, there exists  $Q_0\in \Omega$ such that, up to a subsequence $Q_n\rightarrow Q_0$, and there exists $\delta>0$  such that  $B_{\delta}(Q_n)\subset\Omega$. Let $\rho\geq 2$ such that $\frac{|Q_n|}{\rho}<\delta$, for any $n$.
    Let us define
    \begin{equation}
w_{n}(y):=\frac{1}{u_{n}(Q_n)}u_{n}\left(
\frac{y}{u_{n}(Q_n)^{\frac{p_{n}-1}{2s}
}}+Q_{n}\right).
\end{equation}
Then $w_n(0)=1$ and $w_{n}$ satisfies
\begin{equation}
    \label{eq:wnLi}
\left\{
\begin{array}{ll}
(-\Delta)^s w_{n}=w_{n}^{p_{n}}& \mbox{in }\Omega_{n},\\
\hspace{1cm}w_{n}>0 & \mbox{in }\Omega_{n},\\
\hspace{1cm}w_{n}=0 & \mbox{in }\mathbb R^N\setminus\Omega_{n},
\end{array}
\right.\end{equation}
where $\Omega_{n}:= u_{n}(Q_n)^{\frac{p_{n}-1}{2s}
}\left(\Omega-Q_{n}\right)
$.

Let us consider the balls \[\widetilde B_{n}:=\{|y|\leq  M_n^{\frac{p_{n}-1}{2s}
}/\rho
\}.\]
By definition, $\widetilde B_n$ converges to $\mathbb R^N$ as $n\rightarrow +\infty$. Next we show that  $\widetilde B_n\subset\Omega_n$ and that there exists  a constant $C>0$ such that
\begin{equation}\label{boundinEnlarginset}w_n(y)\leq C\quad \mbox{ for }y\in \widetilde B_{n}.\end{equation}
Indeed for $y\in \widetilde B_n$, $x:=\frac{y}{u_{n}(Q_n)^{\frac{p_{n}-1}{2s}
}}+Q_{n}$ satisfies 
\begin{equation} \label{palla}
|x-Q_n|\leq \frac{|Q_n|}{\rho}<\delta,
\end{equation}
and so, by the definition of $\delta$, the ball $\widetilde B_{n}\subset\Omega_n$. Furthermore,  for $y\in\widetilde B_n$, 
 the triangular inequality implies from \eqref{palla} that \[|x|\geq \left(1-\frac{1}{\rho}\right)|Q_n|\geq\frac{|Q_n|}{\rho},
\]
from which it follows  that
\[\left(\frac{|Q_n|}{\rho}\right)^{\frac{2s}{p_{n}-1}}u_{n}(x)\leq|x|^{\frac{2s}{p_{n}-1}}u_{n}(x)\leq |Q_n|^{\frac{2s}{p_{n}-1}}u_{n}(Q_n),\]
namely that
\[w_n(y)=\frac{u_n(x)}{u_n(Q_n)}\leq \rho^{\frac{2s}{p_{n}-1}},\]
which gives \eqref{boundinEnlarginset}.

Let $j\in \mathbb N$, then, by \eqref{boundinEnlarginset}, we may assume that $w_n\leq C$ in $B_{4j}$, up to a subsequence. Let $W_n$ be the harmonic extension of $w_n$ (see Section~\ref{subsec:ext}). Using Proposition~\ref{2.6} with $a(x)=w_n^{p_n-1}$ (note that $a\in L^\infty (B_{4j})$ and $W_n\in H(t^{1-2s}, Q_{4j})$), we know that 
\begin{equation*}
\sup_{Q_{2j}} W_n \le C \inf_{Q_{2j}} W_n \leq C\inf_{B_{2j}} w_n\leq C_j\qquad \text{ for all }n\in \mathbb N.
\end{equation*}

Next, by Propositions~\ref{2.6} and~\ref{prop:Li2.19}, there is $\alpha_j\in(0,1)$ and $C>0$ independent of $n$ such that 
\[
\|W_n\|_{C^{\alpha_j}(\overline{Q}_{j})} + \|\nabla_x W_n\|_{C^{\alpha_j}(\overline{Q}_{j})} + \|t^{1-2s} \partial_t W_n\|_{C^{\alpha_j}(\overline{Q}_{j})} \le C_j.
\]
Hence, for $\alpha_j'\in (0,\alpha_j)$ using the compact embedding of $C^{\alpha_j}(\overline{Q}_{j})\hookrightarrow C^{\alpha_j'}(\overline{Q}_{j})$, there is $W_0$ such that 
\begin{align*}
    W_n\to W_0,\quad 
    \nabla_x W_n\to \nabla_x W_0,\quad 
    t^{1-2s} \partial_t W_n\to t^{1-2s} \partial_t W_0\qquad \text{ in }C^{\alpha_j'}(\overline{Q}_{j})
\end{align*}
as $n\to\infty.$ With a standard diagonalization argument, we may define $W_0$ in the whole of $\mathbb {R}^{N+1}_+.$ Hence, $W_0 \in H_{\mathrm{loc}}(t^{1-2 s}, \overline{\mathbb{R}^{N+1}_+})$ (see \eqref{Hlocdef}). Furthermore, for $\Psi\in C^\infty_c(\overline{\mathbb R^{N+1}_+})$, $W_n$ satisfies that 
\begin{align*}
    \int_{\mathbb {R}^{N+1}_+}t^{1-2s} \nabla W_n \nabla \Psi\, dxdt = k_s \int_{\rn} W_n(x,0)^{p_n}\Psi(x,0)\, dx.
\end{align*}
Then, letting $n\to\infty$,
\begin{align*}
    \int_{\mathbb {R}^{N+1}_+}t^{1-2s} \nabla W_0 \nabla \Psi\, dxdt = k_s \int_{\rn} W_0(x,0)^{2^*_s-1}\Psi(x,0)\, dx.
\end{align*}
Using that $\Psi$ is arbitrary, this implies that $W_0$ is a weak solution of \eqref{crit}. Hence, Theorem~\ref{1.8} implies that 
 $W_0(\cdot,0)=U_{t,a}$, where $a \in\mathbb R^N$, $t>0$, and  $U_{t,a}$ is one of the bubbles defined in \eqref{tutteLeSoluzEqLimite}. In particular, $w_n\to U_{t,a}$ pointwise in $\rn$ as $n\rightarrow \infty.$ Furthermore, since $w_n(0)=1$, then
\begin{equation}
1=U_{t,a}(0)=
    \left(\frac{t}{t^2+\tau_{N,s}a^2}\right)^{\frac{N-2s}{2}}\end{equation}
then  $t=t^2+\tau_{N,s}a^2$, hence  $t\leq 1$ and $|a|^2=\frac{t(1-t)}{\tau_{N,s}}\leq (\tau_{N,s})^{-1}$.

By \eqref{boundinEnlarginset}, $w_n\leq C$ in $B_1(a)$, then using a priori regularity estimates \cite{RS14,fall2020regularity} (see the proof of Lemma~\ref{blowup:positive:lem}) we have that $w_{n}\rightarrow U_{t,a}$ in $C^{1}(B_1(a))$.  Since $a$ is a strict maximum point for the bubble $U_{t,a}$,  we obtain that there exists a sequence of points $y_n\in B_1(a)$ such that $\nabla w_n (y_n)=0$. 
Furthermore, under assumptions \emph{(H1)} and \emph{(H2)}, we know that $\nabla u_n(x)=0$ if and only if $x=0$ (see Remark \eqref{rem:symmetry}). This implies that
\[0=\frac{y_n}{u_{n}(Q_n)^{\frac{p_{n}-1}{2s}
}}+Q_{n},\]
and then that
\[|y_n|=|Q_n|u_{n}(Q_n)^{\frac{p_{n}-1}{2s}
}=M_n^{\frac{p_{n}-1}{2s}
}\rightarrow +\infty,\]
a contradiction with $y_n\in B_1(a)$.
\end{proof}

\subsubsection{Proof of Theorem~\ref{PositiveAreLow:thm}}

\begin{proof}[Proof of Theorem~\ref{PositiveAreLow:thm}]
With the previous lemmas in this section, we can follow the strategy in \cite[Theorem 2.1]{G00}, in particular, in order to prove the result,  it is enough to show that 
\begin{equation}\label{claimProofAllLowEnergy}
\|u_\eps\|^{p_\eps+1}_{p_\eps+1}\to \left(S_{N,s}\right)^{\frac{N}{2s}} \qquad\mbox{ as } \eps\to 0,
\end{equation}
since \eqref{claimProofAllLowEnergy} implies that
\[
\frac{\|u_\eps\|_{\mathcal D^s}^2}{\left(\int_{\Omega}\left|u_\eps\right|^{2^*_s-\eps}\right)^{\frac{2}{2^*_s-\eps}}} =\|u_\eps\|_{p_\eps+1}^{p_\eps-1}\to  S_{N,s}\qquad \text{ as }\eps\to 0,
\]
namely that \eqref{limSn} holds.

Next, we show \eqref{claimProofAllLowEnergy}. By Lemma~\ref{lem:yyl1}, $\|u_\eps\|_\infty=u_\eps(0)\rightarrow \infty$ as $\eps\rightarrow 0$, namely  $\bar x=0$ is a blow-up point with blowing-up sequence $\bar x_{\eps_n}\equiv 0$. Moreover, by Proposition~\ref{prop:0isIsolatedBU} we know that $\bar x=0$ is an isolated blow-up point.  The general Proposition~\ref{prop:Li4.16} then implies that $\bar x=0$ is an isolated and simple blow-up point. 
As a consequence, the general Lemma~\ref{lemma:Li4.8} and Proposition~\ref{prop:Li4.9} apply, from which we have that
\begin{equation}\label{BigControl}
\|u_\eps\|_\infty^\eps=1+o(1)\quad \text{ as }\eps\to 0;
\end{equation}
and, furthermore, that there exist constants $C,r>0$ independent of $\eps$ such that  
\begin{equation}\label{yyliIneq}
u_\eps(x)\leq \frac{C}{\|u_\eps\|_{\infty}|x|^{N-2s}}\qquad \text{ if } 0<|x|\leq r.
\end{equation}
Let us split
\[\|u_\eps\|^{p_\eps+1}_{p_\eps+1}=\int_{B_\delta(0)} u_\eps^{p_\eps+1}+
\int_{\Omega\backslash B_\delta(0)} u_\eps^{p_\eps+1}=:I_{\eps,1}+I_{\eps,2}
\]
for some $0<\delta<r$. By \eqref{yyliIneq}, $u_\eps(x) \le \frac{C}{\|u_\eps\|_\infty} \frac{1}{\delta^{N-2s}}$ if $|x|=\delta,$
and, since $u_\eps$ is monotone decreasing along the coordinate axis (see Remark~\ref{rem:symmetry}), we obtain that $u_\eps(x) \le \frac{C}{\|u_\eps\|_\infty} \frac{1}{\delta^{N-2s}}$ for $|x|>\delta$. Hence,
\begin{align*}
I_{\eps,2}
&:=
\int_{\Omega\backslash B_\delta(0)} u_\eps^{p_\eps+1}
\leq \frac{C}{\|u_\eps\|_\infty^{p_\eps+1}} \frac{1}{\delta^{(N-2s)(p_\eps+1)}}|\Omega|\longrightarrow 0\qquad \text{ as }\eps\to 0.
\end{align*}

In order to compute $I_{\eps,1}$, let $\Omega_\eps := \|u_\eps\|_\infty^{\frac{p_\eps-1}{2s}} \Omega$ and let $v_\eps(x) := \frac{1}{\|u_\eps\|_\infty} \, u_\eps\!\left( x\|u_\eps\|_\infty^{-\frac{p_\eps-1}{2s}} \right)$ for $x\in \Omega_\eps$.  Arguing as in Lemma~\ref{lemma:wetoW}, we obtain that $v_\eps \to U$ in $C^1(K)$ for every compact set $K$ of $\mathbb{R}^N$, where $U$ is given by \eqref{UIntro}.  Furthermore,  by \eqref{yyliIneq} and \eqref{BigControl},
\[
v_\eps(x) \le
\frac{C}{ \|u_\eps\|_\infty^{\eps\frac{N-2s}{2 s}}}\frac{1}{|x|^{N-2s}}
=\frac{C}{(1+o(1))}\frac{1}{|x|^{N-2s}},\qquad \text{if $0<|x|<\delta \|u_\eps\|_\infty^{\frac{p_\eps-1}{2s}}$}
\]
as $\eps\to 0.$  
Then, by dominated convergence, and using \eqref{BigControl},  we get
\begin{align*}
I_{\eps,1}
&:=
\int_{B_\delta(0)} u_\eps^{p_\eps+1}
=
\|u_\eps\|_\infty^{\eps\frac{N-2s}{2s}}
\int_{\mathbb{R}^N}
v_\eps^{p_\eps+1}
\chi_{|x|<\delta \|u_\eps\|_\infty^{\frac{p_\eps-1}{2s}}}\longrightarrow 
\int_{\mathbb{R}^N} U^{2^*_s}=\left(S_{N,s}\right)^{\frac{N}{2s}}
\end{align*}
as $\eps\to 0$, where the last equality follows from the explicit computation in \eqref{ExactValue2}.

As a consequence, $\|u_\eps\|^{p_\eps+1}_{p_\eps+1}=I_{\eps,1}+I_{\eps,2}\to \left(S_{N,s}\right)^{\frac{N}{2s}}$ as $\eps\to 0$, as claimed.

\end{proof}

\subsection{Uniqueness in symmetric domains}\label{section:uniqueness}
In this section we consider smooth domains $\Omega\subset\R^N$ satisfying \emph{(H1)} and \emph{(H2)} (stated in the Introduction) and we prove our main uniqueness result Theorem~\ref{thm:uniqueness1}.

\begin{proof}[Proof of Theorem~\ref{thm:uniqueness1}]
Suppose, by contradiction, that there exists a sequence $p_{\eps_n}:=\frac{N+2s}{N-2s}-\eps_n\nearrow\frac{N+2s}{N-2s}$ and two sequences of different functions $u_{\eps_n}$, $\bar u_{\eps_n}$ which are solutions of \eqref{equn} with $p$ replaced by $p_{\eps_n}$. By Theorem \ref{PositiveAreLow:thm}, $u_{\eps_n}$ and $\bar u_{\eps_n}$ are families of low energy solutions.  By Remark~\ref{rem:symmetry}, 
\begin{align}\label{xeps0}
x_{\eps_n}=0\qquad \text{for all $n\in \mathbb N$.}
\end{align}

    Assume, without loss of generality, that 
    \begin{equation}\label{wlog}
\|u_{\eps_n}\|_{\infty}\geq\|\bar u_{\eps_n}\|_{\infty}\qquad\text{ for any $n\in\mathbb N$.}
    \end{equation}
 In order to simplify the notation in the rest of the proof we ignore the subscript $n$.

\medskip

\textbf{Step 1.} 
     Let $\widetilde z_\eps:\Omega_\eps\to \R$ be given by
\begin{align*}
\widetilde z_\eps(x):= u_\eps \left( \frac{x}{\|u_\eps\|_{\infty}^\frac{p_\eps-1}{2s}} \right) - \bar u_\eps \left( \frac{x}{\|u_\eps\|_{\infty}^\frac{p_\eps-1}{2s}} \right),
\end{align*}
    where $\Omega_\eps:= {\|u_\eps\|_{\infty}^\frac{p_\eps-1}{2s}} \left(\Omega\right)$. 

Note that $\|u_\eps\|_{\infty}\to+\infty$ and $\Omega_\eps\to\mathbb R^N$, by \eqref{blow-upintro} in Theorem~\ref{fractional:Han}. In addition, by \eqref{xeps0} and \eqref{eps_normueps2},
\begin{equation}\label{ratioLinftynorm}
\frac{\|u_\eps\|_\infty}{\|\bar u_\eps\|_\infty}\to 1\qquad\text{as $\eps\to0$.}
\end{equation}
Let $v_\eps$ be the rescaling of $u_\eps$ defined in \eqref{rescaling_Intro} and let $\bar v_\eps$ denote the corresponding rescaling of $\bar u_\eps$, namely,
\[
\bar v_\eps(x):=\frac{1}{\|\bar u_\eps\|_{\infty}}\bar u_\eps\left(\frac{x}{\|\bar u_\eps\|_{\infty}^{\frac{p_\eps-1}{2s}}}\right).
\]
If we consider the following different rescaling of $\bar u_\eps$,
\[
\bar w_\eps(x):=\frac{1}{\|u_\eps\|_{\infty}}\bar u_\eps\left(\frac{x}{\|u_\eps\|_{\infty}^{\frac{p_\eps-1}{2s}}}\right),
\]
we have that
\[
\bar w_\eps(x)=\frac{\|\bar u_\eps\|_{\infty}}{\|u_\eps\|_{\infty}}\bar v_\eps\left(\frac{\|\bar u_\eps\|_{\infty}^{\frac{p_\eps-1}{2s}}}{\| u_\eps\|_{\infty}^{\frac{p_\eps-1}{2s}}}x\right).
\]
Then, by \eqref{ratioLinftynorm} and \eqref{convvIntro}  applied to $\bar v_\eps$, we deduce that
\begin{equation}\label{whate_to_U}
\bar w_\eps\to U\qquad\text{in $C_{loc}^{2s-\delta}(\R^N)$ for any $\delta>0$ as $\eps\to0$,}
\end{equation}
where $U$ is the function in \eqref{UIntro}. Furthermore, by \eqref{ratioLinftynorm} and \eqref{BuonduIntro} applied to $\bar v_\eps$, and using that $U$ (defined  in \eqref{UIntro}) is radially decreasing, there exists $\bar \Lambda>\Lambda>0$ such that
\begin{align}\label{whate_estimate}
\bar w_\eps(x)
&\leq 2\bar v_\eps\left(\left(\frac{\|\bar u_\eps\|_{\infty}}{\| u_\eps\|_{\infty}}\right)^\frac{p_\eps-1}{2s}x\right)\leq2\Lambda U\left(\left(\frac{\|\bar u_\eps\|_{\infty}}{\| u_\eps\|_{\infty}}\right)^\frac{p_\eps-1}{2s}x\right)
\leq 2\Lambda U\left(\frac{x}{2}\right)\leq\bar \Lambda U(x).
\end{align}

Moreover,
\begin{align*}
(-\Delta)^s \widetilde z_\eps = c_\eps(x) \widetilde z_\eps\quad \text{ in }\Omega_\eps,\qquad
\widetilde z_\eps=0\quad \text{ in }\mathbb R^N\backslash \Omega_\eps,
\end{align*}
where
\begin{eqnarray*}
c_\eps(x)&:=&\frac{p_\eps}{\|u_\eps\|_\infty^{p_\eps-1}}  \int_0^1 \left( 
t\,u_\eps\left(\frac{x}{\|u_\eps\|_{\infty}^{\frac{p_\eps-1}{2s}}}\right)+(1-t)\,\bar u_\eps\left(\frac{x}{\|u_\eps\|_{\infty}^{\frac{p_\eps-1}{2s}}}\right)
\right)^{p_\eps-1}\, dt\\
&=&p_\eps \int_0^1 \left( t\,v_\eps(x)+(1-t)\,\bar w_\eps(x)\right)^{p_\eps-1}\, dt=p_\eps \int_0^1 A_\eps^{p_\eps -1}(t,x)\,dt,
\end{eqnarray*}
where $A_\eps(t,x):= t\,v_\eps(x)+(1-t)\,\bar w_\eps(x)$. Applying \eqref{BuonduIntro} and \eqref{whate_estimate}, 
\begin{equation}\label{heart}
A_\eps(t,x)\leq t\,\Lambda U(x)+(1-t)\,\bar \Lambda U(x)\leq \bar \Lambda U(x),\qquad c_\eps(x)\leq C U^{p_\eps-1}(x).
\end{equation}
Moreover, by \eqref{convvIntro} and \eqref{whate_to_U},
\begin{equation}
    \label{star2}
A_\eps(t,x)\to U(x)\quad\text{and}\quad c_\eps(x)\to(2^*_s-1)\,U^{\frac{4s}{N-2s}}(x)=\frac{2^*_s-1}{(1+\tau_{N,s}|x|^2)^{2s}},
\end{equation}
and the convergence is uniform on compact sets of $\mathbb R^N$.

Now, let $z_\eps:\Omega_\eps\to \R$ be given by $z_\eps(x)= \widetilde z_\eps(x) / \|\widetilde z_\eps\|_{\infty}$. Then,
\begin{align}\label{wneq}
(-\Delta)^s z_\eps = c_\eps(x)z_\eps\quad \text{ in }\Omega_\eps,\qquad
z_\eps=0\quad \text{ in }\Omega_\eps^c,\qquad 
\|z_\eps\|_{\infty}=1.
\end{align}

Using $z_\eps\in \mathcal H_0^s(\Omega_\eps)$ as a test function in \eqref{wneq},
\begin{equation}\label{sea1}
\|z_\eps\|^2_{\mathcal D^s}=\int_{\Omega_\eps}c_\eps(x) z_\eps^2(x)\,dx.
\end{equation}
Let us fix $\delta\in(0,\frac{4s}{N-2s})$, then applying Sobolev and H\"older inequalities,
\begin{eqnarray}\label{fiori}
S_{N,s}\left(\int_{\Omega_\eps}|z_\eps(x)|^{2^*_s}\,dx\right)^{\frac{2}{2^*_s}}&\leq&\|z_\eps\|^2_{\mathcal D^s}=\int_{\Omega_\eps}c_\eps(x) z_\eps^2(x)\,dx\overset{\|z_\eps\|_{\infty}=1}{\leq}\int_{\Omega_\eps}c_\eps(x) |z_\eps(x)|^{2-\delta}\,dx\nonumber\\
&\leq&\left(\int_{\Omega_\eps}|c_\eps(x)|^{\tfrac{2^*_s}{2^*_s-2+\delta}}\,dx\right)^{\tfrac{2^*_s-2+\delta}{2^*_s}} \left(\int_{\Omega_\eps}|z_\eps(x)|^{2^*_s}\,dx\right)^{\frac{2-\delta}{2^*_s}}
\!\!\!\!\!\!\!\!\!\!.
\end{eqnarray}
This, in turn, implies that
\begin{eqnarray}\label{tree}
\left(\int_{\Omega_\eps}|z_\eps(x)|^{2^*_s}\,dx\right)&\leq& C \left(\int_{\Omega_\eps}|c_\eps(x)|^{\tfrac{2^*_s}{2^*_s-2+\delta}}dx\right)^{\tfrac{2^*_s-2+\delta}{\delta}}\nonumber\\
&\overset{\eqref{heart}}{\leq}&C\left(\int_{\mathbb{R}^N}\left[\frac{1}{(1+\tau_{N,s}|x|^2)^{2s-\frac{(N-2s)\eps}{2}}}\right]^{\tfrac{2^*_s}{2^*_s-2+\delta}}dx\right)^{\tfrac{2^*_s-2+\delta}{\delta}}\leq C
\end{eqnarray}
for any $\eps$ sufficiently small, being $\delta\in(0,\frac{4s}{N-2s})$.

From this last chain of inequalities we can deduce that 
\begin{align*}
\|z_\eps\|_{L^{2^*}(\Omega_\eps)}\leq C
\qquad 
\text {and}\qquad 
\|c_\eps\|_{L^{\frac{2^*_s}{2^*_s-2+\delta}}(\Omega_\eps)}\leq C.
\end{align*}
Thus, from \eqref{fiori} and \eqref{sea1}, we obtain the uniform bound $\|z_\eps\|_{\mathcal D^s}\leq C$.

Reasoning as in the proof of Lemma~\ref{lemma:wetoW}, 
\begin{equation}\label{zepstoz}
z_\eps\to z\quad \text{in $C_{loc}^{2s-\delta}(\mathbb R^N)$\qquad  for any $\delta>0$},
\end{equation} 
where $z\in \mathcal  D^{s}(\R^N)$, being $\|z\|^2_{\mathcal D^s}\leq\liminf_{\eps\to0}\|z_\eps\|^2_{\mathcal D^s}\leq C$, and 
$z$ satisfies
\begin{equation}\label{ast}
\left\{
\begin{array}{ll}
(-\Delta)^s z=\frac{2^*_s-1}{(1+\tau_{N,s}|x|^2)^{2s}}z& \mbox{in }\mathbb{R}^N,\\
\hspace{0.35cm}\|z\|_{\infty}\leq1. & 
\end{array}
\right.\end{equation}

\

\textbf{Step 2.} In this step we show that, for any $r>0$, there exists $C_r>0$ such that
\begin{equation}\label{zepsestimate}
|z_\eps(x)|\leq \frac{C_r}{|x|^{N-2s}}\qquad\text{ for } x\in\Omega_\eps\setminus B_r(0)
\end{equation}
and for any $\eps$ sufficiently small. Let us consider the Kelvin transform $\hat{z}_\eps$ of ${z}_\eps$, which is defined as 
\[
\hat{z}_\eps(x)=\frac{1}{|x|^{N-2s}}z_\eps\left(\frac{x}{|x|^2}\right)\qquad \text{ for }x\in\hat\Omega_\eps\setminus\{0\},
\]
where $\hat\Omega_\eps$ is the image of $\Omega_\eps$ under the Kelvin transform.
By \eqref{sharm}, $\hat{z}_\eps$ solves
\begin{equation*}
\left\{
\begin{array}{ll}
(-\Delta)^s \hat{z}_\eps=a_\eps(x)\hat{z}_\eps
& \mbox{in }\hat\Omega_\eps,\\
\hspace{1cm}\hat{z}_\eps=0 & \mbox{in }\mathbb R^N\setminus\hat\Omega_\eps,
\end{array}\right.
\end{equation*}
where $a_\eps(x):=\frac{1}{|x|^{4s}}c_\eps(\frac{x}{|x|^2})$.

It is easy to see, by \eqref{heart}, that
\[
\frac{1}{|x|^{4s}}c_\eps\left(\frac{x}{|x|^2}\right)
\leq
\frac{C}{|x|^{(N-2s)\eps}(|x|^2+\tau_{N,s})^{2s-\frac{N-2s}{2}\eps}}
\leq
\frac{C}{(|x|^2+\tau_{N,s})^{2s-\frac{N-2s}{2}\eps}},
\]
where in the last inequality we have used \eqref{stima|x|} and that, by \eqref{eps_normueps2}, $\|u_\eps\|_\infty^\eps\to1$.
Hence, for $q>\frac{N}{2s}$,
\[
\|a_\eps\|^q_{L^q(\hat\Omega_\eps)}\leq C\int_{\mathbb R^N}\frac{dx}{(|x|^2+\tau_{N,s})^{(2s-\frac{N-2s}{2}\eps)q}}\leq C.
\]
Now, let $Q_\frac{2}{r}=B_\frac{2}{r}(0)\times (0,1)\subset\mathbb R^{N+1}_+$ and let $Z_\eps$ be the extension of $\hat z_\eps$ (see Subsection~\ref{subsec:ext}). Then, in particular,
\begin{align*}
\operatorname{div}\left(t^{1-2s} \nabla Z_\eps(x,t)\right)=0\ \text { for } (x,t)\in Q_\frac{2}{r},
\qquad -\lim _{t \rightarrow 0^{+}} t^{1-2s} \partial_t Z_\eps(x, t)=a_\eps(x) Z_\eps(x, 0) \quad \text { for } x\in B_\frac{2}{r}(0)\cap\hat \Omega_\eps.
\end{align*}
By Proposition~\ref{prop:yyLiProp2.6}, Lemma~\ref{est:lem}, and \eqref{tree}, 
\begin{align*}
\|\hat z_\eps\|_{L^{\infty}(B_\frac{1}{r}(0))}
&\leq C_{1,r}
\left(\int_{Q_\frac{2}{r}}t^{1-2s} |Z_\eps|^{2^*_s}dx dt\right)^{\frac1{2^*_s}}
\leq C_{2,r} \left(\int_{\R^N} |\hat z_\eps(x)|^{2^*_s}\, dx\right)^\frac{1}{2*_s}\\
&= C_{2,r}\left(\int_{\hat\Omega_\eps} |\hat z_\eps(x)|^{2^*_s}\, dx\right)^{\frac1{2^*_s}}= 
C_{2,r}\left(\int_{\Omega_\eps} |z_\eps(x)|^{2^*_s}\, dx\right)^{\frac1{2^*_s}}
\leq C_{3,r}
\end{align*}
for some positive constants $C_{1,r},C_{2,r},$ and $C_{3,r}$ independent of $\eps$.  By definition of $\hat z_\eps$ we obtain the desired estimate \eqref{zepsestimate}.

\

\textbf{Step 3.} The linear span of solutions to \eqref{ast} was described in \cite[Theorem 1.1]{DDS13}, where it was proved that a general solution to \eqref{ast} can be written as 
\[
z(x)=\alpha\frac{1-\tau_{N,s}|x|^2}{(1+\tau_{N,s}|x|^2)^{\tfrac{N-2s+2}{2}}}+\sum_{i=1}^N\beta_i\frac{2\sqrt{\tau_{N,s}}\,x_i}{(1+\tau_{N,s}|x|^2)^{\tfrac{N-2s+2}{2}}},\qquad \text{where $\alpha$, $\beta_i\in\mathbb R$.
}
\]

By definition and Remark~\ref{rem:symmetry}, $z_\eps$ is symmetric with respect to the coordinate hyperplanes, then, by \eqref{zepstoz}, $z$ has the same symmetry and so, necessarily, $\beta_i=0$ for any $i\in\{1,\ldots,N\}$, namely  
\begin{equation}\label{cross}
z(x)=\alpha\frac{1-\tau_{N,s}|x|^2}{(1+\tau_{N,s}|x|^2)^{\tfrac{N-2s+2}{2}}}\qquad \text{for some $\alpha\in\mathbb R$.}
\end{equation}
 In addition, by \eqref{wlog} $u_\eps(0)\geq\bar u_\eps(0)$, then $z_\eps(0)\geq0$, thus $\alpha\geq0$.

\

\textbf{Step 4.} Let us now prove that $\alpha=0$. In order to do so we argue by contradiction assuming that $\alpha>0$.

We start computing
\begin{eqnarray}\label{c0}
\frac{N+2s}{N-2s}\int_{\mathbb R^N}\frac{1-\tau_{N,s}|y|^2}{(1+\tau_{N,s}|y|^2)^{\frac{N+2s+2}{2}}}dy&=&\frac{N+2s}{N-2s}\tau_{N,s}^{-\frac N2}\int_{\mathbb R^N}\frac{1-|w|^2}{(1+|w|^2)^{\frac{N+2s+2}{2}}}dw\nonumber\\
&=&\frac{N+2s}{N-2s}\tau_{N,s}^{-\frac N2}\sigma_N\int_0^{+\infty}\frac{(1-\rho^2)\rho^{N-1}}{(1+\rho^2)^{\frac{N+2s+2}{2}}}d\rho\nonumber\\
&=&\frac{N+2s}{N-2s}\tau_{N,s}^{-\frac N2}\sigma_N\left[\int_0^{1}\frac{(1-\rho^2)\rho^{N-1}}{(1+\rho^2)^{\frac{N+2s+2}{2}}}d\rho-\int_1^{0}\frac{(1-\frac{1}{t^2})}{(1+\frac{1}{t^2})^{\frac{N+2s+2}{2}}t^{N+1}}dt\right]\nonumber\\
&=&\frac{N+2s}{N-2s}\tau_{N,s}^{-\frac N2}\sigma_N\int_0^{1}\frac{(1-\rho^2)\rho^{2s-1}(\rho^{N-2s}-1)}{(1+\rho^2)^{\frac{N+2s+2}{2}}}d\rho=:-c_0,
\end{eqnarray}
where $c_0$ is a positive constant, being the last integral negative and finite.  Also, let 
\begin{align*}
k_1:=\frac{N+2s}{N-2s}\int_{\mathbb R^N}\frac{|1-\tau_{N,s}|y|^2|}{(1+\tau_{N,s}|y|^2)^{\tfrac{N+2s+2}{2}}}\, dy
\end{align*}
and observe that $k_1>c_0.$

Let $\Omega'$ be any neighborhood of $\partial\Omega$, such that $\Omega'$ does not contain the origin and let us fix $\sigma>0$ such that $\overline{\Omega'}\cap \overline{B_{\sigma}(0)}=\emptyset$. We claim that
\begin{equation}\label{claim_Step4_0}
\|u_\eps\|^2_{\infty}\frac{(u_\eps(x)-\bar u_\eps(x))}{\|u_\eps-\bar u_\eps\|_{\infty} d^s(x)}\to -\alpha c_0\frac{\mathcal G(x,0)}{d^s(x)}\qquad \text{uniformly in $\Omega'$,}
\end{equation}
namely that, for any $\eta>0$,
\begin{equation}\label{claim_Step4}
\sup_{x\in\Omega'}\left|\|u_\eps\|^2_{\infty}\frac{(u_\eps(x)-\bar u_\eps(x))}{\|u_\eps-\bar u_\eps\|_{\infty} d^s(x)}+\alpha c_0\frac{\mathcal G(x,0)}{d^s(x)}\right|<\eta
\end{equation}
for any $\eps$ sufficiently small.

 Let us fix $\eta>0$. 
Then, from the uniform continuity of $(x,y)\mapsto \frac{\mathcal G(x,y)}{d^s(x)}$ in $\overline{\Omega'}\times\overline{B_{\sigma}(0)}$, we deduce that there exists $r_\eta\in(0,\sigma)$ such that, for any $x\in\Omega'$ and for any $y\in B_{r_\eta}(0)$,
\begin{equation}\label{Gx0}
\left|\frac{\mathcal G(x,y)}{d^s(x)}-\frac{\mathcal G(x,0)}{d^s(x)}  \right|\leq \frac{\eta}{2\alpha k_1}.    
\end{equation}
In order to simplify the notation we set
\[
\phi_\eps(x):=\|u_\eps\|^2_{\infty}\frac{u_\eps(x)-\bar u_\eps(x)}{\|u_\eps-\bar u_\eps\|_{\infty}},
\]
and we observe that
\begin{equation}\label{flower}
\left\{
\begin{array}{ll}
(-\Delta)^s \phi_\eps=f_\eps 
& \mbox{in }\Omega,\\
\hspace{1cm}\phi_\eps=0 & \mbox{in }\R^N\backslash \Omega,
\end{array}\right.   
\end{equation}
where
\[
f_\eps:=\frac{\|u_\eps\|^2_{\infty}}{\|u_\eps-\bar u_\eps\|_{\infty}}   g_\eps(x)
(u_\eps(x)-\bar u_\eps(x))\quad \text{with}\quad g_\eps(x)=p_\eps \int_0^1(t u_\eps(x)+(1-t)\bar u_\eps(x))^{p_\eps-1} dt.
\]
Since 
\begin{align*}
g_\eps(x)=\|u_\eps\|^{p_\eps-1}_{\infty}c_\eps(\|u_\eps\|^\frac{p_\eps-1}{2s}_{\infty} x),    
\end{align*}
using \eqref{heart} we obtain
\begin{equation}\label{smile}
g_\eps(x)\leq C\|u_\eps\|_{\infty}^{p_\eps-1}U^{p_\eps-1}(\|u_\eps\|_{\infty}^{\tfrac{p_\eps-1}{2s}}x)\leq \frac{C\|u_\eps\|_{\infty}^{-\tfrac{4s}{N-2s}+3\eps-\tfrac{N-2s}{2s}\eps^2}}{|x|^{4s-\eps(N-2s)}}.
\end{equation}
Next we estimate $f_\eps$ in $\Omega\setminus B_{r_\eta}(0)$.
By \eqref{zepsestimate}, \eqref{smile}, and being $\|u_\eps\|_\infty^\eps\to1$, by \eqref{eps_normueps2}, for $x\in\Omega\setminus B_{r_\eta}(0)$,
\begin{eqnarray}\label{feps_estimate}
|f_\eps(x)|&=&
\|u_\eps\|_{\infty}^2 z_\eps\left({x}{\|u_\eps\|_{\infty}^{\frac{p_\eps-1}{2s}}}\right)g_\eps(x)\leq
C \|u_\eps\|_{\infty}^2 \frac{g_\eps(x)}{\|u_\eps\|_{\infty}^{2-\tfrac{N-2s}{2s}\eps}|x|^{N-2s}}\nonumber\\
&\leq& C\frac{\|u_\eps\|_{\infty}^{-\tfrac{4s}{N-2s}}}{|x|^{N+2s-\eps(N-2s)}}\leq \frac{C}{{r_\eta}^{N+2s-\eps(N-2s)}}\|u_\eps\|_{\infty}^{-\tfrac{4s}{N-2s}}\leq C_{\Omega'}\|u_\eps\|_{\infty}^{-\tfrac{4s}{N-2s}}
\end{eqnarray}
for a constant $C_{\Omega'}>0$ independent of $\eps$.  Moreover,
\begin{equation}\label{Subclaim2}
\lim_{\eps\to0}\int_{B_{r_\eta}(0)} f_\eps(x)\,dx=-\alpha c_0.
\end{equation}
Indeed, by the change of variable $x=\tfrac{y}{\|u_\eps\|_{\infty}^{\tfrac{p_\eps-1}{2s}}}$, using \eqref{heart}, \eqref{star2}, \eqref{zepstoz}, \eqref{zepsestimate}, \eqref{cross}, \eqref{c0}, and the dominated convergence theorem,
\begin{eqnarray*}
\int_{B_{r_\eta}(0)} f_\eps(x)\, dx
&=&\frac{\|u_\eps\|_{\infty}^2}{\|u_\eps-\bar u_\eps\|_{\infty}}\int_{B_{r_\eta}(0)} g_\eps(x)(u_\eps(x)-\bar u_\eps(x))\, dx\notag\\
&=&\|u_\eps\|_{\infty}^2\int_{B_{r_\eta}(0)}
(\|u_\eps\|^{p_\eps-1}_{\infty}c_\eps(\|u_\eps\|^\frac{p_\eps-1}{2s}_{\infty} x))\frac{u_\eps(x)-\bar u_\eps(x)}{\|u_\eps-\bar u_\eps\|_{\infty}}\, dx\nonumber\\
&=&
\|u_\eps\|_{\infty}^{2+p_\eps-1-\frac{p_\eps-1}{2s}N}
\int_{B_{r_\eta\|u_\eps\|_\infty^{\frac{p_\eps-1}{2s}}}(0)}c_\eps(y) z_\eps(y)\, dy\\
&=&
\|u_\eps\|_{\infty}^{\eps \left(\frac{N-2s}{2s}\right)}
\int_{B_{r_\eta\|u_\eps\|_\infty^{\frac{p_\eps-1}{2s}}}(0)}c_\eps(y) z_\eps(y)\, dy
\to\alpha\frac{N+2s}{N-2s}\int_{\mathbb R^N}\frac{1-\tau_{N,s}|y|^2}{(1+\tau_{N,s}|y|^2)^{\tfrac{N+2s+2}{2}}}\, dy=-\alpha c_0
\end{eqnarray*}
as $\eps\to 0.$ A similar argument also shows that 
\begin{align}\label{Subclaim2aux}
\lim_{\eps\to0}\int_{B_{r_\eta}(0)} |f_\eps(x)|\, dx = \alpha k_1.
\end{align}

In order to prove \eqref{claim_Step4} we use the Green's representation formula of $\phi_\eps$, obtaining
\begin{eqnarray}\label{twoterms}
\sup_{x\in\Omega'}\left|\frac{\phi_\eps(x)}{d^s(x)}+\alpha c_0\frac{\mathcal G(x,0)}{d^s(x)}\right|&=&
\sup_{x\in\Omega'}\left|\int_\Omega\frac{\mathcal G(x,y)f_\eps(y)}{d^s(x)}dy+\alpha c_0\frac{\mathcal G(x,0)}{d^s(x)}\right|\nonumber\\
&\leq&\sup_{x\in\Omega'}\left|\int_{\Omega\setminus B_{r_\eta}(0)}\!\!\!\frac{\mathcal G(x,y)f_\eps(y)}{d^s(x)}dy\right|+
\sup_{x\in\Omega'}\left|\int_{B_{r_\eta}(0)}\!\!\frac{\mathcal G(x,y)f_\eps(y)}{d^s(x)}dy+\alpha c_0\frac{\mathcal G(x,0)}{d^s(x)}\right|.
\end{eqnarray}
On one hand, by \eqref{feps_estimate}, \eqref{Green1}, the integrability of $|x-y|^{s-N}$ in $B_{2\diam(\Omega)}$, and \eqref{blow-upintro}, we can estimate the first term in the right hand side of \eqref{twoterms} as follows
\begin{equation}\label{firstterm}
\sup_{x\in\Omega'}\left|\int_{\Omega\setminus B_{r_\eta}(0)}\!\!\!\frac{\mathcal G(x,y)f_\eps(y)}{d^s(x)}dy\right|\leq C_{\Omega'}\|u_\eps\|_{\infty}^{-\tfrac{4s}{N-2s}}\sup_{x\in\Omega'}\left|\int_{\Omega\setminus B_{r_\eta}(0)}\frac{C}{|x-y|^{N-s}}dy\right|\leq C_{\Omega'}\|u_\eps\|_{\infty}^{-\tfrac{4s}{N-2s}}=o_\eps(1).
\end{equation}

On the other hand,
\begin{eqnarray}\label{secondterm}
\sup_{x\in\Omega'}\left|\int_{B_{r_\eta}(0)}\!\!\frac{\mathcal G(x,y)f_\eps(y)}{d^s(x)}dy+\alpha c_0\frac{\mathcal G(x,0)}{d^s(x)}\right|&\overset{\eqref{Gx0}}{\leq}& 
\frac{\eta}{2\alpha k_1}\int_{B_{r_\eta}(0)} |f_\eps(y)|\,dy+
\sup_{x\in\Omega'}\left|
\frac{\mathcal G(x,0)}{d^s(x)}
\left[
\int_{B_{r_\eta}(0)} f_\eps(y)\,dy+\alpha c_0
\right]
\right|\nonumber\\
&\overset{\eqref{Subclaim2}, \eqref{Subclaim2aux}}{\underset{\eqref{Green2}}{\leq}}& \frac{\eta}{2}(1+o_\eps(1))+K o_\eps(1)
=\frac\eta2+o_\eps(1)\qquad \text{ as }\eps\to0.
\end{eqnarray}
Combining \eqref{twoterms}, \eqref{firstterm}, and \eqref{secondterm} we conclude that \eqref{claim_Step4} holds true for any $\eps$ sufficiently small.

\

In order to complete the proof of Step 4, we apply the Pohozaev identity \eqref{Pohozaev} to $u_\eps$ and $\bar u_\eps$ and we get
\[
\left[\frac{2N}{p_\eps+1}-(N-2s)\right]\int_{\Omega}u_\eps^{p_\eps+1}\,dx=\Gamma^2(1+s)\int_{\partial\Omega}\left(\frac{u_\eps(x)}{d^s(x)}\right)^2\langle x,\nu (x)\rangle d\sigma_x,
\]
\[
\left[\frac{2N}{p_\eps+1}-(N-2s)\right]\int_{\Omega}\bar u_\eps^{p_\eps+1}\,dx=\Gamma ^2(1+s)\int_{\partial\Omega}\left(\frac{\bar u_\eps(x)}{d^s(x)}\right)^2\langle x,\nu (x)\rangle d\sigma_x.
\]
Subtracting one from another and multiplying by $\frac{\|u_\eps\|_{\infty}^3}{\|u_\eps-\bar u_\eps\|_{\infty}}$ in both sides yields
\begin{equation}\label{fromPohozaev}
\begin{aligned}
&\Gamma^2(1+s)\int_{\partial\Omega}\frac{\|u_\eps\|_{\infty}^2(u_\eps(x)-\bar u_\eps(x))}{\|u_\eps-\bar u_\eps\|_{\infty}d^s(x)}\frac{(u_\eps(x)+\bar u_\eps(x))\|u_\eps\|_{\infty}}{d^s(x)}\langle x,\nu (x)\rangle d\sigma_x\qquad\hfill\\
&\qquad\qquad=\left[\frac{2N}{p_\eps+1}-(N-2s)\right](p_\eps+1)\int_{\Omega}\frac{\|u_\eps\|_{\infty}^3}{\|u_\eps-\bar u_\eps\|_{\infty}}(u_\eps(x)-\bar u_\eps(x))\int_0^1(tu_\eps(x)+(1-t)\bar u_\eps(x))^{p_\eps}dt dx.
\end{aligned}
\end{equation}
By doing the change of variable $x=y\|u_\eps\|_{\infty}^{-\tfrac{p_\eps-1}{2s}}$, the right hand side of \eqref{fromPohozaev} reduces to
\begin{align*}
[\text{RHS of \eqref{fromPohozaev}}]
&=\eps(N-2s)\|u_\eps\|_{\infty}^{2+\frac{N-2s}{2s}\eps}\int_{\Omega_\eps}z_\eps(y)\left[\int_0^1\left(t\,v_\eps(y)+(1-t)\,\bar w_\eps(y)\right)^{p_\eps}dt\right]dy\\
&=\eps(N-2s)\|u_\eps\|_{\infty}^{2+\frac{N-2s}{2s}\eps}\int_{\Omega_\eps}z_\eps(y)\left[\int_0^1\left(A_\eps(t,y)
\right)^{p_\eps}dt\right]dy.
\end{align*}
 By \eqref{eps_normueps2}, we have that $\|u_\eps\|_\infty^\eps\to1$ and, up to a subsequence, $\lim_{\eps\to0} \eps(N-2s)\|u_\eps\|_{\infty}^{2+\frac{N-2s}{2s}\eps}=C_1$, for some constant $C_1>0$.
Therefore, by \eqref{heart}, \eqref{star2}, \eqref{zepstoz}, \eqref{zepsestimate}, \eqref{cross}, and the dominated convergence theorem,
\[
\lim_{\eps\to0}\text{RHS}= C_1 \alpha \int_{\R^N}\frac{1-\tau_{N,s}|x|^2}{(1+\tau_{N,s}|x|^2)^{N+1}}dx.
\]
Performing the same computations as in \eqref{c0} we obtain that
\begin{align*}
\int_{\R^N}\frac{1-\tau_{N,s}|x|^2}{(1+\tau_{N,s}|x|^2)^{N+1}}dx=0.
\end{align*}

Hence,
\begin{equation}\label{RHS>0}
\lim_{\eps\to0}[\text{RHS of \eqref{fromPohozaev}}]=0.
\end{equation}
On the other hand, by \eqref{gamma0_G_delta}, \eqref{ratioLinftynorm}, \eqref{claim_Step4_0}, and Lemma~\ref{lemma:PohozaevGreen_DiebIanni},
\begin{equation}\label{LHS<0}
\lim_{\eps\to0}[\text{LHS of \eqref{fromPohozaev}}]=-2c_0\alpha \gamma_{N,s}\Gamma^2(1+s)
\int_{\partial\Omega}\left(\frac{\mathcal G(x,0)}{d^s(x)}\right)^2\langle x,\nu (x)\rangle \,d\sigma_x=-2 c_0\alpha\gamma_{N,s}(N-2s)\mathcal{R}(0) <0.
\end{equation}
Combining \eqref{fromPohozaev}, \eqref{RHS>0}, and \eqref{LHS<0} we get the desired contradiction and Step 4 is concluded.

\

\textbf{Step 5.} Let $y_\eps\in\mathbb R^N$ be such that
\begin{equation}\label{last}
|z_\eps(y_\eps)|=\|z_\eps\|_{\infty}=1.
\end{equation}
The previous Step implies that $z\equiv0$. Therefore, by \eqref{zepstoz}, $z_\eps\to0$ in $K$ for every compact set $K$ in $\R^N$ and this implies that $|y_\eps|\to+\infty$ as $\eps\to0$. This, together with \eqref{last}, gives a contradiction against \eqref{zepsestimate} and concludes the proof.

\end{proof}

\subsection{Non-degeneracy in symmetric domains}
\label{section:nondeg}
\begin{theorem} Let $N\geq2 $ and let $\Omega\subset \R^N$ be a smooth bounded domain satisfying (H1) and (H2).
Let $\eps_0>0$ be such that the least energy solution $u_\eps$ of \eqref{problema} is unique for any $\eps\in(0,\eps_0)$ and consider the eigenvalue problem 
\begin{equation}\label{EigenvalueProblem}
\left\{
\begin{array}{ll}
(-\Delta)^s \varphi -p_\eps u_\eps^{p_\eps-1}\varphi =\mu\varphi& \mbox{in }\Omega,\\
\hspace{2.85cm}\varphi=0 & \mbox{in }\mathbb R^N\setminus\Omega.
\end{array}
\right.\end{equation}
Then $\mu=0$ is not an eigenvalue for \eqref{EigenvalueProblem}, namely $u_\eps$ is non-degenerate.
\end{theorem}
\begin{proof}
For $\eps\in(0,\eps_0)$, let $u_\eps$ be the unique least energy solution to \eqref{problema} and, by contradiction, assume that there exists $\varphi_\eps\not\equiv0$ so that
\begin{equation*}
\left\{
\begin{array}{ll}
(-\Delta)^s \varphi_\eps -p_\eps u_\eps^{p_\eps-1}\varphi_\eps =0& \mbox{in }\Omega,\\
\hspace{3cm}\varphi_\eps=0 & \mbox{in }\mathbb R^N\setminus\Omega.
\end{array}
\right.\end{equation*}
Set
\[
\tilde\psi_\eps(x):=\varphi_\eps\left(\frac{x}{\|u_\eps\|_\infty^{\frac{p_\eps-1}{2s}}}\right)\quad \text{ for } x\in\Omega_\eps:= \|u_\eps\|_{\infty}^{\frac{p_\eps-1}{2s}}\Omega,\qquad \text{and}\qquad \psi_\eps:=\frac{\tilde\psi_\eps\;\:}{\|\tilde\psi_\eps\|_\infty}.
\]
By easy computations, 
\begin{equation}\label{equationpsieps}
\left\{
\begin{array}{ll}
(-\Delta)^s \psi_\eps =p_\eps v_\eps^{p_\eps-1}\psi_\eps& \mbox{in }\Omega_\eps,\\
\hspace{1cm}\psi_\eps=0 & \mbox{in }\mathbb R^N\setminus\Omega_\eps,
\end{array}
\right.\end{equation}
where $v_\eps$ is the rescaled function defined in \eqref{defvepsilon} with $x_\eps=0$. By Lemma~\ref{lemma:Lemma3_HAN}, we have that $p_\eps v_\eps^{p_\eps-1}\leq C U^{p_\eps-1}$ in $\Omega_\eps$. Reasoning as in Step 1 of the proof of Theorem~\ref{thm:uniqueness1} we can obtain the uniform bound $\|\psi_\eps\|_{\mathcal D^s}\leq C$ and that
\begin{equation}\label{psiepstoPsi}
\psi_\eps\longrightarrow\Psi\quad\text{in $C_{loc}^{2s-\delta}(\mathbb R^N)$\quad for any $\delta>0$,}
\end{equation}
where $\Psi\in \mathcal  D^{s}(\R^N)$, being $\|\Psi\|_{\mathcal D^s}^2\leq\liminf_{\eps\to0}\|\psi_\eps\|_{\mathcal D^s}^2\leq C$, and $\Psi$ satisfies
\begin{equation}\label{equationPsi}
\left\{
\begin{array}{ll}
(-\Delta)^s \Psi =\frac{2^*_s-1}{(1+\tau_{N,s}|x|^2)^{2s}}\,\Psi& \mbox{in }\R^N,\\
\hspace{0.35cm}\|\Psi\|_\infty\leq1. &
\end{array}
\right.\end{equation}
Moreover, following the proof of Step 2 of Theorem~\ref{thm:uniqueness1}, but using \eqref{equationpsieps}, we can show that, for any $r>0$, there exists $C_r>0$ such that
\begin{equation}\label{psiepsestimate}
|\psi_\eps(x)|\leq \frac{C_r}{|x|^{N-2s}}\qquad \text{ for all }x\in\Omega_\eps\setminus B_r(0).
\end{equation}

By \cite[Theorem 1.2]{DiebIanniSaldanaCalcVar} we deduce that $\varphi_\eps$ is symmetric with respect to the coordinate hyperplanes. Then, by \eqref{equationPsi}, we can proceed exactly as in Step 3 of the proof of Theorem~\ref{thm:uniqueness1} and  conclude that
\[
\Psi(x)=\alpha\frac{1-\tau_{N,s}|x|^2}{(1+\tau_{N,s}|x|^2)^{\tfrac{N-2s+2}{2}}}\qquad \text{for some $\alpha\in\mathbb R$. }
\]
Now, we suppose by contradiction that $\alpha\neq 0.$ Without loss of generality, assume that $\alpha>0$.  Let $c_0>0$ be as in \eqref{c0}. 

Arguing as in Step 4 of the proof of Theorem~\ref{thm:uniqueness1}, we obtain 
\begin{equation}\label{conv_nondeg}
    \lim_{\eps\rightarrow 0} \|u_\eps\|_{\infty}^2\frac{\varphi_\eps(x)}{\|\varphi_\eps\|_\infty d^s(x)}=  
 -\alpha c_0\frac{\mathcal G(x,0)}{d^s(x)}
 \qquad\mbox{ in  }\ C_{loc}(\overline\Omega\setminus \{0\}).
\end{equation}
 By Lemma~\ref{Lemma:boundaryTerm=0},
\begin{equation}\label{last1}
\int_{\partial\Omega}\|u_\eps\|_\infty\frac{u_\eps(x)}{d^s(x)}\|u_\eps\|_\infty^2\frac{\varphi_\eps(x)}{\|\varphi_\eps\|_\infty d^s(x)}\langle x,\nu\rangle\, d\sigma=0.
\end{equation}
But, by \eqref{laSeconda}, \eqref{ExactVAlueGamma0}, \eqref{last1}, and Lemma~\ref{lemma:PohozaevGreen_DiebIanni},
\begin{equation*}
\int_{\partial\Omega}\|u_\eps\|_\infty\frac{u_\eps(x)}{d^s(x)}\|u_\eps\|_\infty^2\frac{\varphi_\eps(x)}{\|\varphi_\eps\|_\infty d^s(x)}\langle x,\nu\rangle\, d\sigma\underset{\eps\to0}{\to}-\gamma_{N,s}\alpha c_0\int_{\partial\Omega}\left(\frac{G(x,0)}{d^s(x)}\right)^2\langle x,\nu\rangle \,d\sigma=-\gamma_{N,s}\alpha c_0\frac{(N-2s)R(0)}{\Gamma^2(1+s)}<0.
\end{equation*}
We have reached a contradiction. Hence, $\alpha=0$ and $\Psi\equiv0$.

Finally, let $y_\eps\in\mathbb R^N$ be such that
\begin{equation}\label{last3}
|\psi_\eps(y_\eps)|=\|\psi_\eps\|_{\infty}=1.
\end{equation}
Being $\Psi\equiv0$, by \eqref{psiepstoPsi}, $\psi_\eps\to0$ in $K$ for every compact set $K$ in $\R^N$ and this implies that $|y_\eps|\to+\infty$ as $\eps\to0$. This, together with \eqref{last3}, gives a contradiction against \eqref{psiepsestimate} and concludes the proof.
\end{proof}

\section{Multiplicity results}\label{Sec:MR}

In this section we show that there are low energy solutions that are not least energy solutions both in the fractional and in the local case (corresponding formally to the case $s=1$). 

\subsection{Examples in the fractional case}
Let us recall some notation. For $U\subset \rn$, let $\mathcal G_U$ denote the Green's function of the fractional Laplacian in $U$ and let $H_{U}$ be its regular part, namely, for $x \in U,$ $H_{U}(x,\cdot)$ is the solution of 
\begin{align*}
    (-\Delta)^s H_U(x, \cdot)=0\quad \text{ in }U,\qquad 
    H_U(x,y)=\Phi(x-y)\quad \text{ for }y\in \rn\backslash U,
\end{align*}
where $\Phi$ is the fundamental solution of the fractional Laplacian in $\rn$.  The Robin function of $U$ is given by  $\mathcal R_{U}(x):=H_U(x,x)$ for $x\in U.$

\begin{lemma} \label{robin:lem}
Let $s\in(0,1)$, $2s<N$, $V\subset \rn\backslash \overline{B_2}$ be a smooth open bounded set, and let
\begin{align}\label{omegadef}
\Omega_r:= B_r\cup V.    
\end{align}
Then $\mathcal R_{\Omega_r}$ has at least two stable critical points (which are local minima) $\xi_r\in B_r(0)$ and $\widetilde \xi_{r}\in V$ and there is $r_0\in(0,1)$ such that 
\begin{align}\label{ecgoal}
\mathcal R_{\Omega_r}(\xi_r)>\mathcal R_{\Omega_r}(\widetilde \xi_r)\qquad \text{for every $r\in (0,r_0)$}.
\end{align}
\end{lemma}
\begin{proof}
Let $r\in(0,1)$ and let $\Omega$ be as in \eqref{omegadef}. Since $\mathcal R_{\Omega_r}(x)\to \infty$ as $x\to \partial \Omega_r$, we have that $\mathcal R_{\Omega_r}$ has at least two local minima $\xi_r\in B_r(0)$ and $\widetilde \xi_{r}\in V$.     By scalings and \eqref{R0},
\begin{align}\label{ec1}
\mathcal R_{B_r}(x)=r^{2s-N}\mathcal R_{B_1}(x/r)
\geq r^{2s-N}\mathcal R_{B_1}(0)
=r^{2s-N}C_0\qquad  \text{ for }x\in B_r,\qquad C_0:=\frac{2^{1-2 s} \Gamma \left(\frac{N}{2}\right)}{(N-2 s) \Gamma (s)^2\pi ^{\frac{N}{2}}}.
\end{align}

For $x\in B_r(0)$ and $y\in \R^N$ let $W(x,y):=H_{B_r}(x,y)-H_{\Omega_r}(x,y)$. Then, 
\begin{align*}
    (-\Delta)^s W(x,\cdot)=0\quad \text{ in }B_r(0),\qquad 
    W(x,y)=\chi_{V}(y)\mathcal G_{\Omega_r}(x,y)\quad \text{ for }y\in \rn\backslash B_r(0). 
\end{align*}
We can represent $W$ with the fractional Poisson kernel $\Gamma_{B_r}$ for the ball $B_r$ (see, for instance, \cite{N15,B99,AJS18}), namely,
\begin{align*}
W(x,y)
=\int_{\rn\backslash B_r}\Gamma_{B_r}(y,z)\chi_{V}(z)\mathcal G_{\Omega_r}(x,z)\, dz
=\widetilde\gamma_{N,s}\int_{V}\frac{(r^2-|y|^2)^s}{|y-z|^N(|z|^2-r^2)^s}\mathcal G_{\Omega_r}(x,z)\, dz
\end{align*}
for $x,y\in B_r$ with $\widetilde\gamma_{N,s}=\frac{2}{\Gamma(s)\Gamma(1-s)\sigma_N}$ and $\sigma_N=|\partial B_1|=\frac{2\pi^\frac{N}{2}}{\Gamma(\frac{N}{2})}$. Observe that 
\begin{align*}
\max_{x\in B_r,\, z\in V}  \mathcal G_{\Omega_r}(x,z)
=\max_{x\in B_r,\, z\in V} (\Phi(x-z)-H_{\Omega_r}(x,z))
\leq \max_{x\in B_1,\, z\in V} \Phi(x-z)=:C_1 
\end{align*}
Hence,
\begin{align*}
|W(x,y)|\leq C_1\widetilde\gamma_{N,s}\int_{\rn\backslash B_r}\frac{(r^2-|y|^2)^s}{|y-z|^N(|z|^2-r^2)^s}\, dz=C_1
\qquad \text{ for all }x,y\in B_r,
\end{align*}
where we used that $\int_{\rn\backslash B_r}\Gamma_{B_r}(x,z)\, dz=1$ for all $x\in B_r.$  By \eqref{ec1}, for $x\in B_r(0)$,
\begin{align*}
\mathcal R_{{\Omega_r}}(x)
=H_{\Omega_r}(x,x)
\geq H_{B_r}(x,x)-|W(x,x)|
\geq \min_{B_r}\mathcal R_{B_r}-C_1\to \infty\qquad \text{ as }r\to  0.
\end{align*}

Next, we show that $\mathcal R_{\Omega_r}(\widetilde \xi_r)$ remains bounded as $r\to 0.$  Set $\Psi(x,y):=H_{\Omega_r}(x,y)-H_V(x,y)$ for $x\in V$ and $y\in \R^N$. Then,
\[
(-\Delta)^s \Psi(x,\cdot)=0 \ \text{ in } V, \qquad 
\Psi(x,y)=-\chi_{B_r}(y)\,\mathcal G_{\Omega_r}(x,y)\ \text{ in } \mathbb R^N\setminus V.
\]
As before, using the fractional Poisson kernel $\Gamma_V(x,y)=-(-\Delta)^s{\mathcal G}_V(x,y)$ (see \cite[Proposition 1.2.2]{N15}), we have that $\Psi(x,x)=-\int_{B_r}\Gamma_V(x,z)\,\mathcal G_{\Omega_r}(x,z)\,dz$ for $x\in V$ and
\[
0 \le \mathcal G_{\Omega_r}(x,z)\le \max_{x'\in V,\, z'\in B_1}\Phi(x'-z')=C_1 \qquad \text{for all } x\in V,\ z\in B_r.
\]
Hence, $|\Psi(x,x)| \le C \int_{B_r}\Gamma_V(x,z)\,dz \to 0$ as $r\to 0$ uniformly for $x$ in compact subsets of $V$.  Therefore, $\mathcal R_{\Omega_r}(x)=\mathcal R_V(x)+o(1)$ as $r\to 0$ for $x\in V$. Since $\mathcal R_V(x)\to\infty$ as $x\to \partial V$, the minimum points $\widetilde \xi_r$ of $\mathcal R_{\Omega_r}$ in $V$ remain in a compact subset of $V$, and consequently
\begin{align}\label{ec2}
\mathcal R_{\Omega_r}(\widetilde \xi_r)
=\min_{V}\mathcal R_{\Omega_r}(x)=\min_{V}\mathcal R_V+o(1)\qquad \text{ as }r\to 0.
\end{align}
Finally, \eqref{ec1} and \eqref{ec2} yield \eqref{ecgoal} for some $r_0\in (0,1)$. 
\end{proof}

\begin{theorem}
Let $s\in(0,1)$, $2s<N$, $V\subset \rn\backslash \overline{B_2}$ be a smooth open bounded set, and let
\begin{align}\label{omegadef2}
\Omega_r:= B_r\cup V.    
\end{align}
Then, for $r$ sufficiently small, there is $\eps_0>0$ such that \eqref{equn}-\eqref{defp_eps} admits a family of low energy solutions which are not least energy solutions for $\eps\in(0,\eps_0)$.
\end{theorem}
\begin{proof}
Let $\Omega_r$ be as in \eqref{omegadef2}, $\eps>0$ small, and let $J_{\eps,r}:\mathcal H^s_0(\Omega_r)\to \r$ be the energy functional associated to \eqref{equn}-\eqref{defp_eps}, namely,
\begin{align*}
    J_{\eps,r}(u):=\frac{1}{2}\|u\|_{\mathcal D^s}^2 - \frac{1}{p_\eps+1}\|u\|_{p_\eps+1}^{p_\eps+1}.
\end{align*}

Let $\xi_r$ and $\widetilde \xi_r$ be as in Lemma~\ref{robin:lem} and let $\Lambda_r>0$ and $\widetilde \Lambda_r>0$ be such that $(\xi_r,\Lambda_r)$ and $(\widetilde \xi_r,\widetilde \Lambda_r)$ are local minima of the function 
\begin{align*}
\Psi_r(\xi,t)
:=
\frac{1}{2}\,H_{\Omega_r}(\xi,\xi)\,t^2
-
\log t,
\qquad
\xi \in \Omega_r,\ t > 0.
\end{align*}
 
Now we use the Lyapunov-Schmidt reduction method \cite[Theorem 1.2]{DRS17} to construct two families of one-bubble solutions  of  \eqref{equn}-\eqref{defp_eps}. 
Indeed, by \cite[Proposition 7.3]{DRS17}
we may find 2 critical points of the reduced functional  of the form $(\xi_{r,\eps},\Lambda_{r,\eps})=(\xi_r,\Lambda_r)+o(1)$ and $(\widetilde \xi_{r,\eps},\widetilde \Lambda_{r,\eps})=(\widetilde \xi_r,\widetilde \Lambda_r)+o(1)$, as $\eps\rightarrow 0$. Hence by \cite[Lemma 7.2]{DRS17} we obtain two families  $u_{\eps,r}$ and $\widetilde u_{\eps,r}$ of solutions of  \eqref{equn}-\eqref{defp_eps} each given by a suitable rescaled bubble $U$ centered at $\xi_{r,\eps}$ and $\widetilde \xi_{r,\eps}$ respectively, plus an error term that goes to zero as $\eps\rightarrow 0$, see \cite[eq. (1.9)]{DRS17}.
\\
Furthermore, following the energy expansions in \cite{DRS17}\footnote{We note that the expansion in
\cite[Section~7.1 and Proposition~7.3]{DRS17}
is for the rescaled functional and returning to the original variables gives
$J_{\varepsilon,r}(u)
=\left(1+\frac{N-2s}{4s}\varepsilon|\log\varepsilon|
+o(\varepsilon)\right)J_{\varepsilon,r}^{\mathrm{res}}(v).
$}, it is not difficult to make explicit the constant appearing in the first order term, and so obtain that \begin{align}\label{exps}
J_{\eps,r}(u_{\eps,r})
&= \frac{s}{N}\left(S_{N,s}\right)^{\frac{N}{2s}}
+ \bigl[\hat\gamma_{N,s} + \omega_{N,s}\,\Psi_r(\xi_r,\Lambda_r)\bigr]\eps
+ \frac{N-2s}{4N}\left(S_{N,s}\right)^{\frac{N}{2s}}
  \eps|\log\eps|
+o(\eps),\\\label{exps2}
J_{\eps,r}(\widetilde u_{\eps,r})
&= \frac{s}{N}\left(S_{N,s}\right)^{\frac{N}{2s}}
+ \bigl[\hat\gamma_{N,s} + \omega_{N,s}\,\Psi_r(\widetilde \xi_r,\widetilde \Lambda_r)\bigr]\eps
+ \frac{N-2s}{4N}\left(S_{N,s}\right)^{\frac{N}{2s}}
  \eps|\log\eps|
+ o(\eps),
\end{align}
as $\eps\to 0$, where $S_{N,s}$ is the best fractional Sobolev constant (see Section~\ref{section:SobolevConstant_and_bubbles}), and  $\hat\gamma_{N,s},$ and  $\omega_{N,s}>0$ are constants independent of $\eps$ and $r$.  Observe that    $\|u\|_{\mathcal D^s}^2=\|u\|^{p_\eps+1}_{p_\eps+1}$ and 
$J_{\eps,r}(u)=\left(\frac{1}{2}-\frac{1}{p_\eps+1}\right)\|u\|_{\mathcal D^s}^2= \left(\frac{s}{N}+o(1)\right)\|u\|_{\mathcal D^s}^2$, for solutions $u$ of \eqref{equn}-\eqref{defp_eps}.
Hence by \eqref{exps}-\eqref{exps2} in particular it follows that 
\[\lim_{\eps\rightarrow 0}\| u_{\eps,r}\|_{\mathcal D^s}^2=\lim_{\eps\rightarrow 0}\|\widetilde u_{\eps,r}\|_{\mathcal D^s}^2=\left(S_{N,s}\right)^{\frac{N}{2s}},
\]
and so that
\[
\lim_{\eps\rightarrow 0}\frac{\|u_{\eps,r}\|_{\mathcal D^s}^2}{\|u_{\eps,r}\|^{2}_{p_\eps+1}}= 
\lim_{\eps\rightarrow 0}\frac{\|\widetilde u_{\eps,r}\|_{\mathcal D^s}^2}{\|\widetilde u_{\eps,r}\|^{2}_{p_\eps+1}}=
\frac{\left(S_{N,s}\right)^{\frac{N}{2s}}}{\left(S_{N,s}\right)^{\frac{N}{2s}\frac{2}{2^*_s}}}=S_{N,s},
  \]
namely   $u_{\eps,r}$ and $\widetilde u_{\eps,r}$ are both  families of \emph{low energy solutions}.
Furthermore, by Lemma~\ref{robin:lem}, there is $r_0>0$ such that $H_{\Omega_r}(\xi_r,\xi_r)>H_{\Omega_r}(\widetilde \xi_r,\widetilde \xi_r)$ for $r\in (0,r_0)$. Therefore, $\Psi_r(\xi_r,\Lambda_r)>\Psi_r(\widetilde \xi_r,\widetilde \Lambda_r)$.  Fix $r\in (0,r_0),$ then, by \eqref{exps}-\eqref{exps2}, there is $\eps_0>0$ such that $J_{\eps,r}(u_{\eps,r})>J_{\eps,r}(\widetilde u_{\eps,r})$ for $\eps\in (0,\eps_0)$. Using again that $\|u_{\eps,r}\|_{\mathcal D^s}^2=\|u_{\eps,r}\|_{p_\eps+1}^{p_\eps+1}$ and $\|\widetilde u_{\eps,r}\|_{\mathcal D^s}^2=\|\widetilde u_{\eps,r}\|_{p_\eps+1}^{p_\eps+1},$ we have that 
\begin{align*}
\frac{\|u_{\eps,r}\|_{\mathcal D^s}^2}{\|u_{\eps,r}\|_{p_\eps+1}^{2}}
=\|u_{\eps,r}\|_{\mathcal D^s}^{2(1-\frac{2}{p_\eps+1})}
=\left(\left(\frac{1}{2}-\frac{1}{p_\eps+1}\right)^{-1}J_{\eps,r}(u_{\eps,r})\right)^{\frac{p_\eps-1}{p_\eps+1}}
>\left(\left(\frac{1}{2}-\frac{1}{p_\eps+1}\right)^{-1}J_{\eps,r}(\widetilde u_{\eps,r})\right)^{\frac{p_\eps-1}{p_\eps+1}}
=\frac{\|\widetilde u_{\eps,r}\|_{\mathcal D^s}^2}{\|\widetilde u_{\eps,r}\|_{p_\eps+1}^{2}}.    
\end{align*}
Hence, for $\eps\in (0,\eps_0)$, $(u_{\eps,r})$ is a family of low energy solutions that does not coincide with a family of least energy solutions. 
\end{proof}

\begin{remark}
Let us notice that as a byproduct of the asymptotic analysis we can also derive multiplicity of least energy solutions of \eqref{equn}-\eqref{defp_eps} in disconnected symmetric smooth open sets. As an example, let us consider $P\in\R^N$, with $|P|>1$, and let
\(
\Omega=B_1(-P)\cup B_1(P),
\)
and let $u_\eps$ be a family of least energy solutions of \eqref{equn}-\eqref{defp_eps}. By virtue of Remark~\ref{rem:least_are_low} and applying  Theorem~\ref{fractional:Han} there exist a subsequence $\eps_j\to0$, as $j\to+\infty$, and a point $x_0$ belonging to one of the two disjoint unit balls, such that $\|u_{\eps_j}\|_\infty=u_{\eps_j}(x_{\eps_j})\to+\infty$, with $x_{\eps_j}\to x_0$. Suppose without loss of generality that $x_0\in B_1(-P)$, then by \eqref{gamma0_Ga} and \eqref{Green1} we can deduce that $\|u_{\eps_j}\|_{L^\infty(B_1(P))}\to0$ as $j\to+\infty$. These considerations show that $v_{\eps_j}(x):=u_{\eps_j}(-x)$ is still a least energy solution of \eqref{equn} with $p=p_{\eps_j}$, but $v_{\eps_j}\not=u_{\eps_j}$ for any $j$ sufficiently large.
\end{remark}

\subsection{Examples in the local case}
Let $N\geq 3$, $\Omega$ be a smooth open bounded set, let $\tilde p_\eps=2^*-1-\eps=\frac{N+2}{N-2}-\eps$, $\eps>0$, and let us consider the local Lane-Emden problem
\begin{equation}\label{local_eq}
-\Delta u =u^{p}\quad \text{ in }\Omega,
\qquad u>0\quad \text{ in }\Omega,
\qquad u= 0\ \mbox{ on }\mathbb \partial\Omega.
\end{equation}
We refer to solutions of \eqref{local_eq} with $p=\tilde p_\eps$, which are minimizers for the following quotient
\[
S_{N,1,\eps}:=\inf_{w\in \mathcal H_0^1(\Omega)\setminus\{0\}}\frac{\|\nabla w\|^2_{L^2(\Omega)}}{\| w\|^2_{L^{2^*-\eps}(\Omega)}}
\]
as to \emph{least energy solutions} of \eqref{local_eq} with $p=\tilde p_\eps$. It is very well known that
$S_{N,1,\eps}\to S_N$ as $\eps\to0$, where 
\[
S_{N}:=\inf_{w\in \mathcal H_0^1(\Omega)\setminus\{0\}}\frac{\|\nabla w\|^2_{L^2(\Omega)}}{\| w\|^2_{L^{2^*}(\Omega)}}.
\]
So in particular any family $w_\eps$ of least energy solutions to \eqref{local_eq} with $p=\tilde p_\eps$ verifies
\begin{equation}\label{local_limit_lowenergy}
\frac{\|\nabla w_\eps\|^2_{L^2(\Omega)}}{\|w_\eps\|^2_{L^{2^*-\eps}(\Omega)}}\to S_N\qquad \text{as $\eps\to0$.}
\end{equation}
We say that a family $u_\eps$ of solutions to \eqref{local_eq} with $p=\tilde p_\eps$ satisfying \eqref{local_limit_lowenergy} is a \emph{family of low energy solutions}.

We will give examples of contractible sets (of dumbbell type) where \eqref{local_eq} admits families of low energy solutions which are not least energy solutions.

Let $\Omega_1$, $\Omega_2$ be two smooth bounded domains such that for any $i=1,2$  
\begin{equation}\label{Omega_i}
\Omega_i\subset\{(x_1,x')\in\R\times \R^{N-1}\,:\,a_i\leq x_1\leq b_i\},\qquad \Omega_i\cap\{x'=0\}\neq\emptyset,
\end{equation}
with $a_i<b_i$ and $b_1<a_2$.  For any $\delta>0$ let 
\begin{equation}\label{C_delta}
C_\delta=\{(x_1,x')\in\R\times\R^{N-1}\,:\, x_1\in(a_1,b_2),\,|x'|\leq \delta\}.
\end{equation}

For a domain $D\subset\R^N$ we denote by ${R}_D$ the Robin function of the Laplacian with Dirichlet boundary conditions in $D$.

\begin{theorem}\label{thm:local_problem}
Let $\Omega_1$, $\Omega_2\subset\R^N$ be two smooth contractible domains satisfying \eqref{Omega_i}. Let us assume that the Robin function $R_{\Omega_i}$ admits a strict minimum point $P_i\in\Omega_i$, for any $i=1,2$, and that $R_{\Omega_1}(P_1)<R_{\Omega_2}(P_2)$. Then there exist $\delta>0$ and a smooth contractible domain $\Omega_\delta$ verifying $\Omega_1\cup \Omega_2\subset\Omega_\delta\subset\Omega_1\cup\Omega_2\cup C_\delta$ such that problem \eqref{local_eq} in $\Omega=\Omega_\delta$ admits a family of low energy solutions which are not least energy ones.
\end{theorem}

\begin{remark}\label{rem:local}
Notice that it is easy to find a pair of domains $\Omega_1$, $\Omega_2$ fulfilling the assumptions of Theorem~\ref{thm:local_problem}. Indeed, for any smooth bounded domain $D\subset\R^N$ and for any $\alpha>0$ 
\[
R_{\alpha D}(x)=\alpha^{2-N} R_D(\frac{x}{\alpha}) \qquad \text{ for all }\:x\in \alpha D, 
\]
thus $\min_{\alpha D} R_{\alpha D}=\alpha^{2-N}\min_D R_{D}$. Moreover, by \cite{michpist}, for a generic domain the corresponding Robin function has only non degenerate critical points, in particular global minimum points are strict minimum points.  In particular, $\Omega_1=B_1(0)$ and $\Omega_2=B_\alpha(Q)$, with $Q=(q,0)\in\r\times\r^{N-1}$, $q\geq2$, and $\alpha\in(0,1)$, are an eligible pair.
\end{remark}

\begin{proof}
Let $\Omega_0:=\Omega_1\cup\Omega_2$, then being $\overline{\Omega_1}\cap\overline{\Omega_2}=\emptyset$,
\[
R_{\Omega_0}(x)=
\begin{cases}
R_{\Omega_1}(x) &\text{if $x\in\Omega_1$},\\
R_{\Omega_2}(x) &\text{if $x\in\Omega_2$},
\end{cases}
\]
so $P_i$'s are strict minimum points of $R_{\Omega_0}$ and, by assumption, $R_{\Omega_0}(P_1)<R_{\Omega_0}(P_2)$. 
By \cite[Lemma 3.2]{mussopistoia}, if $\Omega_\delta$ is a smooth contractible domain such that $\Omega_0\subset \Omega_\delta\subset \Omega_0\cup C_\delta$ and $\delta$ is small enough, then $R_{\Omega_\delta}$ has at least one strict minimum point $P_{1,\delta}$, close to $P_1$, and a strict minimum point $P_{2,\delta}$, close to $P_2$, such that $R_{\Omega_\delta}(P_{1,\delta})<R_{\Omega_\delta}(P_{2,\delta})$. As a consequence, by \cite[Theorem 4.1]{mussopistoia} there exist two families of solutions $u_{i,\eps}$, $i=1,2$, to \eqref{local_eq} with $p=\tilde p_\eps$ and $\Omega=\Omega_\delta$, concentrating at $P_{i,\delta}$ respectively, which are low energy solutions (see also \cite{BLR} and \cite{Rey}). More precisely, denoting by $\tilde J_{\eps}:\mathcal H^1_0(\Omega_\delta)\to \r$ the energy functional associated to \eqref{local_eq} with $p=\tilde p_\eps$, namely
\begin{align*}
    \tilde J_{\eps}(u):=\frac{1}{2}\|\nabla u\|^2_{L^2(\Omega_\delta)} - \frac{1}{\tilde p_\eps+1}\|u\|_{L^{\tilde p_\eps+1}(\Omega_\delta)}^{\tilde p_\eps+1},
\end{align*}
then it holds
\[
\tilde J(u_{\eps,i})=\frac{S^{\frac N2}}{N}+\gamma_1\eps|\log\eps|+\gamma_2\eps+2\gamma_1(\frac12+\frac12\log(R(P_{i,\delta})))\eps+o(\eps),
\]
with $\gamma_1>0$, $\gamma_2\in\R$ depending only on $N$. Therefore,
\begin{equation}\label{last11}
\tilde J(u_{\eps,2})>\tilde J(u_{\eps,1})\qquad\text{for $\eps$ sufficiently small.}
\end{equation}
Then, being $u_{\eps,i}$, $i=1,2$, solutions to \eqref{local_eq} with $p=\tilde p_\eps$
\begin{equation}\label{last21}
\|\nabla u_{\eps,i}\|_{L^2(\Omega_\delta)}^2=\|u_{\eps,i}\|^{\tilde p_\eps+1}_{L^{\tilde p_\eps+1}(\Omega_\delta)} \quad\text{and}\quad
\tilde J_{\eps}(u_{\eps,i})=\left(\frac{1}{2}-\frac{1}{\tilde p_\eps+1}\right)\|\nabla u_{\eps,i}\|_{L^2(\Omega_\delta)}^2,
\end{equation}
thus, combining \eqref{last11} and \eqref{last21}, for $\eps$ small enough
\[
S_{N,1,\eps}\leq \frac{\|\nabla u_{\eps,1}\|_{L^2(\Omega_\delta)}^2}{\|u_{\eps,1}\|_{L^{\tilde p_\eps+1}(\Omega_\delta)}^2}= \left(\left(\frac12-\frac1{\tilde p_\eps+1}\right)^{-1}\tilde J(u_{\eps,1})\right)^{\frac{\tilde p_\eps-1}{\tilde p_\eps+1}}<\left(\left(\frac12-\frac1{\tilde p_\eps+1}\right)^{-1}\tilde J(u_{\eps,2})\right)^{\frac{\tilde p_\eps-1}{\tilde p_\eps+1}}=\frac{\|\nabla u_{\eps,2}\|_{L^2(\Omega_\delta)}^2}{\|u_{\eps,2}\|_{L^{\tilde p_\eps+1}(\Omega_\delta)}^2},
\]
and this immediately implies that $u_{\eps,2}$ is not a least energy solution for $\eps$ small.
\end{proof}

\section*{Acknowledgments}
A. Saldaña is supported by the 2021 Visiting Professor Program of La Sapienza University (Italy), by
SECIHTI grant CBF2023-2024-116 (Mexico) and by UNAM-DGAPA-PAPIIT grant IN102925 (Mexico). A. Saldaña thanks the Dipartimento di Scienze di Base e Applicate per l’Ingegneria (Sapienza Università di Roma) for the kind hospitality.
F. De Marchis is partially supported by PRIN Project 2022AKNSE4 “Variational and Analytical aspects of Geometric PDE”. I. Ianni is partially supported by the PRIN-20227HX33Z project “Pattern formation in nonlinear phenomena”,
 F. De Marchis and I. Ianni are also supported by the GNAMPA-INdAM  group for the financial support through the projects E5324001950001-2025 and E53C25002010001-2026. 
 F. De Marchis and I. Ianni thank the Institute for Mathematics of the National Autonomous University of Mexico for the kind hospitality.

\subsection*{AI Use}
All mathematical ideas, methods, and proof strategies were developed entirely by the authors, who also wrote the manuscript in full without AI assistance. Upon completion, the authors used ChatGPT-6 Astra solely for proofreading. Its suggestions were reviewed critically, and only those deemed appropriate were incorporated.

\bibliographystyle{plain}
\bibliography{references}

\end{document}